\documentclass[a4paper]{article}

\usepackage{amsthm,amsmath,amssymb}
\usepackage{mathrsfs}
\usepackage{graphicx,graphics}
\usepackage{verbatim}
\usepackage{threeparttable,booktabs}
\usepackage{algorithm, algorithmicx}

\usepackage{multirow}
\usepackage{makecell}
\usepackage[noend]{algpseudocode}
\usepackage[nohead,margin=1.0in]{geometry}
\usepackage{color}
\usepackage{hyperref}\hypersetup{hypertex=true,colorlinks=true,linkcolor=blue,anchorcolor=blue,citecolor=blue}
\usepackage{epstopdf}
\graphicspath{{./figures/}}

\algdef{SE}[DOWHILE]{Do}{doWhile}{\algorithmicdo}[1]{\algorithmicwhile\ #1}
\newcommand\outerp[2]{\left\langle #1, #2 \right\rangle_{X,X^{\prime}}}
\newcommand\xnorm[1]{\left\| #1 \right\|_{X}}
\newcommand\xpnorm[1]{\left\| #1 \right\|_{X^{\prime}}}

\newtheorem{Theorem}{Theorem}[section]
\newtheorem{Assumption}{Assumption}[section]
\newtheorem{Lemma}{Lemma}[section]
\newtheorem{Proposition}{Proposition}[section]
\newtheorem{Example}{Example}[section]
\newtheorem{Corollary}{Corollary}[section]
\newtheorem{Remark}{Remark}[section]
\newtheorem{Definition}{Definition}[section]
\numberwithin{equation}{section}
\allowdisplaybreaks

\title{An Iterative Deep Ritz Method for Monotone Elliptic Problems\thanks{The work of B. Jin is supported by Hong Kong RGC General Research Fund (14306423 and 14306824), a Direct Grant for Research (2024/25) and a start-up fund, both from The Chinese University of Hong Kong.}}
\author{Tianhao Hu\thanks{Department of Mathematics, The Chinese University of Hong Kong, Shatin, N.T., Hong Kong (email: \texttt{thhu@link.cuhk.edu.hk, b.jin@cuhk.edu.hk, fengruwang@cuhk.edu.hk})} \and Bangti Jin\footnotemark[2] \and Fengru Wang\footnotemark[2]}
\date{}

\begin{document}
\maketitle
\begin{abstract}
In this work, we present a novel iterative deep Ritz method (IDRM) for solving a general class of elliptic problems. It is inspired by the iterative procedure for minimizing the loss during the training of the neural network, but at each step encodes the geometry of the underlying function space and incorporates a convex penalty to enhance the performance of the algorithm. The algorithm is applicable to elliptic problems involving a monotone operator (not necessarily of variational form) and does not impose any stringent regularity assumption on the solution. It improves several existing neural PDE solvers, e.g., physics informed neural network and deep Ritz method, in terms of the accuracy for the concerned class of elliptic problems. Further, we establish a convergence rate for the method using tools from geometry of Banach spaces and theory of monotone operators, and also analyze the learning error. To illustrate the effectiveness of the method, we present several challenging examples, including a comparative study with existing techniques.\\
\textbf{Key words}: deep Ritz method, deep neural network, iterative method, monotone elliptic problem, convergence
\end{abstract}

\section{Introduction}
General elliptic partial differential equations (PDEs) arise in many scientific areas, e.g., static advection-diffusion equation \cite{masud2004multiscale, Garth2011}, $p$-Laplace equation \cite{Liu2023,garcia1987existence}, and Stokes equation \cite{HUGHES198785}, and represent a very important class of mathematical models in engineering, physics and biology. Their efficient and accurate numerical solution is of enormous practical importance.
These problems have been extensively studied in the context of classical numerical PDE solvers, e.g., finite difference method and finite element method, and the literature on their efficient implementation and error analysis is very extensive (see, e.g., \cite{feistauer1986finite,doi:10.1137/120887655,chow1989finite}). However, these traditional methods can be computationally intensive and are often challenging to apply in complex scenarios, particularly in high-dimensional spaces.

The advent of deep learning has led to the development of a completely new paradigm for solving PDEs, including physics informed neural network (PINN) \cite{RAISSI2019686} and deep Ritz method (DRM) \cite{yu2018deep}, deep Galerkin method (DGM) \cite{sirignano2018dgm}, weak adversarial network (WAN) \cite{Zang2020}, deep least-squares method \cite{CaiChenLiu:2021} and inf-sup neural network \cite{HuoLiu:2024}.
These neural PDE solvers harness the excellent expressivity of deep neural networks (DNNs) to approximate (high-dimensional) PDE solutions and the technique of automatic differentiation \cite{Garth2011, ad-article} to efficiently evaluate the derivatives with respect to the DNN input (spatial and temporal variables) and DNN parameters, and have shown very promising empirical results in many applications, especially in the high-dimensional setting and the model-data hybrid setting \cite{Karniadakis2021,EHanJentzen2022,TanyuMaass:2023,BeckJentzen:2023}. All these neural solvers employ DNNs as the ansatz space (and possibly test space), and construct a suitable loss function using a proper formulation of the governing PDE and boundary conditions. One class of losses is based on the strong form of PDEs (PINN, DGM), and the other class of losses is based on the Ritz or weak form of PDE (DRM, WAN).
Then one minimizes the loss function with respect to the DNN parameters using gradient type methods to find an optimal parameter configuration. The learning process of all these methods progressively approximates the solution of a PDE by adjusting the DNN parameters.

Despite the impressive progresses over last few years, these methods still face outstanding challenges when dealing with general elliptic PDEs.
First, the methods based on the strong form of PDEs (e.g., PINN and DGM) require evaluating high-order derivatives when formulating the loss, as well as high regularity assumptions on the underlying exact solutions. These methods likely fail when the exact solution has only low regularity (and the population loss may not even be well defined), e.g., point sources and geometric singularities  \cite{Krishnapriyan:2021,HuJinZhou:2022,hu2024solving}, and moreover, the computational cost of evaluating high-order derivatives grows rapidly with the spatial dimensionality \cite{lu2021deepxde}.
Meanwhile, for methods based on the weak form, DRM can only be applied to problems with a variational potential, which seriously restricts the scope of its applications to self-adjoint problems.
WAN is more flexible than DRM in that it is based on the weak formulation and employs two DNNs to approximate the test (adversarial) function and trial (primal) function separately. In particular, it imposes less stringent regularity condition on the solution $u$, and can outperform PINN when dealing with problems involving weakly singular solutions.
However, its training is often much more challenging due to the saddle point structure of the loss formulation and it tends to be unstable and may exhibit strong oscillations in practice \cite{bertoluzza2024best}. We refer interested readers to Section \ref{ssec:challenge} for further discussions on these challenges.

In view of these limitations of existing new neural PDE solvers, there is an imperative need to develop neural approaches for solving high-dimensional general elliptic PDEs that can deal with weak solutions without requiring the Ritz potential yet enjoy robust convergence behavior. Generally, the variational potential may not exist or be highly nontrivial to identify. Nonetheless, we can construct a surrogate loss based at the current approximate minimizer, and then obtain a better approximation by minimizing the surrogate loss. In this work, we develop a novel
neural solver, termed as iterative deep Ritz method (IDRM), for a broad spectrum of general elliptic
problems. It fully exploits the minimization structure of the training procedure to design an iterative algorithm
that respects the geometry of the underlying function space, and maintains the strong convexity of the loss at each step by incorporating a suitable penalty term, thereby ensuring the well-definedness of the local update. By carefully designing
the step size schedule, we can ensure that each step is performed within a neighborhood of the current approximation, and the global convergence of the approach. By its very construction, IDRM is applicable to any second-order elliptic problems involving a monotone operator with a weak
regularity assumption, which may be beyond the reach of standard PINN and DRM.

Furthermore, we analyze the convergence of the method in the following two aspects: the convergence of the inexact iteration (with the tolerance given by the learning error), and the learning error at each iteration, in terms of the DNN architecture (depth, width and parameter bound) and the number of sampling points in the domain. The former relies heavily on geometric properties of the underlying Banach space, e.g., convexity, duality mapping and Bregman distance, and the latter employs the approximation theory of DNNs and Rademacher complexity and continuity property of DNN functions for the statistical error. Together, these results provide solid theoretical underpinnings of IDRM for solving a broad class of general elliptic problems, including nonlinear ones. We conduct several high- and low-dimensional numerical experiments, including stationary convection diffusion equation, $p$-Laplace equation, diffusion equation and Navier-Stokes system, to illustrate the performance of the algorithm. These numerical results show the convergence behavior of IDRM and demonstrate its practical applicability and robustness in comparison with existing neural solvers, e.g., PINN, WAN and DRM. The results consistently show that the proposed IDRM can outperform existing approaches in terms of accuracy, and thus holds big potentials for solving general monotone elliptic problems.

The rest of the paper is organized as follows. In Section \ref{sec:pre}, we describe the specific mathematical model and standing assumptions, as well as the challenges associated with existing neural PDE solvers. In Section \ref{sec:alg}, we develop a novel neural algorithm for solving general elliptic problems. In Section \ref{sec:theorem}, we discuss the convergence of the algorithm, including the model error (by an iterative algorithm) and the learning error. We present a set of numerical experiments to showcase the performance of IDRM, including a comparative study with two existing solvers (PINN and WAN) in Section \ref{sec:experiment}. All the technical proofs are given in Section \ref{sec:proof}. Throughout, we use the notation $C$, with or without a subscript, to denote a constant which depends on its argument and whose value may differ at different occurrences.

\section{Preliminaries}\label{sec:pre}
In this section, we describe the problem formulation and discuss the challenges faced by existing neural PDE solvers.
\subsection{Mathematical formulation}
Let $\Omega \subset\mathbb{R}^d$ be an open bounded domain, and $X=W^{1,p}(\Omega)$ with $1<p\leq 2$ be the standard Sobolev space on $\Omega$. We denote its dual space by $X'$, and $\langle\cdot,\cdot\rangle_{X,X'}$ denotes the duality pairing between $X$ and $X'$.
Let $\mathcal{M}(X,X^{\prime})$ be the set of all bounded operators from the space $X$ to $X^{\prime}${and $L(X,X^{\prime})\subset \mathcal{M}(X,X^{\prime})$ consists of all bounded linear operators}.
This work focuses on the following general elliptic problem:
\begin{equation}\label{maineq}
\mathcal{A}(u)=f,
\end{equation}
where the elliptic operator $\mathcal{A}\in \mathcal{M}(X,X^{\prime})$ and the source $f\in X^{\prime}$ are respectively defined in the weak form by
\begin{align*}
\outerp{v}{\mathcal{A}(u)}=&\int_{\Omega}v(x)a_{0}(x,u(x),\nabla u(x))+\nabla v(x)\cdot a_{1}(x,u(x),\nabla u(x))\,{\rm d}x,\\
\outerp{v}{f}=&\int_{\Omega}v(x)f_{0}(x)+\nabla v(x)\cdot f_{1}(x)\,{\rm d}x,
\end{align*}
where $a_m:\mathbb{R}^{d+2}\rightarrow\mathbb{R}^{d^m}$ and $f_m:\mathbb{R}\rightarrow\mathbb{R}^{d^m}$ for $m=0,1$ are given functions. (The second term in $\langle v,f\rangle_{X,X'}$ represents the weak form of a distributional source.)
The notation $B_X(0,r)$ denotes a ball in $X$ centered at $0$ of radius $r$, i.e., $B_X(0,r)=\{u\in X, \xnorm{u}\leq r\}$. The G\^ateaux derivative $D\mathcal{A}(u)\in{L}(X,X')$ of $\mathcal{A}$ is a linear operator defined by
\begin{equation*}
\lim_{\xnorm{w}\rightarrow 0^{+}}\frac{\xpnorm{\mathcal{A}(u+w)-\mathcal{A}(u)-D\mathcal{A}(u)w}}{\xnorm{w}}=0,
\end{equation*}
and $D\mathcal{A}(u)^{*}\in L(X,X')$ denotes its adjoint operator (since $X$ is a reflexive Banach space). We denote the asymmetric part of $D\mathcal{A}(u)$ by $\Sigma(u)$, i.e., $\Sigma(u)=D\mathcal{A}(u)-D\mathcal{A}(u)^{*}$.
Throughout we make the following assumptions on the mapping $\mathcal{A}$.

\begin{Assumption}The following conditions hold on the elliptic operator $\mathcal{A}$.\label{eqassp}
\begin{itemize}
\item[{\rm(i)}] $\mathcal{A}$ is continuous and G\^ateaux differentiable, and there exists a monotone function $M$ such that
\begin{equation}\label{def:cts}
\begin{aligned}
\xpnorm{\mathcal{A}(u)-\mathcal{A}(v)}&\leq M(r)\xnorm{u-v}^{p-1},\quad \forall u,v\in B_X(0,r).
\end{aligned}
\end{equation}
\item[{\rm(ii)}] $\mathcal{A}$ is coercive:
\begin{equation*}
\lim_{\xnorm{u}\rightarrow +\infty}\frac{\outerp{u}{\mathcal{A}(u)}}{\xnorm{u}}=+\infty.
\end{equation*}
\item[{\rm(iii)}] $\mathcal{A}$ is locally uniformly monotone: for some $\rho\geq 2$ and any $R>0$, there exists a positve constant $C$ such that $\outerp{u-v}{\mathcal{A}(u)-\mathcal{A}(v)}\geq C\xnorm{u-v}^{\rho}$ holds for $u,v\in B_X(0,R)$.
\item[{\rm(iv)}] The asymmetric part $\Sigma(u)$ of the operator $\mathcal{A}(u)$ is Lipschitz continuous:
\begin{equation*}
\left\|\Sigma(u)-\Sigma(v)\right\|_{L(X,X')}\leq N\xnorm{u-v},\quad \forall u,v\in B_X(0,r).
\end{equation*}
\end{itemize}
\end{Assumption}

\begin{Remark}\label{remark:assump}
Assumption \ref{eqassp}(i)--(iii) ensure that for any $f\in X^{\prime}$, there exists a unique solution $\mathcal{A}^{-1}(f)\in X$, and the inverse mapping $\mathcal{A}^{-1}$ is uniformly continuous \cite[Theorem 2.14]{2013Roubiek}: \begin{equation}\label{eqn:Acoer}
\xnorm{\mathcal{A}^{-1}(f_{1})-\mathcal{A}^{-1}(f_{2})}\leq C\xpnorm{f_{1}-f_{2}}^{\frac{1}{\rho-1}}.
\end{equation}
Assumption \ref{eqassp}(iv) is about the asymmetric part $\Sigma(u)$,
and it is satisfied if the asymmetric part is linear, e.g., linear convection term, or quadratic, e.g. Stokes equation. Thus the assumption covers a range of elliptic problems.
\end{Remark}

To illustrate the validity of Assumption \ref{eqassp}, we give three concrete examples: a linear equation in Hilbert space, a non-linear equation in Banach space ($p\neq 2$) and a non-linear equation with a non-trivial asymmetric term. Below, the notation $|\cdot|$ denotes the Euclidean norm of vectors.

The first example is about linear second-order elliptic PDEs.
\begin{Example}
This example is about a general second-order linear elliptic equation with a Neumann boundary condition:
\begin{equation}\label{eqn4}
\left\{\begin{aligned}
-\sum_{i,j=1}^d(a^{ij}(x)u_{x_i})_{x_j}+\sum_{i=1}^db^{i}(x)u_{x_i}+c(x)u&=f,&\text{in }\Omega,\\
\sum_{i,j=1}^d a^{ij}\frac{\partial u}{\partial x_i}n_j&=g,&\text{on }\partial\Omega,\\
\end{aligned}\right.
\end{equation}
where the symmetric matrix-valued coefficient $A=[a^{ij}]_{i,j=1}^d$ satisfies $\sum_{i,j=1}^da^{ij}(x)\xi_i\xi_j\geq\theta|\xi|^2$ for some $\theta>0$ and all $\xi\in \mathbb{R}^d$, $a^{ij}\in C^{1}(\Omega)$, $b^i,c\in L^\infty(\Omega)$, $c(x)\geq c_0> \frac{(\sum_{i=1}^d\|b^i\|_{L^\infty(\Omega)})^2}{2\theta}$, and $n$ denotes the unit outward normal vector to the boundary $\partial\Omega$.
Let $X=H^{1}(\Omega)$, and
the weak form of the problem follows directly from Green's identity:
\begin{equation*}
\begin{aligned}
{\left\langle v , \mathcal{A}(u) \right\rangle_{X,X^{\prime}}}&=\sum_{i,j=1}^d\int_{\Omega}a^{ij}(x)u_{x_i}v_{x_j}\,{\rm d}x+\sum_{i=1}^d\int_\Omega b^i(x)u_{x_i}v\,{\rm d}x+\int_{\Omega}c(x)uv\,{\rm d}x,\\
\outerp{v}{f}&=\int_{\Omega}f(x)v\,{\rm d}x+\int_{\partial \Omega}g(x)v\,{\rm d}x.
\end{aligned}
\end{equation*}
Trivially, the G\^ateaux derivative $D\mathcal{A}(u)\in L(X,X')$ is $D\mathcal{A}(u)v=\mathcal{A}(v)$.
By the boundedness of the coefficients, we have $$\|\mathcal{A}(u)-\mathcal{A}(v)\|_{X'}=\sup_{w\in H^{1}(\Omega), \|w\|_{H^1(\Omega)}=1}{\left\langle w, \mathcal{A}(u-v) \right\rangle_{X,X^{\prime}}}\leq C\|u-v\|_{X}.$$
Using Young's inequality, there holds
\begin{align*}
{\left\langle u , \mathcal{A}(u) \right\rangle_{X,X^{\prime}}}\geq&\int_\Omega \theta{|\nabla u|^2}-\left(\frac{\theta}{2}|\nabla u|^2+\frac{(\sum_{i=1}^d\|b^i\|_{L^\infty(\Omega)})^2}{2\theta}u^2\right)+c_0u^2\,{\rm d}x
\geq C\|u\|_{X}^2,
\end{align*}
which gives the desired coercivity. This and the linearity of $\mathcal{A}$ imply also the uniform monotonicity of $\mathcal{A}$ with $\rho=2$. The asymmetric part $\Sigma(u)$ satisfies $$\Sigma(u)-\Sigma(v)=D\mathcal{A}(u)-D\mathcal{A}(v)-D\mathcal{A}(u)^*+D\mathcal{A}(v)^*=0,$$ which gives property (iv).
\end{Example}

The next example is about the $p$-Laplace equation,  a second-order quasi-linear elliptic problem, and it has found applications in image processing, Bingham fluids and landslides.
\begin{Example}
This example considers the $p$-Laplace equation with a zero Dirichlet boundary condition:
\begin{equation*}
\left\{
\begin{aligned}
-\nabla\cdot(|\nabla u|^{p-2}\nabla u)&= f,\quad \text{in }\Omega,\\
u&= 0,\quad \text{on }\partial\Omega,
\end{aligned}
\right.
\end{equation*}
with $1<p\leq 2$. Let $X=W_0^{1,p}(\Omega)$. Then the weak form of the $p$-Laplace equation is given by:
\begin{equation*}
\outerp{v}{\mathcal{A}(u)}=\int_{\Omega}|\nabla u|^{p-2}\nabla u\cdot\nabla v\,{\rm d}x\quad\mbox{and}\quad\outerp{v}{f}=\int_{\Omega}fv\,{\rm d}x.
\end{equation*}
Recall the following fundamental inequalities \cite[Chapter 12]{Lindqvist:2019}:
\begin{subequations}
\begin{equation}\label{ineq:conti}
||b|^{p-2}b-|a|^{p-2}a|\leq 2^{2-p}|b-a|^{p-1},
\end{equation}
\begin{equation}\label{ineq:monomtone}
(|b|^{p-2}b-|a|^{p-2}a)\cdot(b-a)\geq (p-1)|b-a|^2(1+|a|^2+|b|^2)^{\frac{p-2}{2}},
\end{equation}
\end{subequations}
for $1<p\leq 2$ and $a,b\in\mathbb{R}^d$. Let $q$ be the conjugate exponent of $p$, i.e., $\frac{1}{p}+\frac{1}{q}=1$. The continuity of $\mathcal{A}$ follows from the estimate \eqref{ineq:conti}:
\begin{align*}
&\outerp{w}{\mathcal{A}(u)-\mathcal{A}(v)}=\int_{\Omega}(|\nabla u|^{p-2}\nabla u-|\nabla v|^{p-2}\nabla v)\cdot\nabla w\,{\rm d}x\\
\leq& \int_{\Omega}||\nabla u|^{p-2}\nabla u-|\nabla v|^{p-2}\nabla v||\nabla w|\,{\rm d}x\leq 2^{2-p}\int_{\Omega}|\nabla u-\nabla v|^{p-1}|\nabla w|\,{\rm d}x\\
\leq& 2^{2-p}\left(\int_{\Omega}|\nabla u-\nabla v|^{(p-1)q}\,{\rm d}x\right)^{\frac{1}{q}}\left(\int_{\Omega}|\nabla w|^{p}\,{\rm d}x\right)^{\frac{1}{p}}\leq 2^{2-p}\|u-v\|_{W^{1,p}(\Omega)}^{p-1}\|w\|_{W^{1,p}(\Omega)}.
\end{align*}
Using Poincar\'e's inequality
$\outerp{u}{\mathcal{A}(u)}=\int_\Omega|\nabla u|^p{\rm d}x\geq c\|u\|_{X}^p$, we can deduce the coercivity.
With the help of inequality \eqref{ineq:monomtone}, we get
$$\left\langle u-v , \mathcal{A}(u)-\mathcal{A}(v) \right\rangle_{X,X^{\prime}}\geq (p-1)\int_{\Omega}|\nabla u-\nabla v|^2(1+|\nabla u|^2+|\nabla v|^2)^{\frac{p-2}{2}}\,{\rm d}x.$$
When {\color{blue}$p=2$,} we can deduce the desired monotonicity with $\rho=2$ by Poincar\'e's inequality.
If $p<2$, by H\"older's inequality, we have
\begin{align*}
\int_\Omega|\nabla u-\nabla v|^p{\rm d}x=&\int_\Omega|\nabla u-\nabla v|^p(1+|\nabla u|^2+|\nabla v|^2)^{-(1-\frac{p}{2})\frac{p}{2}}(1+|\nabla u|^2+|\nabla v|^2)^{(1-\frac{p}{2})\frac{p}{2}}{\rm d}x\\
\leq &\left(\int_\Omega|\nabla u-\nabla v|^2(1+|\nabla u|^2+|\nabla v|^2)^{-(1-\frac{p}{2})}{\rm d}x\right)^{\frac{p}{2}}\left(\int_\Omega(1+|\nabla u|^2+|\nabla v|^2)^{\frac{p}{2}}{\rm d}x\right)^{\frac{2-p}{2}}.
\end{align*}
In view of the inequality $(a+b)^t\leq a^t+b^t$ when $0<t<1$ and $a,b\geq 0$, we can deduce
\begin{align*}
\int_\Omega|\nabla u-\nabla v|^p{\rm d}x\leq& \left(\frac{1}{p-1}\left\langle u-v , \mathcal{A}(u)-\mathcal{A}(v)\right\rangle_{X,X^{\prime}}\right)^\frac{p}{2}\left(\int_\Omega1+|\nabla u|^p+|\nabla v|^p{\rm d}x\right)^{\frac{2-p}{2}}\\
\leq &C(p,R,\Omega)\left\langle u-v , \mathcal{A}(u)-\mathcal{A}(v) \right\rangle_{X,X^{\prime}}^\frac{p}{2},
\end{align*}
with $\xnorm{u}\leq R$, $\xnorm{v}\leq R$. So the local monotinicity of $\mathcal{A}$ holds with $\rho=2$.
Last, since the elliptic operator $\mathcal{A}$ has a potential $\Phi(u)=\frac{1}{p}\int_{\Omega}|\nabla u|^{p}\,{\rm d}x$, i.e., $D\Phi=\mathcal{A}$, its G\^ateaux derivative is symmetric:
\begin{align*}
\outerp{w}{D\mathcal{A}(u)v}&=\lim_{t\to 0}\frac{\Phi(u+tw)+\Phi(u+tv)+\Phi(u-tw)+\Phi(u-tv)-4\Phi(u)}{t^{2}}=\outerp{v}{D\mathcal{A}(u)w}.
\end{align*}
Then we deduce that $\Sigma(u)=0$, and that Assumption \ref{eqassp}(iv) also holds.
\end{Example}

The last example is about the steady-state Navier-Stokes system.
\begin{Example}
Consider the Navier-Stokes system
\begin{equation*}
\left\{
\begin{aligned}
-\mu\Delta \boldsymbol{u} + \nabla p + (\boldsymbol{u}\cdot\nabla) \boldsymbol{u}&= \boldsymbol{f}&\text{ in }\Omega,\\
\nabla\cdot\boldsymbol{u}&=0&\text{ in }\Omega,\\
\boldsymbol{u}&=0&\text{ on }\partial \Omega,
\end{aligned}
\right.
\end{equation*}
where $\mu>0$ is the viscosity coefficient.
Let $X=(H_{{\rm div},0}^1(\Omega))^3=\{{\boldsymbol{u}}\in (H_{0}^1(\Omega))^3, \nabla\cdot\boldsymbol{u}=0\}$. Then
the weak form of the system is given by
$$\outerp{\boldsymbol{v}}{\mathcal{A}(\boldsymbol{u})}=\int_\Omega\mu\nabla\boldsymbol{u}:\nabla\boldsymbol{v}+ (\boldsymbol{u}\cdot\nabla) \boldsymbol{u}\cdot\boldsymbol{v}\,{\rm d}{x}\quad\mbox{and}\quad\outerp{\boldsymbol{v}}{\boldsymbol{f}}=\int_{\Omega}\boldsymbol{f}\cdot\boldsymbol{v}\,{\rm d}{x}$$
We denote the convection term by $n_{{\rm c}}(\boldsymbol{u},\boldsymbol{v},\boldsymbol{w})=\int_\Omega(\boldsymbol{u}\cdot\nabla) \boldsymbol{v}\cdot\boldsymbol{w}\,{\rm d}x$. By  \cite[(6.24), Lemma 6.11]{GLOWINSKI20033},  we have
\begin{equation}\label{stokesquality}
n_{{\rm c}}(\boldsymbol{u},\boldsymbol{v},\boldsymbol{w})=-n_{{\rm c}}(\boldsymbol{u},\boldsymbol{w},\boldsymbol{v})\quad \mbox{and}\quad|n_{{\rm c}}(\boldsymbol{u},\boldsymbol{v},\boldsymbol{w})|\leq C_B\|\boldsymbol{u}\|_X\|\nabla \boldsymbol{v}\|_{(L^2(\Omega))^3}\|\boldsymbol{w}\|_X,
\end{equation}
with $C_B>0$. Then, we can verify Assumption \ref{eqassp} (i)-(iii):
\begin{align*}
&\outerp{\boldsymbol{w}}{\mathcal{A}(\boldsymbol{u})-\mathcal{A}(\boldsymbol{v})}=\mu\int_\Omega\nabla\boldsymbol{w}:\nabla(\boldsymbol{u}-\boldsymbol{v})\,{\rm d}{x}+ n_{{\rm c}}(\boldsymbol{u},\boldsymbol{u},\boldsymbol{w})-n_{{\rm c}}(\boldsymbol{v},\boldsymbol{v},\boldsymbol{w})\\
=&\mu\int_\Omega\nabla\boldsymbol{w}:\nabla(\boldsymbol{u}-\boldsymbol{v})\,{\rm d}{x}+ n_{{\rm c}}(\boldsymbol{u},\boldsymbol{u}-\boldsymbol{v},\boldsymbol{w})+n_{{\rm c}}(\boldsymbol{u}-\boldsymbol{v},\boldsymbol{v},\boldsymbol{w})\leq C\xnorm{\boldsymbol{u}-\boldsymbol{v}}\xnorm{\boldsymbol{w}}.
\end{align*}
By Poincar\'e's inequality and \eqref{stokesquality}, we have
\begin{equation*}
\outerp{\boldsymbol{u}}{\mathcal{A}(\boldsymbol{u})}=\mu\int_\Omega\nabla\boldsymbol{u}:\nabla\boldsymbol{u}\,{\rm d}{x}+ n_{{\rm c}}(\boldsymbol{u},\boldsymbol{u},\boldsymbol{u})=\mu\int_\Omega\nabla\boldsymbol{u}:\nabla\boldsymbol{u}\,{\rm d}\boldsymbol{x}\geq C\mu\xnorm{\boldsymbol{u}}^2.
\end{equation*}
Next, direct computation gives
\begin{align*}
\outerp{\boldsymbol{u}-\boldsymbol{v}}{\mathcal{A}(\boldsymbol{u})-\mathcal{A}(\boldsymbol{v})}&=\mu\int_\Omega\nabla(\boldsymbol{u}-\boldsymbol{v}):\nabla(\boldsymbol{u}-\boldsymbol{v})\,{\rm d}{x}+ n_{{\rm c}}(\boldsymbol{u},\boldsymbol{u},\boldsymbol{u}-\boldsymbol{v})-n_{{\rm c}}(\boldsymbol{v},\boldsymbol{v},\boldsymbol{u}-\boldsymbol{v})\\
&=\mu\int_\Omega\nabla(\boldsymbol{u}-\boldsymbol{v}):\nabla(\boldsymbol{u}-\boldsymbol{v})\,{\rm d}{x}+n_{{\rm c}}(\boldsymbol{u},\boldsymbol{u},-\boldsymbol{v})-n_{{\rm c}}(\boldsymbol{v},\boldsymbol{v},\boldsymbol{u})\\  &=\mu\int_\Omega\nabla(\boldsymbol{u}-\boldsymbol{v}):\nabla(\boldsymbol{u}-\boldsymbol{v})\,{\rm d}{x}-n_{{\rm c}}(\boldsymbol{u}-\boldsymbol{v},\boldsymbol{u},\boldsymbol{v})\\
&=\mu\int_\Omega\nabla(\boldsymbol{u}-\boldsymbol{v}):\nabla(\boldsymbol{u}-\boldsymbol{v})\,{\rm d}{x}-n_{{\rm c}}(\boldsymbol{u}-\boldsymbol{v},\boldsymbol{u},\boldsymbol{v}-\boldsymbol{u})\\
&\geq (C(\Omega)\mu-C_BR)\xnorm{\boldsymbol{u}-\boldsymbol{v}}^2.
\end{align*}
So the desired monotonicity holds with $\rho=2$ when $\mu>C(\Omega)^{-1}C_BR$. Last, we verify Assumption \ref{eqassp}(iv). Indeed, we have
\begin{align*}
\outerp{\boldsymbol{w}}{\Sigma(\boldsymbol{u})\boldsymbol{v}}=&\outerp{\boldsymbol{w}}{D\mathcal{A}(\boldsymbol{u})\boldsymbol{v}}-\outerp{\boldsymbol{v}}{D\mathcal{A}(\boldsymbol{u})\boldsymbol{w}}\\
=&n_{{\rm c}}(\boldsymbol{v},\boldsymbol{u},\boldsymbol{w})+n_{{\rm c}}(\boldsymbol{u},\boldsymbol{v},\boldsymbol{w})-n_{{\rm c}}(\boldsymbol{w},\boldsymbol{u},\boldsymbol{v})-n_{{\rm c}}(\boldsymbol{u},\boldsymbol{w},\boldsymbol{v})\\
=&n_{{\rm c}}(\boldsymbol{v},\boldsymbol{u},\boldsymbol{w})+2n_{{\rm c}}(\boldsymbol{u},\boldsymbol{v},\boldsymbol{w})+n_{{\rm c}}(\boldsymbol{w},\boldsymbol{v},\boldsymbol{u}).
\end{align*}
Consequently, we have
\begin{align*}
&\outerp{\boldsymbol{w}}{\Sigma(\boldsymbol{u}_1)\boldsymbol{v}-\Sigma(\boldsymbol{u}_2)\boldsymbol{v}}\\
=&n_{{\rm c}}(\boldsymbol{v},\boldsymbol{u_1}-\boldsymbol{u_2},\boldsymbol{w})+2n_{{\rm c}}(\boldsymbol{u_1}-\boldsymbol{u_2},\boldsymbol{v},\boldsymbol{w})+n_{{\rm c}}(\boldsymbol{w},\boldsymbol{v},\boldsymbol{u_1}-\boldsymbol{u_2})\\
\leq& 4C_B\xnorm{\boldsymbol{u_1}-\boldsymbol{u_2}}\xnorm{\boldsymbol{v}}\xnorm{\boldsymbol{w}}.
\end{align*}
This shows the Lipschitz continuity of $\Sigma$.
\end{Example}

These examples show that Assumption \ref{eqassp} is not very restrictive and does cover a wide range of important applications. Hence, it is important to develop an algorithm that can effectively treat this class of elliptic problems.
\subsection{Challenges with existing neural PDE solvers}\label{ssec:challenge}
Although problem \eqref{maineq} is well-defined at the continuous level, existing neural PDE solvers face certain challenges.
We discuss the challenges associated with three representative neural PDE solvers, specifically, PINN, DRM and WAN, separately.
All three methods transform problem \eqref{maineq} into an equivalent optimization problem, and then approximate the solution $u$ by a DNN, i.e., DNN as the ansatz space, and then
use automatic differentiation and suitable optimization algorithms, e.g., ADAM \cite{KingmaBa:2015}, to solve the resulting parameterized optimization problem (by DNN parameters), which yields a neural network approximation of problem \eqref{maineq}.
To illustrate the point, consider the following model elliptic problem:
\begin{equation}\label{eqn2}
\left\{
\begin{aligned}
-\Delta u+\vec{b}\cdot \nabla u + cu&=f,\quad\text{in }\Omega,\\
 u&=g,\quad\text{on }\partial\Omega.
\end{aligned}
\right.
\end{equation}
Below we denote the class of DNN functions by $\mathcal{N}$; see Section \ref{subsec_pracIDRM} for
the precise definition.

First, PINN employs the $L^{2}(\Omega)$-norm of the PDE residual and boundary fitting as the loss:
\begin{equation}\label{eqn:PINN}
u_{\text{PINN}}\in \mathop{\arg\min}_{u\in\mathcal{N}} \{L_{\text{PINN}}(u)=\|-\Delta u+\vec{b}\cdot \nabla u + cu-f\|_{L^{2}(\Omega)}^2+\sigma \|u-g\|_{L^{2}(\partial \Omega)}^2\},
\end{equation}
where $\sigma>0$ is the penalty parameter balancing the residuals of the governing PDE in the domain and and boundary condition. This formulation imposes only minimal assumptions on the form of the PDE, and hence it is very flexible and applicable to a broad class of PDEs.
Nonetheless, there are two major drawbacks of the approach.
First, this formulation restricts the PDE residual $-\Delta u+{\vec{b}}\cdot \nabla u + cu$ and source term $f$ to be functions in $L^{2}(\Omega)$.
Thus, it requires high smoothness of the solution $u$ (i.e.,  $u\in H^{2}(\Omega)$ and $f\in L^{2}(\Omega)$).
Moreover, PINN also cannot naturally handle various physical boundary conditions, e.g., Robin or Neumann types. For example, for the Poisson equation with the Robin boundary condition
\begin{equation}\label{eqn1}
\left\{
\begin{aligned}
-\Delta u + cu&=0,\quad \text{in }\Omega,\\
\frac{\partial u}{\partial n} + u&=g,\quad\text{on }\partial \Omega,
\end{aligned}
\right.
\end{equation}
the problem has a Ritz variational formulation, given by
\begin{equation*}
u\in\mathop{\arg\min}_{v\in H^{1}(\Omega)}\,\left\{\frac12\int_\Omega |\nabla u|^{2}+cu^2\,{\rm d} x + \frac12\int_{\partial\Omega}u^2\,{\rm d} s -\int_{\partial\Omega}ug\,{\rm d} s\right\}.
\end{equation*}
The variational formulation encodes directly the Robin boundary condition, which physically corresponds to the surface pressures acting as the external loading on the boundary $\partial\Omega$.
However, in PINN, all the boundary conditions, including Robin or Neumann type, can only be enforced approximately by penalizing their $L^{2}(\partial\Omega)$-norm (in order to ensure the ease of implementation).
Note that the imposition of the boundary condition via the penalty method as in \eqref{eqn:PINN} inherently limits the accuracy of the approximation \cite{HuJinZhou:2022}, whereas the exact imposition is often very challenging albeit highly desirable \cite{Zeinhofer:2022exact}.

Second, DRM \cite{yu2018deep} is based on the Ritz variational formulation of the elliptic problem \eqref{maineq}.
It avoids the above two issues of PINN, i.e., it does not require the existence of strong solutions and can handle the Robin / Neumann boundary condition easily.
Instead, it requires the existence of a potential $\Phi:X\rightarrow\mathbb{R}$ whose G\^ateaux derivative coincides with the elliptic operator $\mathcal{A}$, i.e., $D\Phi(u) =\mathcal{A}(u)$.
For example, for problem \eqref{eqn1}, the potential $\Phi$ coincides with its Ritz variational functional.
Generally, for problem \eqref{maineq}, one may define a candidate potential $\Phi_0$ by
\begin{equation}\label{eqn3}
\Phi_{0}(u)=\int_{0}^{1}\outerp{u}{\mathcal{A}(tu)}\,{\rm d}t,
\end{equation}
and then DRM for \eqref{maineq} is equivalent to minimizing
\begin{equation}\label{ritzmin}
u_{\text{DRM}}\in\mathop{\arg\min}_{u\in\mathcal{N}}\,\{\Phi_{0}(u)-\outerp{u}{f}\},
\end{equation}
which is the Ritz variational formulation of problem \eqref{maineq}, provided that it does exist.
In practice, there are two major difficulties when applying the DRM.
First, the G\^ateaux derivative $D\Phi_0$ of  $\Phi_0$ is not necessarily $\mathcal{A}$. This issue arises for problem
\eqref{eqn2} with a nonzero convection term $\vec{b}$.
Indeed, upon letting $X=H_0^1(\Omega)$, then
\begin{equation*}
\begin{aligned}
{}&\Phi_{0}(u)=\int_{0}^{1}\outerp{u}{\mathcal{A}(tu)}\,{\rm d}t=\frac{1}{2}\int_\Omega |\nabla u|^{2}+u\vec{b}\cdot\nabla u+cu^2\,{\rm d} x,\\
{}&\outerp{v}{D\Phi_0(u)}=\int_\Omega \nabla u\cdot\nabla v+\frac{1}{2}\left(u\vec{b}\cdot \nabla v + v\vec{b}\cdot\nabla u\right)+cuv\,{\rm d} x.
\end{aligned}
\end{equation*}
For a G\^ateaux differentiable operator $\mathcal{A}$, its potential exists if and only if its G\^ateaux derivative $D\mathcal{A}$ is symmetric, otherwise the minimization problem \eqref{ritzmin} only gives a solution of the equation $D\Phi_{0}(u)=f$ instead of the solution of problem \eqref{maineq}.
Second, the potential $\Phi_0$ may not always be convex even if $\mathcal{A}$ is monotone \cite[Theorem 1]{rockafellar1966characterization}, whereas non-convex potentials may lead to unstable or even non-convergent training in practice.

Zang et al \cite{Zang2020} proposed WAN to deal with these issues. WAN is based on transforming the weak formulation of problem \eqref{eqn2} to the following min-max problem
\begin{equation}\label{eqn:WAN}
\min_{u\in \mathcal{N}:u|_{\partial\Omega}=g}\max_{v\in \mathcal{N}:v|_{\partial\Omega}=0}\dfrac{\int_{\Omega}(\nabla u \cdot\nabla v+(\vec{b}\cdot\nabla u)v+cuv-fv)\,{\rm d}x}{\Big(\int_\Omega |\nabla v|^2{\rm d}x\Big)^{1/2}}.
\end{equation}
It is more broadly applicable than DRM in that it does not require the existence of the Ritz potential, and hence can be applied on any problem with a (stable) weak formulation.
However, the resulting min-max problem in \eqref{eqn:WAN} is very challenging to optimize, due to the nonlinearity and non-convexity of the parameter space of the DNN set $\mathcal{N}$, and the training is often very delicate and unstable \cite{bertoluzza2024best}. The training issue is attributed to the fact that
the Lipschitz continuity constant of the test maximizers with respect to the trial
functions can become arbitrarily large when approaching the exact solution; and hence the corresponding test maximizer is highly non-unique in the limit, leading to an inherent lack of numerical stability of the method \cite{Uriarte:2023}.
To address the challenge, Uriarte et al \cite{Uriarte:2023} proposed a novel deep double Ritz method (D$^2$RM), by reformulating the residual minimization as an equivalent minimization of a Ritz functional fed by optimal test functions computed from another Ritz functional. It combines two neural networks for approximating trial functions and optimal test functions using a nested double Ritz minimization strategy so as to overcome the training instability of WAN.

Now we summarize the limitations with the existing neural PDE solvers in Table \ref{table:limitations}. These observations motivate the development of new neural PDE solvers that can overcome the drawbacks.
\begin{table}[ht]
\centering
\caption{Limitations of three existing neural PDE solvers.}\label{table:limitations}
\begin{tabular}{c|c|c|c}
\toprule
method&solution regularity&applicability&convexity of formulation\\
\midrule
PINN&\textbf{strong}&non-symmetric&convex\\
DRM&weak&\textbf{symmetric}&\textbf{non-convex}\\
WAN&weak&non-symmetric&\textbf{non-convex}\\
\bottomrule
\end{tabular}
\end{table}

\section{Iterative deep Ritz method (IDRM)}\label{sec:alg}

\subsection{The derivation of IDRM loss}

To address the challenges in Section \ref{ssec:challenge}, we propose a novel iterative deep Ritz method (IDRM).
The method requires only weak regularity, is applicable to both symmetric and non-symmetric equations, and the loss function is convex with respect to the DNN function (but not the DNN parameters).
The method extends the DRM by transforming the asymmetric problem into a sequence of symmetric problems and incorporating an iterative mechanism.
At each iteration, we construct a convex surrogate loss locally on the fly to guide the training. (Note that the convexity is with respect to $u$, and the surrogate loss is still nonconvex with respect to the DNN parameters, due to the nonlinearity of the DNN functions in the DNN parameters.)
The convergence is expected to be stable and fast due to the initialization of the current iteration using the result from last one.
This alleviates the need of globally recasting a monotone equation into a Ritz formulation, which is generally infeasible. For example, for the model problem \eqref{eqn2}, at each iteration, we may recast the problem into a symmetric one with the external source including both the function $f$ and the asymmetric part $\vec{b}\cdot \nabla u$:
\begin{equation}\label{iter-prob}
\left\{
\begin{aligned}
-\Delta u_{k+1} + cu_{k+1}&=f-\vec{b}\cdot \nabla u_{k},\quad\text{in }\Omega,\\
u_{k+1}&=g,\quad\text{on }\partial\Omega,
\end{aligned}
\right.
\end{equation}
and solve the modified problem using  DRM:
\begin{equation*}
u_{k+1}\in \mathop{\arg\min}_{u\in \mathcal{N}: u|_{\partial\Omega} =g} \frac{1}{2}\int_{\Omega}|\nabla u|^{2}+cu^{2}\,{\rm d}x-\int_{\Omega}fu-\vec{b}\cdot \nabla u_{k}u\,{\rm d}x.
\end{equation*}
Obviously this is simply a fixed point iteration for solving \eqref{eqn2}.
Under suitable conditions, we expect the convergence of the iteration to the unique solution of problem \eqref{eqn2}.

Similar to the model problem \eqref{eqn2},  the symmetric part of the general case in \eqref{maineq} can be defined by the G\^ateaux derivative of the potential $\Phi_{0}(u)$ in \eqref{eqn3}, which indicates that the asymmetric part is $\mathcal{A}(u)-D\Phi_{0}(u)$.
Hence for \eqref{maineq}, the modified equation reads
\begin{equation*}
D\Phi_{0}(u_{k+1})-D\Phi_{0}(u_{k})=f-\mathcal{A}(u_{k}).
\end{equation*}
Then DRM can be applied to the modified equation.
Namely, at step $k$, we solve the following minimization problem:
\begin{equation}\label{eqn:Ritz-modified}
u_{k+1}\in\mathop{\arg\min}_{u\in \mathcal{N}}\int_{0}^{1}\outerp{u-u_k}{D\Phi_{0}(t(u-u_k)}\,{\rm d}t-\outerp{u-u_k}{f-\mathcal{A}(u_k)}.
\end{equation}
In view of the identity $$\Phi_{0}(u)=\int_{0}^{1}\outerp{u}{\mathcal{A}(tu)}\,{\rm d}t=\int_{0}^{1}\outerp{u}{D\Phi_{0}(tu)}\,{\rm d}t,$$
we may rewrite the minimization problem \eqref{eqn:Ritz-modified} as
\begin{equation}\label{eqn:modified-gen}
u_{k+1}\in\mathop{\arg\min}_{u\in \mathcal{N}}\Phi_{0}(u-u_k)+\outerp{u-u_k}{\mathcal{A}(u_k)-f}.
\end{equation}
This functional can be viewed as an extension of the Ritz functional to non-linear and asysmetric problems. In fact, if $\mathcal{A}$ is a linear self-adjoint operator (i.e., $\mathcal{A}(u)=\mathcal{A}u$ and $\mathcal{A}=\mathcal{A}^*$), we can obtain
\begin{align*}
{}&\Phi_{0}(u-u_k)+\outerp{u-u_k}{\mathcal{A}(u_k)-f}\\
=&\tfrac{1}{2}\outerp{u}{\mathcal{A}u}-\outerp{u}{\mathcal{A}u_{k}}+\tfrac{1}{2}\outerp{u_{k}}{\mathcal{A}u_{k}} + \outerp{u-u_{k}}{\mathcal{A}u_{k}-f}\\
{}=&\tfrac{1}{2}\outerp{u}{\mathcal{A}u}-\outerp{u}{f}+C_k,
\end{align*}
where the constant term $C_k$ depends only on $u_k$ and plays no role on the minimizer. Thus, the loss essentially coincides with the Ritz variational function.
%\cite{LiuZhu:2024}.

In practice, the potential $\Phi_{0}$ may be non-convex, and the scheme \eqref{eqn:modified-gen} may not always converge with a constant step size.
Thus, we augment the optimization function with two extra parameters $\mu$ and $\lambda$ to enhance the convexity and control the step size schedule.
Specifically, we incorporate an extra penalty term $\xnorm{u}^{2}$ to enhance its convexity as well as the stability in the  neural network training:
\begin{equation*}
\Phi_{\mu}(u)=\int_{0}^{1}\outerp{u}{\mathcal{A}(tu)}\,{\rm d}t + \mu\xnorm{u}^{2},
\end{equation*}
where $\mu>0$ controls the amount of smoothing.
The potential $\Phi_\mu$ enjoys the following nice property: $\Phi_\mu$ is locally convex in the neighborhood
\begin{equation}\label{eqn:Vmu}
V_\mu=\{u\in X, \xnorm{u}\leq\tfrac{3\mu(p-1)}{2N}\},
\end{equation}
where the constant $N$ is given in Assumption \ref{eqassp} (iv); see Proposition \ref{prop:minimizer} for the precise statement. In addition, for the term $\langle u-u_k,\mathcal{A}(u_k)-f\rangle_{X,X'}$ in the loss in \eqref{eqn:modified-gen}, one may incorporate an adjustable step size $\lambda$.
Since the step size $\lambda$ may vary with the iteration $k$ (i.e., the step size at the $k$th iteration is $\lambda_k$), we propose a surrogate loss $L^k(u)$ based at the current approximation $u_k$ defined by
\begin{equation}\label{eqn:surrogate}
{L}^{k}(u) = \Phi_{\mu}\left[\lambda_k(u-u_{k})\right] + \outerp{\lambda_k(u-u_{k})}{\mathcal{A}(u_{k})-f}.
\end{equation}
The first term penalizes moving too far away from the current iterate $u_k$ (i.e., local update), while preserving the geometry of the underlying space, and the second term recovers the gradient direction $\mathcal{A}(u_k)-f$.
By the coercivity of the elliptic operator $\mathcal{A}$ (cf. Assumption \ref{eqassp} (ii)), we have $\lim_{\lambda\rightarrow\infty}\frac{\Phi_{\mu}\left[\lambda(u-u_{k})\right]}{\lambda}=\infty$, which indicates that a larger $\lambda_k$ leads to a smaller update $u-u_{k}$.

The overall construction admits the following interpretation in the lens of preconditioning: using the definition of the dual potential $\Phi_\mu^*$ (cf. Definition \ref{def:dual}), minimizing the loss $L^k$ is equivalent to minimizing $\Phi^{*}_{\mu}\left[\mathcal{A}(u_{k})-f\right]$. By the continuity and convexity of $\Phi_{\mu}$, we deduce that $\Phi_{\mu}^{*}:X'\to \mathbb{R}_+$ is a convex penalty. Hence, the formulation \eqref{eqn:surrogate} requires lower regularity on the solution $u$ when compared with directly minimizing the $L^{2}(\Omega)$-norm of the PDE residual, which is used  in neural PDE solvers based on the strong form like PINN.
Then we can formally write the minimizer of \eqref{eqn:surrogate} into:
\begin{equation*}
u_{k+1}=u_{k}-\lambda_{k}^{-1}D\Phi_{\mu}^{*}\left(\mathcal{A}(u_{k})-f\right).
\end{equation*}
By the construction of the functional $\Phi_{0}$ and the definition of the dual function $\Phi_\mu^*$, we deduce that $D\Phi_{0} \approx \mathcal{A}$ and $D\Phi_{\mu}^{*}=\left(D\Phi_{\mu}\right)^{-1}$.
Therefore, it follows that $D\Phi_{\mu}^{*}\approx \mathcal{A}^{-1}$ is a good preconditioner for the problem.
Thus, a better approximation by the IDRM can be expected when compared with neural solvers based on the strong form, e.g., PINN, especially for solutions involving singularities. Numerically, proper preconditioning can accelerate the training of the resulting loss.

\subsection{Practical aspects of the IDRM solver}\label{subsec_pracIDRM}
In the IDRM solver, we employ a fully connected feedforward DNN $u_{\theta}:\mathbb{R}^{d}\rightarrow\mathbb{R}$ to approximate the solution $u$, which is defined recursively by
$\boldsymbol{y}_{0}(\boldsymbol{x})=\boldsymbol{x}$,
$\boldsymbol{y}_{\ell}(\boldsymbol{x})=\sigma\left(A_{\ell} \boldsymbol{y}_{\ell-1}(\boldsymbol{x})+b_{\ell}\right)$, for $\ell=1,2, \ldots, L-1$ and $u_{\theta}=y_{L}(\boldsymbol{x})=A_{L} \boldsymbol{y}_{L-1}(\boldsymbol{x})+b_{L}$, where the weight matrices $A_{\ell}\in\mathbb{R}^{N_{\ell}\times N_{\ell-1}}$ and bias vectors $b_{\ell}\in\mathbb{R}^{N_{\ell}}$ are trainable, collectively stacked into one vector  $\theta$ of DNN parameters, and $\sigma:\mathbb{R}\to \mathbb{R}$ is a nonlinear activation function and applied componentwise to a vector. In this work, it is taken to be the hyperbolic tangent function $t\mapsto \frac{e^t-e^{-t}}{e^t+e^{-t}}$. The numbers $L$ and $W=\max_{\ell=1,\ldots,L}N_\ell$ are called the depth and width of the DNN $u_{\theta}$. We denote by $\mathcal{N}(L,W,B)$ the collection of DNNs $u_{\theta}$ which satisfies (i) The depth and width are less than $L$ and $W$, respectively; and (ii) All parameters are bounded by $B$. Below we abbreviate it to $\mathcal{N}$ wherever there is no confusion.

In practice, the integrals in the surrogate loss $L^k(u)$ has to be further discretized, commonly using Monte Carlo methods. Let $U(\Omega\times(0,1))$ be the uniform distribution
over the set $\Omega\times (0,1)$, and let
\begin{align*}
 I_{1}^{k}((x, t), u) &= (u-u_k)a_0(x,t\lambda_k(u-u_k), t\lambda_k\nabla (u-u_k))+\nabla(u-u_k)\cdot a_1(x,t\lambda_k(u-u_k), t\lambda_k\nabla (u-u_k)),\\
I_{2}^{k}(x, u) &= |u-u_k|^p + |\nabla u-\nabla u_k|^p,\\
I_{3}^{k}(x, u) &= (u-u_k)(a_0(x,u_k,\nabla u_k)-f_0)+\nabla(u-u_k)\cdot(a_1(x,u_k,\nabla u_k)-f_1).
\end{align*}
Then the population loss $L^{k}(u)$ can be equivalently rewritten as
\begin{align*}
&{L}^{k}(u)=\lambda_k\int_{\Omega\times(0,1)}I_{1}^{k}((x, t), u)\,{\rm d}(x,t) + \mu\lambda_k^2\left[\int_{\Omega}I_{2}^{k}(x, u)\,{\rm d}x\right]^{2/p}+ \lambda_k\int_{\Omega}I_{3}^{k}(x, u)\,{\rm d}x\\
=&\lambda_k|\Omega|\mathbb{E}_{ U(\Omega\times(0,1))} [I_{1}^{k}((X,T), u(X))] + \mu\lambda_k^2\left[|\Omega|\mathbb{E}_{ U(\Omega)}[I_{2}^{k}(X, u(X))]\right]^{2/p}+ \lambda_k|\Omega|\mathbb{E}_{U(\Omega)}[I_{3}^{k}(X, u(X))],
\end{align*}
where $\mathbb{E}_U[\cdot]$ denotes taking expectation with respect to the distribution $U$.
 Let $\{(X_i,T_i)\}_{i=1}^N$ be $N$ independent and identically
distributed (i.i.d.) random samples drawn from the uniform distribution $U(\Omega\times(0,1))$. Then the empirical loss $\widehat{L}^k(u)$ is given by
\begin{equation*}
\widehat{L}^{k}(u)=\lambda_k\frac{|\Omega|}{N}\sum_{i=1}^{N}I_{1}^{k}((X_{i}, T_{i}), u(X_i)) +\mu\lambda_k^2\left[\frac{|\Omega|}{N}\sum_{i=1}^{N}I_{2}^{k}(X_{i}, u(X_i))\right]^{2/p} + \lambda_k\frac{|\Omega|}{N}\sum_{i=1}^{N}I_{3}^{k}(X_{i}, u(X_i)).
\end{equation*}
The training of the empirical loss $\widehat{L}^k(u)$ reads
\begin{equation*}
\widehat{u}_{k+1}\in \mathop{{\rm argmin}}\limits_{\mathcal{N}\cap \widehat{u}_{k}+V_{\mu}}\widehat{L}^{k}(u).
\end{equation*}
The existence of a minimizer $\widehat{u}_{k+1}$ to the empirical loss $\widehat{L}^k$ follows from the maximum bound $B$ on the DNN parameters $\theta$, the finite dimensionality of the parameter space, and the smoothness of the activation function $\sigma$.
However, the loss $\widehat{L}^k$ is non-convex due to the nonlinearity of the DNN $u_\theta$.
Thus, generally there are multiple global minima, and it is challenging to find a global minimizer.
Nonetheless, in practice, the loss $\widehat{L}^k$ can be effectively minimized by any stand-alone optimizer, e.g., ADAM \cite{KingmaBa:2015} and L-BFGS \cite{ByrdLu:1995}.
Note that the gradient of the empirical loss $\widehat{L}^k(u)$ with respect to the input
$\boldsymbol{x}$ (i.e., spatial derivative) and to the DNN parameters $\theta$ can be computed efficiently using
automatic differentiation \cite{ad-article}. Moreover, more specialized training schemes, e.g., greedy training \cite{SiegelXu:2023} and Gauss-Newton method \cite{jnini2024gaussnewton,Zeinhofer:2023natural}, can also be applied. Note that the loss $\widehat{L}^k(u)$ changes at each outer loop, and we require solving one optimization problem per loop. In practice, we start the training afresh, using the DNN parameter from the last loop to warm-start the current optimization process, so as to speed up the convergence of the overall training.

Now we can present the details of the proposed neural PDE solver in Algorithm \ref{alg:alg1}, where the rule for selecting the step size $\lambda_k$ is given in Theorem \ref{thm:modelconv} and Remark \ref{rmk:step-size} below.
In Algorithm \ref{alg:alg1}, we indicate only the minimization of the population loss $L^k$ instead of the empirical loss $\widehat{L}^k$. This is to isolate the learning error $\epsilon_M$, which include the approximation error and quadrature error.

It is worth noting that the implementation of IDRM (along with the convergence analysis in Section \ref{sec:theorem}) is not limited to the class $\mathcal{N}$ of feedforward neural networks and Monte Carlo sampling.
Other architectures (including different activation functions) can also be employed, provided that the neural network class $\mathcal{N}$ belongs to $X$, i.e. $\mathcal{N}\subset X$ and maintains the property that the loss is differentiable with respect to the neural network parameters $\theta$. This is a very mild condition, but does exclude the standard $\mathrm{ReLU}$ neural network. Similarly, for the discretization of the integrals, tradition quadrature rules, e.g., Gaussian quadrature, can be employed if the domain $\Omega$ is not high-dimensional. Indeed, in Section \ref{sec:experiment}, we employ the trapezoidal rule on a uniformly distributed grid to solve the three-dimensional Navier-Stokes system. These methods are often direct to implement and their accuracy and stability have also been thoroughly analyzed. For high-dimensional problems, more efficient quadrature rules (than standard Monte Carlo) have been proposed \cite{BeckerJentzen:2024,RuschRus:2024}, and may also be employed for implementing IDRM. We refer readers to the interesting work \cite{Rivera:2022} for an in-depth investigation on the impact of quadrature rules on the accuracy of neural solvers for PDEs, and detailed discussions on the pros and cons of various quadrature rules.

\begin{algorithm}[hbt!]
\caption{IDRM for general elliptic problems}
\label{alg:alg1}
\begin{algorithmic}[1]
\State Initialize the networks $u_{0}:\Omega \rightarrow \mathbb{R}$, and set the step size $\lambda_{0}$ and tolerance $\epsilon_{\rm tol}$.
\Do
\State Initialize $u_{k+1}=u_{k}$.
\State Define the surrogate loss \eqref{eqn:surrogate}.
\State Find $u_{k+1}\in \mathcal{N}$ such that
\begin{equation}\label{optimization}
{L}^{k}(u_{k+1})\leq \min_{u\in u_{k}+V_{\mu}}{L}^{k}(u)+\epsilon_M.
\end{equation}
\State Compute the step size $\lambda_{k+1}$.
\State Repeat step 5 with $\lambda_{k+1}$.
\doWhile{$-L^{k}(u_{k+1})\leq \epsilon_{\rm tol}$.}
\end{algorithmic}
\end{algorithm}

\section{Convergence results for IDRM}\label{sec:theorem}

We now present convergence results for the proposed IDRM. In the analysis, we first treat the learning error $\epsilon_M$ as a constant, and have the following model convergence result as well as preasymptotic convergence rates. Below we use two exponents:
\begin{equation*}
\alpha=-\frac{\rho^2(\rho-1)-p^2(p-1)}{p\rho[p(p-1)-(\rho-1)^2]}\quad \mbox{and}\quad \beta=\frac{(p-1)(\rho^2-p(\rho-1))}{\rho(p(p-1)-(\rho-1)^2)}.
\end{equation*}
The proof requires several preliminary tools: we first prove the local convexity and continuity of the G\^ateaux derivative of the potential $\Phi_\mu$, which ensures that the surrogate loss $L^k(u)$ has one unique minimizer at each step. One key tool in the analysis is the dual potential $\Phi_\mu^*$.
\begin{Definition}\label{def:dual}
The dual potential $\Phi_\mu^*:X'\to\mathbb{R}$ of a potential $\Phi_\mu:X\to\mathbb{R}$ is defined by
\begin{equation*}
\Phi_{\mu}^{*}(f)=\max_{u\in V_\mu}\left\{\outerp{u}{f}-\Phi_{\mu}(u)\right\},
\end{equation*}
where the set $V_\mu$ is defined in \eqref{eqn:Vmu}.
\end{Definition}

Throughout we make the following assumption on the sequence $\{u_{k}\}_{k=0}^\infty$.
\begin{Assumption}\label{assump:bound}
The sequence $\{u_{k}\}_{k=0}^\infty$ is uniformly bounded in $X$: $\{u_{k}\}_{k=0}^\infty\subset B(0,R)$.
\end{Assumption}

Using the dual potential $\Phi_\mu^*$, we can establish the convergence as well as the preasymptotic decay rate of $\xnorm{u_k-u^*}$. The proof relies heavily on the geometry of the underlying Banach space, e.g., duality map and Bregman distance. The lengthy proof is given in Section \ref{ssec:modelconv}.
\begin{Theorem}\label{thm:modelconv}
Let $u^*$ be the solution to problem \eqref{maineq}, let $c_{\rm p}$, $c_{\rm r}$ and $c_{\rm s}$ be three constants that depend on $p,\rho,R$ and $\mu$ {\rm(}given in the proof{\rm)} and let the sequence $\{u_{k}\}_{k=0}^\infty$ be generated by Algorithm \ref{alg:alg1} with the following step size schedule
\begin{equation*}
\lambda_k=c_{\rm s}(-L^k(v_{k+1}))^{\alpha},
\end{equation*}
Then, under Assumptions \ref{eqassp} and \ref{assump:bound}, the following statements hold with $c_{\rm t}=c_{\rm p}^{-1}c_{\rm r}c_{\rm s}^{-1}$.
\begin{itemize}
\item [{\rm(i)}] Let $\Delta_c = \max\{(2c_{\rm t})^{\frac{\rho(\rho-1)}{p(p-1)}}\epsilon_M,(2c_{\rm t}\epsilon_M^{\frac{1}{\rho}})^{\frac{1}{\alpha+\beta}}\}$
 and $k^*=\max\{k\in\mathbb{N}:-L^j(v_{j+1})>\Delta_c, \forall 1\leq j\leq k\}$. Then the sequence $\{\Phi_\mu^*(p_k)\}_{k=1}^{k^*}$ is monotonically decreasing.
\item [{\rm(ii)}] Fix $\eta\in(0,1)$, let $\Delta_r = \max\{(2c_{\rm t}\eta^{-1})^{\frac{\rho(\rho-1)}{p(p-1)}}\epsilon_M,(2c_{\rm t}\eta^{-1}\epsilon_M^{\frac{1}{\rho}})^{\frac{1}{\alpha+\beta}}\}$ and  $k^{**}=\max\{k\in\mathbb{N}:-L^j(v_{j+1})\geq \Delta_r,\forall 0\leq j\leq k\}$. Then the following convergence rates when $0\leq k\leq k^{**}$.
\begin{itemize}
\item[{\rm(a)}] For $p=\rho=2$, there exists  $q=C(p,\rho,R,\mu)\in(0,1)$ such that
$\xnorm{u_k-u^*}\leq (\eta q)^{k}\xnorm{u_0-u^*}$.
\item[{\rm(b)}] For {$0<p(p-1)-(\rho-1)^2<1$}, there holds
$\xnorm{u_k-u^*}\leq C(p,\rho,R,\mu,\xnorm{u_0-u^*})((1-\eta)k)^{-\frac{1}{\beta-1}}.$
\end{itemize}
\end{itemize}
\end{Theorem}
\begin{Remark}\label{rmk:step-size}
Theorem \ref{thm:modelconv} chooses the step size $\lambda_k=c_{\rm s}(-L^k(v_{k+1}))^{\alpha}$. At first sight, $v_{k+1}$ depends on $\lambda_k$, since $v_{k}$ is the minimizer of $L^k(u)$, and the rule seems impractical. Nonetheless,
by the definition of the dual potential $\Phi_\mu^*(\mathcal{A}(u_k)-f)$, we have
\begin{align}
-L^k(v_{k+1})&=-\min_{u\in\lambda_k^{-1}(u_k+V_{\mu})}L^k(u)=-\min_{v\in V_{\mu}}\{\Phi_{\mu}(v)+\outerp{v}{{\mathcal{A}(u_{k})}-f}\}\nonumber\\
&=\max_{v\in V_{\mu}}\{\outerp{v}{f-{\mathcal{A}(u_{k})}}-\Phi_{\mu}(v)\}=\Phi_\mu^*(f-{\mathcal{A}(u_{k})}).\label{eqn:lam=Phi*}
\end{align}
i.e., $\lambda_{k} = c_{\rm s}\Phi_\mu^*(f-{\mathcal{A}(u_{k})})^\alpha$, which is independent of the step size $\lambda_k$. Moreover, it can be computed directly using the current approximation $u_k$.
\end{Remark}

Theorem \ref{thm:modelconv} shows that the error $\|u_k-u^*\|_X$ is monotonically decreasing until the $k^{*}$th step, and the convergence rate can be guaranteed before $k^{**}\leq k^{*}$.
The convergence of IDRM is determined by two factors.
The first factor is the convergence rate, and Theorem \ref{thm:modelconv} gives an exponential rate $(\eta q)^{k}$ and an algebraic rate $((1-\eta)k)^{-\frac{1}{\beta-1}}$ in the context of Hilbert and Banach spaces, respectively.
The second factor {arises from the stopping criteria $k^{*}$ and $k^{**}$, which is determined by the learning error} $\epsilon_M$. This factor is related to the approximation
capacity of the DNN class $\mathcal{N}$ and the quadrature error (e.g., via Monte Carlo). Due to the presence of the learning error $\epsilon_M$,
the difference $\xnorm{u_{k}-u_{k-1}}$ between two successive approximations does not tend to zero.
Indeed, once the error $\xnorm{u_{k}-u^*}$ reaches a level comparable with $\epsilon_M$ (i.e., $k\geq k^{*}$), the iterative process will no longer improve the accuracy of the approximation $u_k$.

Next we discuss the learning error $\epsilon_M$. For neural PDE solvers, the learning error refers to the gap between the DNN prediction $\widehat{u}_k$
and the exact solution $v_{k+1}$, and is represented by  $\epsilon_M$.
It arises from several different sources, including approximation capability of DNNs, quadrature rules to compute the loss and inexact minimization procedure used
in the training. The theoretical  analysis of neural PDE solvers has been carried out in many works (see, e.g., \cite{ShinDarbonKarniadakis:2020,JiaoLai:2022cicp,MishraMolinaro:2023,LuChenLu:2021,HuJinZhou:2022,Zeinhofer:2023unified} and the review \cite{Mishra:2024review} and the references therein). Nonetheless, the theoretical studies on general elliptic problems is still scarce, and we adapt the existing analysis to the current context. Below we denote the learning error $\epsilon_M$ of IDRM to be the difference between the population loss $L^k(v_{k+1})$ at the empirical minimizer $\widehat{u}_{k+1}\in\mathcal{N}$ and the true minimizer $v_{k+1}$, i.e., $\epsilon_M = L^{k}(\widehat{u}_{k+1})-L^{k}(v_{k+1})$. To bound this quantity, we need to make one more assumption on $\mathcal{A}$ and $f$.
\begin{Assumption}\label{assump:boundaf}
$a_{m}$ is a Nemytskii mapping: $\left|a_{m}(x, u(x), \nabla u(x))\right|\leq \ell_{a}(|u(x)|^{p-1}+|\nabla u(x)|^{p-1})$, and the sources $f_m\in L^\infty(\Omega)$, $m=0,1$.
\end{Assumption}

Below the notation $\mathbb{E}$ denotes taking expectation with respect to the random samples $\mathbb{X}=\{X_i\}_{i=1}^N$ and $\mathbb{T}=\{T_i\}_{i=1}^N$ used in the construction of the empirical loss $\widehat{L}^k$.

\begin{Theorem}\label{thm:modelerr}
Let the unique minimizer $v_{k+1}$ to the surrogate loss $L^k(u)$ over $X$ belong to $W^{1+r,p}(\Omega)$ with $r>1$. Then for any $\varepsilon>0,\nu>0$, there exists a DNN in the class
$$\mathcal{N}(C\log(d + 1 + r),C\varepsilon^{-\frac{d}{r-\nu}},C\varepsilon^{-\frac{d}{r-\nu}(3(1+\frac{1}{p})+\frac{2r+2}{d})}),$$
provided that the number $N$ of sampling points is taken to be $N\sim\mathcal{O}(\varepsilon^{-\frac{2d}{r-\nu}(\max\{3L-2,4\}+9L(1+\frac{1}{p})+\frac{6r+6}{d}L)-1})$, the learning error can be bounded by
\begin{equation}
\mathbb{E}_{\mathbb{X},\mathbb{T}}[\epsilon_M]\le C\lambda_k\varepsilon.
\end{equation}
\end{Theorem}

Theorem \ref{thm:modelerr} indicates that the learning error $\epsilon_M$ can be made arbitrarily small by taking the DNN $u_\theta\in\mathcal{N}$ sufficiently wide and the number $N$ of sampling points in the domain sufficiently large. The proof of Theorem \ref{thm:modelerr} relies on decomposing the learning error $\epsilon_M$ into the approximation and statistical errors in Lemma \ref{lem:decomp}. The approximation error is bounded using DNN approximation theory (cf. Proposition \ref{prop:approx}). The statistical error is bounded using Rademacher complexity, Lipschitz continuity of the DNN function class $\mathcal{N}$, and Dudley's entropy formula in Proposition \ref{prop:sta}. See Section \ref{ssec:learningerror} for the detail.

\section{Numerical experiments and discussions}\label{sec:experiment}

We now present numerical examples to illustrate IDRM. Unless otherwise specified, we select uniformly at random $N_r=10,000$ points in the domain $\Omega$ and $N_b=800$ points on the boundary $\partial\Omega$ to form the empirical loss $\widehat{L}^k(u_\theta)$. In all examples, we enforce the Dirichlet boundary condition $u=g$ on $\partial\Omega$ using the standard penalty formulation: $
L^k_\sigma(u)=L^k(u) + \sigma\|u- g\|_{L^2(\partial\Omega)}^2$,
where $\sigma>0$ is the penalty parameter. The numerical experiments focus on elliptic PDEs with nonsmooth solutions, but neither very smooth nor strongly singularly solutions. Although not presented, for PDEs with smooth solutions, both PINN and IDRM work fairly well, whereas for strongly singular solutions, special treatment is required in order to make the methods effective. To measure the accuracy of an approximation $\widehat u$ of the exact solution $u$, we use the relative $L^2(\Omega)$-error  $e=\|u-\widehat u\|_{L^2(\Omega)}/\|u\|_{L^2(\Omega)}$. We minimize the empirical loss $\widehat{L}_\sigma^k(u_\theta)$ using the Adam algorithm \cite{KingmaBa:2015}, with the default setting for the momentum parameters. For IDRM, the algorithm is applied to the loss (for each fixed $k$), with the initial guess set to the DNN parameter from the last loop for $k-1$ (so that there is no direction interaction between the iterates from different loops). The Python source codes for reproducing all the numerical experiments will be made available at the github repository \url{https://github.com/hhjc-web/IDRM}. First, we apply IDRM to a non-self-adjoint second-order elliptic problem, adapted from \cite[Sections 4.2.2 and 4.2.4]{Zang2020}.
\begin{Example}\label{exam1:elli}
Let $\Omega=(0,1)^{10}$. Consider the following boundary value problem:
\begin{align}\label{eqn:exam1:e}
\left\{\begin{aligned}
-\Delta u + \vec{b}\cdot \nabla u + \tfrac{\pi^2}{4}u &= f,&& \mbox{in }\Omega,\\
u&=g,&& \mbox{on }\partial\Omega,
\end{aligned}\right.
\end{align}
with the vector $\vec{b}=(1,0,\cdots,0)\in\mathbb{R}^{10}$.  The source $f$ is chosen such that the exact solution $u$ is given by $u(\boldsymbol{x})={\sum_{i=1}^{10}\min\{\sin(\frac{\pi}{2}x_i),0.9\}}$.
\end{Example}

\begin{figure}[hbt!]
\centering\setlength{\tabcolsep}{2pt}
\begin{tabular}{ccc}
\includegraphics[height=4cm]{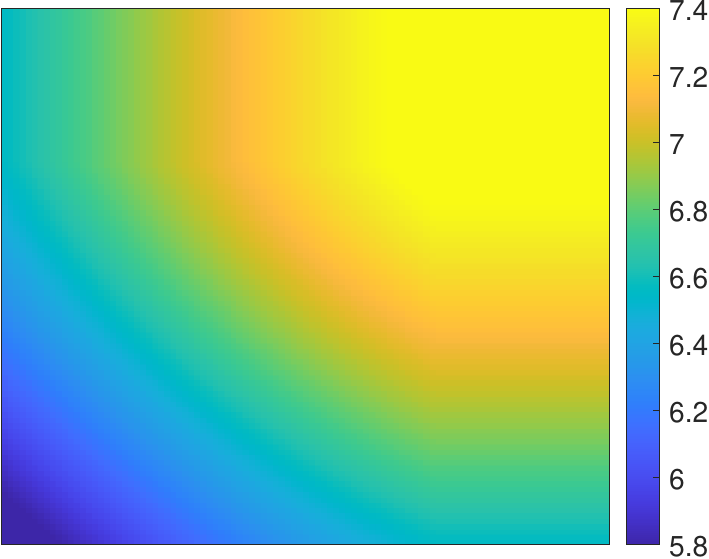} & \includegraphics[height=4cm]{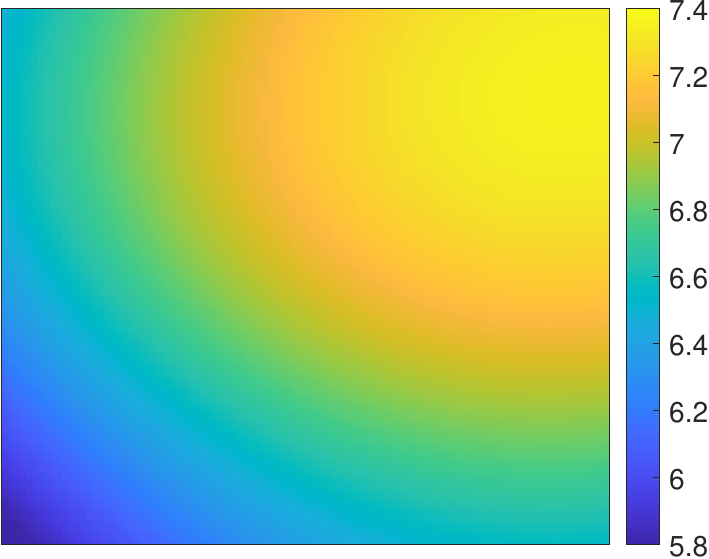} & \includegraphics[height=4cm] {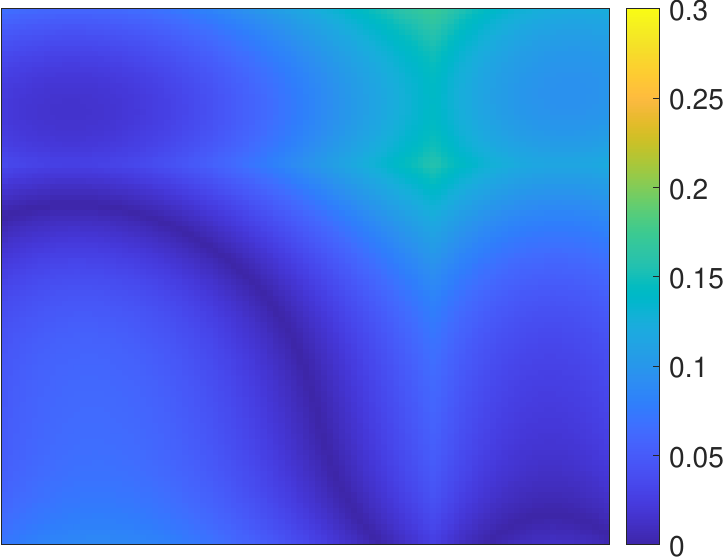}\\
\includegraphics[height=4cm]{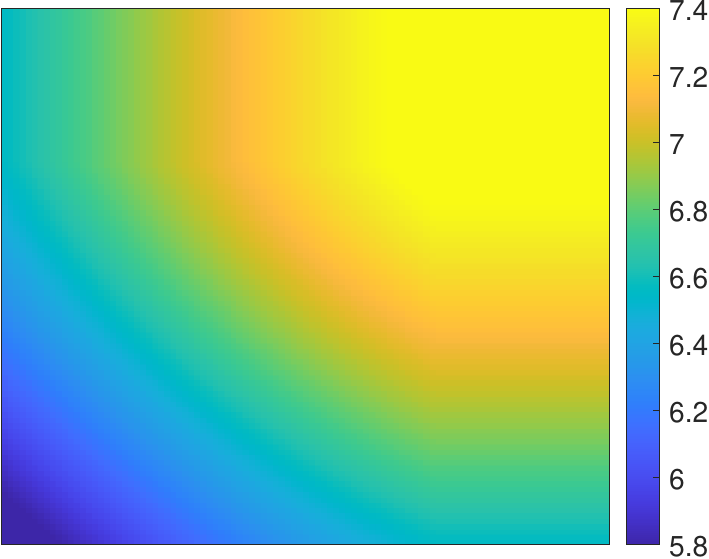} & \includegraphics[height=4cm]{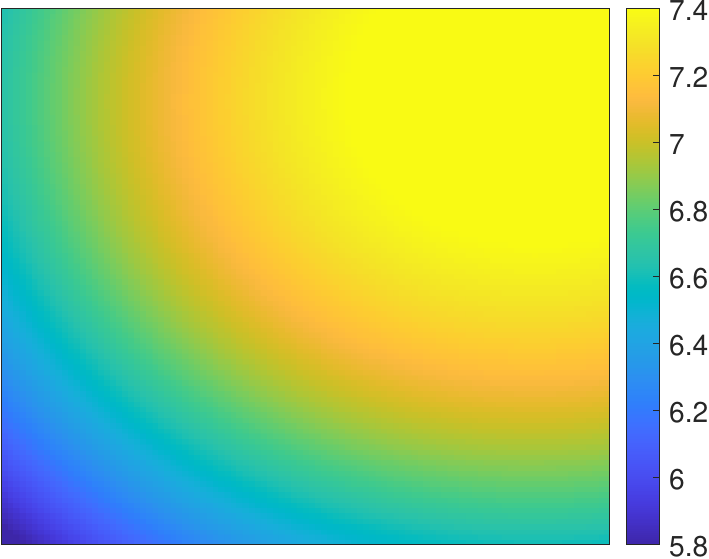} & \includegraphics[height=4cm]  {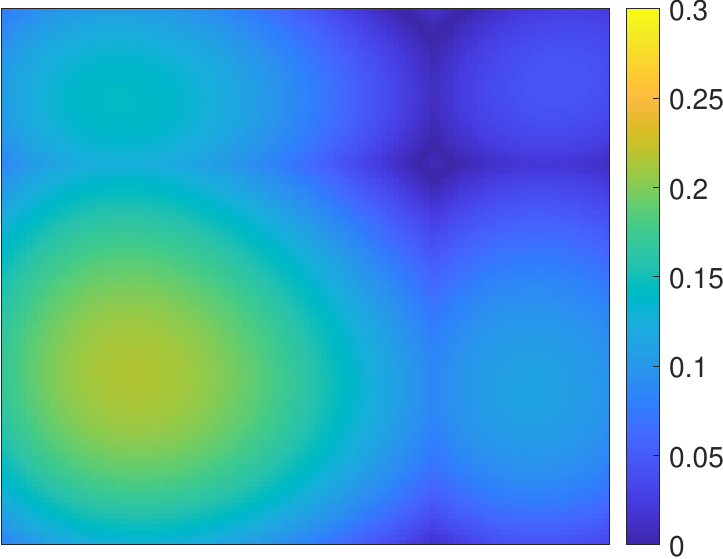}\\
\includegraphics[height=4cm]{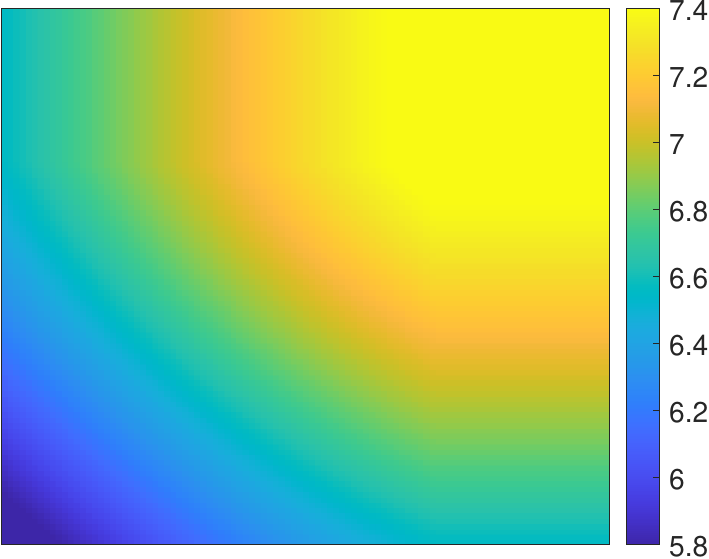} & \includegraphics[height=4cm]{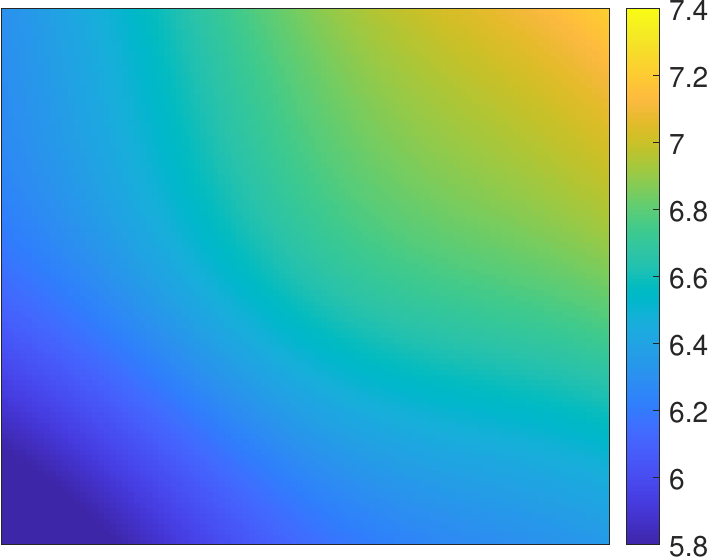} & \includegraphics[height=4cm]  {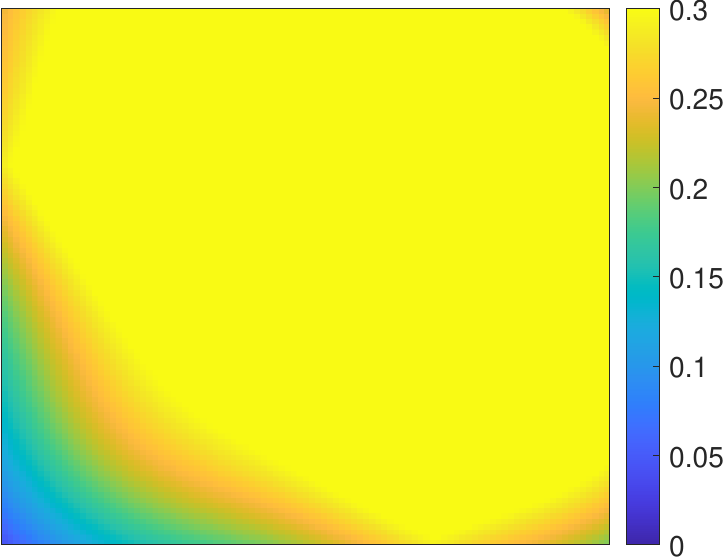}\\
(a) exact & (b) predicted  & (c) error
\end{tabular}
\caption{\label{fig:exam:elli} The DNN approximations for Example \ref{exam1:elli},  slices at $x_i=\frac{1}{2}\ (i=3,4,\cdots,10)$. (top: IDRM, middle: WAN, bottom: PINN)}
\end{figure}

Due to the presence of the convection term $\vec{b}\cdot\nabla u$, this example does not have a Ritz formulation and hence the standard deep Ritz method \cite{yu2018deep} cannot be applied directly. The exact solution $u$ belongs to $H^1(\Omega)$ but not
$H^2(\Omega)$, which makes the PINN loss ill-defined (at the population level), but IDRM and WAN require only the
 existence of the first weak derivative and both are well defined. In IDRM, we fix $\lambda_0=1$ and $\mu=0$ in the loss $\widehat{L}^k$, and employ a DNN architecture 10-20-20-20-1 to approximate the solution $u$. The parameters and relative errors by IDRM, WAN and PINN are shown in Table \ref{table:exam1:elli}, where the results correspond to the last iteration of the training process (instead of the minimum error along the training trajectory). Fig. \ref{fig:exam:elli} (top) shows that IDRM can give an accurate approximation of the solution $u$, where we also show the results by PINN and WAN. The numerical results also show that PINN fails to approximate accurately the non-smooth part of the solution $u$, cf. Fig. \ref{fig:exam:elli} (bottom), and both IDRM and WAN can yield more accurate approximations.
However, the training of the WAN loss is more tricky, as was observed in early works \cite{Uriarte:2023,bertoluzza2024best}.

\begin{table}[hbt!]
\centering
\begin{threeparttable}
\caption{\label{table:exam1:elli} The parameters and results of different methods for Example \ref{exam1:elli}. The reported results are computed from the last iteration of the training for each method.}
\centering
\begin{tabular}[5pt]{c|c|c|c|c|c}
\toprule
Method &  Network structure & Learning rate & Penalty $\sigma$ & Iterations & Relative error \\
\midrule
IDRM & 10-20-20-20-1 & 3.0e-3 & 100.0 & 4.1k & 9.70e-3 \\
\hline
WAN & \makecell{trial net: 10-20-20-20-1\\
test net: 10-20-20-20-1} & 5.0e-3 & 100.0 & 3.4k & 1.80e-2\\
\hline
PINN & 10-20-20-20-1 & 5.0e-3 & 100.0 & 5.0k & 5.62e-2\\
\bottomrule
\end{tabular}
\end{threeparttable}
\end{table}

\begin{figure}[hbt!]
\centering\setlength{\tabcolsep}{0pt}
\begin{tabular}{cc}
\includegraphics[height=6cm]  {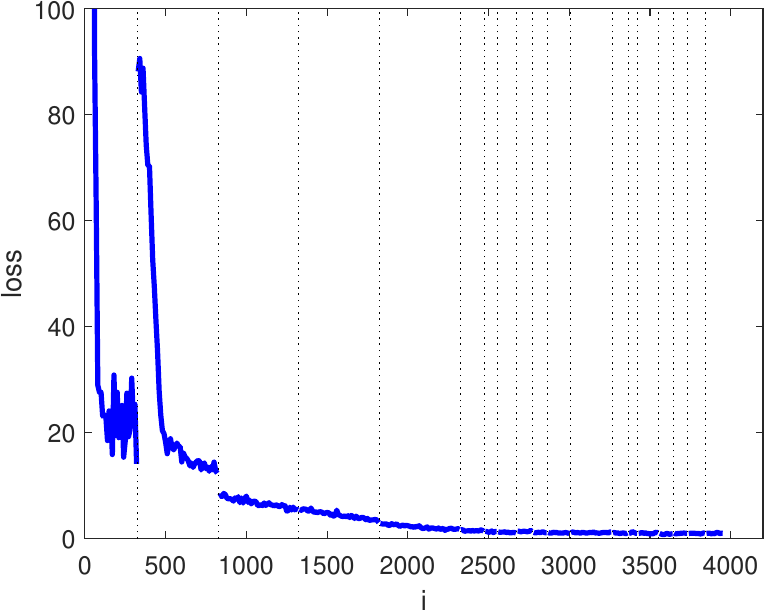} &  \includegraphics[height=6cm]  {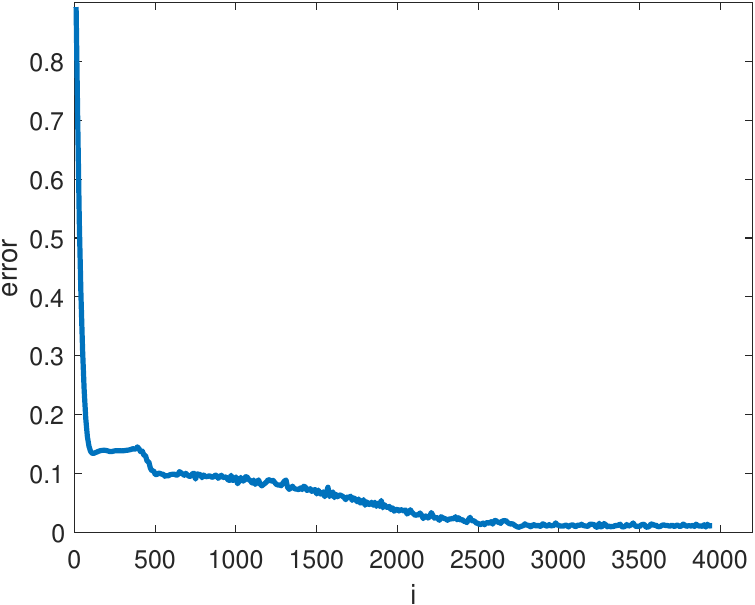}\\
(a) $L_\sigma^k$ vs $i$ & (b) $e$ vs $i$
\end{tabular}
\caption{\label{fig:exam:elli-ite} The training dynamics of IDRM for Example \ref{exam1:elli}: {\rm(a)} the loss $L_\sigma^k$ versus the iteration index $i$, {\rm(b)} the error $e$ versus the iteration index $i$.}
\end{figure}

To shed further insights, we show in Fig. \ref{fig:exam:elli-ite} the training dynamics of the loss $\widehat{L}^k(u_\theta)$ and the relative error $e_{\rm idrm}$, where $i$ denotes the total iteration index. Note that the loss $\widehat{L}^k(u_\theta)$ changes with the iteration index $k$, and thus it undergoes jumps as $ k $ changes, cf. Fig. \ref{fig:exam:elli-ite}(a). Moreover, its value decreases most during the first two loops, and as the loop further proceeds, it does not decrease much. This is attributed to the following fact: In the first few loops, the prediction $ \hat{u}_k$ is far away from the exact one (cf. Fig. \ref{fig:exam:elli-ite}(b)) and training can greatly reduce the error $e_{\rm idrm}$; at the third loop, the approximation is already close to the exact one (cf. Fig. \ref{fig:exam:elli-ite}(b)), so the decay of the loss $\widehat{L}^k(u_\theta)$ is small, and the loss stays at the same level.
From Fig. \ref{fig:exam:elli-ite}(b), the error $e_{\rm idrm}$ drops very quickly during the initial loops,  but as the loop proceeds, it exhibits some oscillations but stays at a relatively low level \cite{RAISSI2019686,CuomoSchiano:2022}. While the overall trend of the loss and the error is comparable, a smaller loss does not indicate a smaller error. This is attributed to the following facts: the loss is computed only on the training points, whereas the error is evaluated on the testing points; the loss combines PDE and boundary residuals, whereas the error is only about the approximation of the solution in the domain.

The next example taken from \cite{Liu2023} is about the $p$-Laplace equation.
\begin{Example}\label{exam2:pLap}
Let $\Omega=(0,1)^{10}$. Consider the following $p$-Laplace problem:
\begin{align}
\left\{\begin{aligned}
-\nabla\cdot(|\nabla u|^{p-2}\nabla u) + \vec{b}\cdot \nabla u &= f,&& \mbox{in }\Omega,\\
u&=g,&& \mbox{on }\partial\Omega,
\end{aligned}\right.
\end{align}
and the following two problem settings.
\begin{itemize}
\item[{\rm(i)}] $p=2.5$. The source $f({x})=r^{\frac{1}{2}}-(\frac{2}{5})^{\frac{2}{3}}x_1$, with $r=\|{x}\|_2$ and the vector $\vec{b}=(1,0,\cdots,0)\in\mathbb{R}^{10}$.  The exact solution $u$ is given by $u(\boldsymbol{x})=\frac{1}{\sqrt[3]{50}}(1-r^2)$.
\item[{\rm(ii)}] $p=1.5$. The source $f({x})=\tfrac{9-x_1}{r}$, with $r=\|{x}\|_2$ and the vector $\vec{b}=(1,0,\cdots,0)\in\mathbb{R}^{10}$.  The exact solution $u$ is given by $u({x})=1-r$.
\end{itemize}
\end{Example}

\begin{figure}[hbt!]
\centering\setlength{\tabcolsep}{2pt}
\begin{tabular}{ccc}
\includegraphics[height=4cm]  {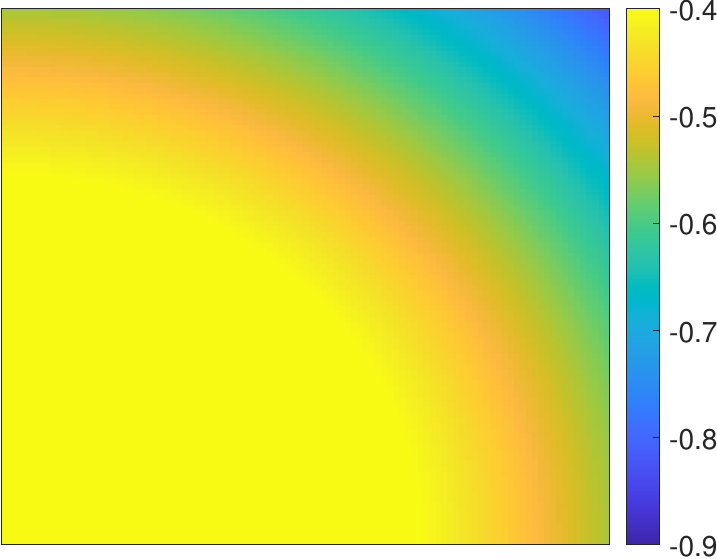} & \includegraphics[height=4cm]  {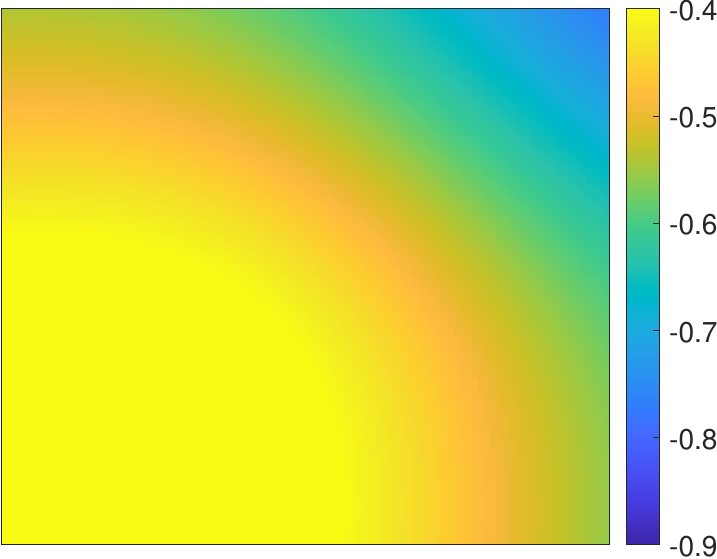} & \includegraphics[height=4cm]  {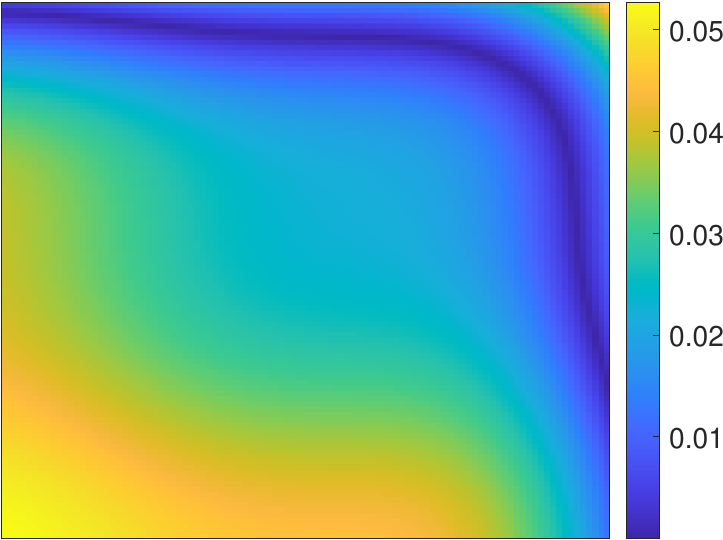}\\
(a) exact & (b) predicted  & (c) error
\end{tabular}
\caption{\label{fig:exam2:pLap0} The DNN approximation for Example \ref{exam2:pLap}(i),  slices at $x_i=\frac{1}{2}(i=3,4,\cdots,10)$.}
\end{figure}

\begin{figure}[hbt!]
\centering\setlength{\tabcolsep}{0pt}
\begin{tabular}{cc}
\includegraphics[height=6cm]  {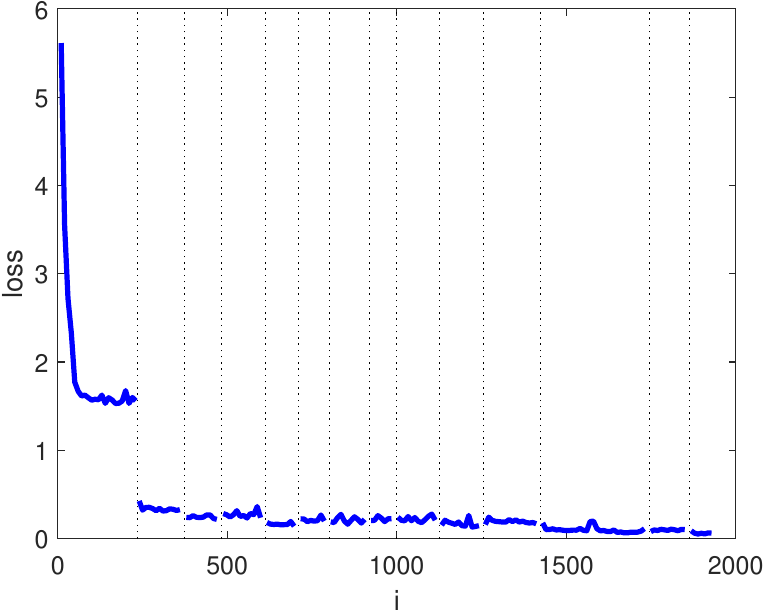} &  \includegraphics[height=6cm]  {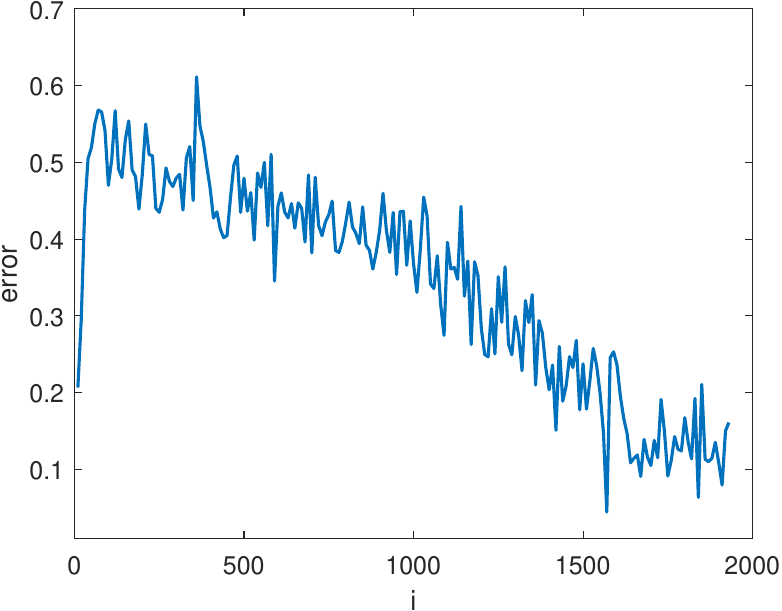}\\
(a) $L_\sigma^k$ vs $i$ & (b) $e$ vs $i$
\end{tabular}
\caption{\label{fig:exam:pLap-ite0} The training dynamics of IDRM for Example \ref{exam2:pLap}(i): {\rm(a)} the loss $L_\sigma^k$ versus the iteration index $i$ and {\rm(b)} the error $e$ versus the iteration index $i$.}
\end{figure}

We use case (i) to illustrate the performance of IDRM for $X=W^{1,p}(\Omega)$, with $p>2$, which is covered by Corollary \ref{cor:plarge}. We employ a DNN architecture 10-50-50-50-50-50-50-1 to approximate the solution $u$, and set $\lambda_0=1.0$ and $\mu=0$ in the loss $\widehat{L}^k(u_\theta)$. The DNN approximation has a relative error $e_{\rm idrm}=\text{4.95e-2}$; See Fig. \ref{fig:exam2:pLap0} for the approximation and Fig. \ref{fig:exam:pLap-ite0} for the training dynamics.

\begin{figure}[hbt!]
\centering\setlength{\tabcolsep}{2pt}
\begin{tabular}{ccc}
\includegraphics[height=4cm]  {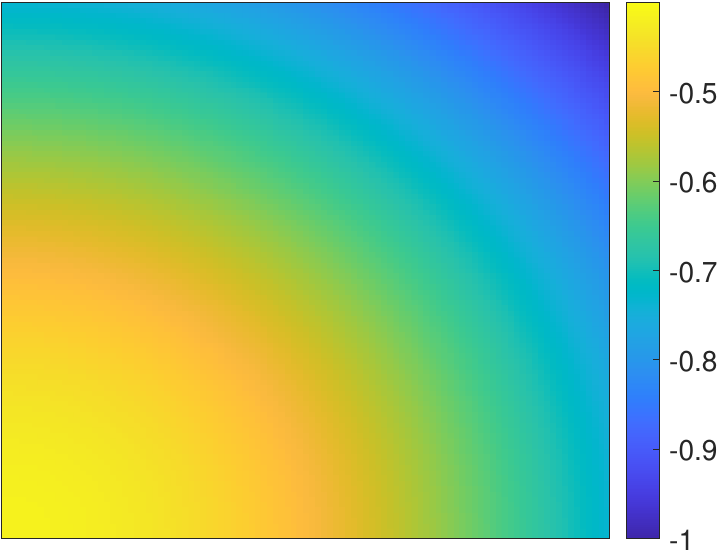} & \includegraphics[height=4cm]  {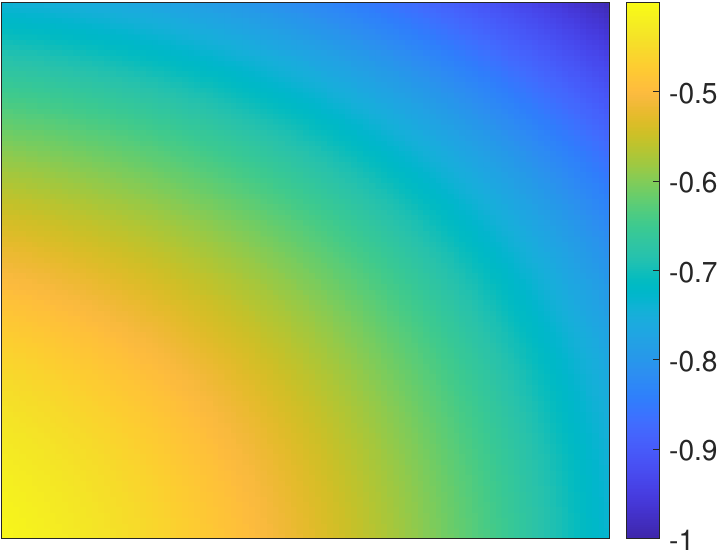} & \includegraphics[height=4cm]  {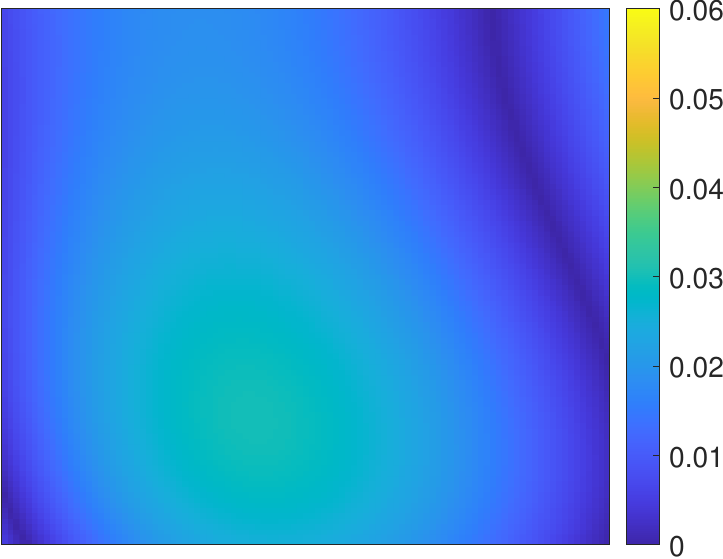}\\
\includegraphics[height=4cm]  {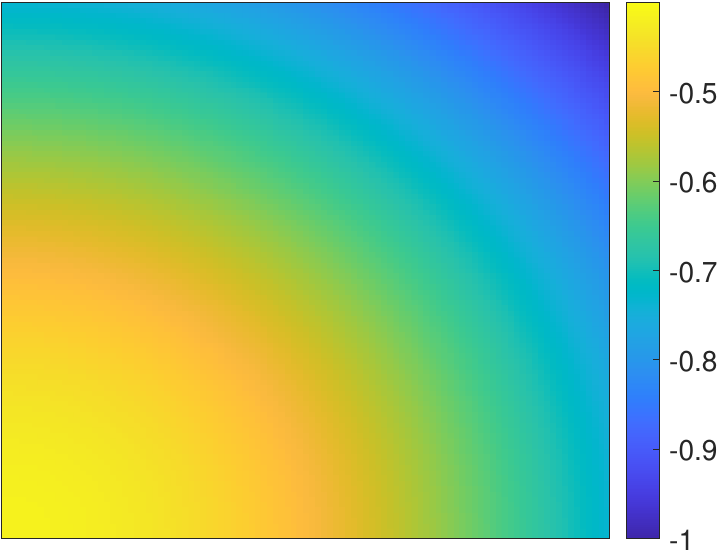} & \includegraphics[height=4cm]  {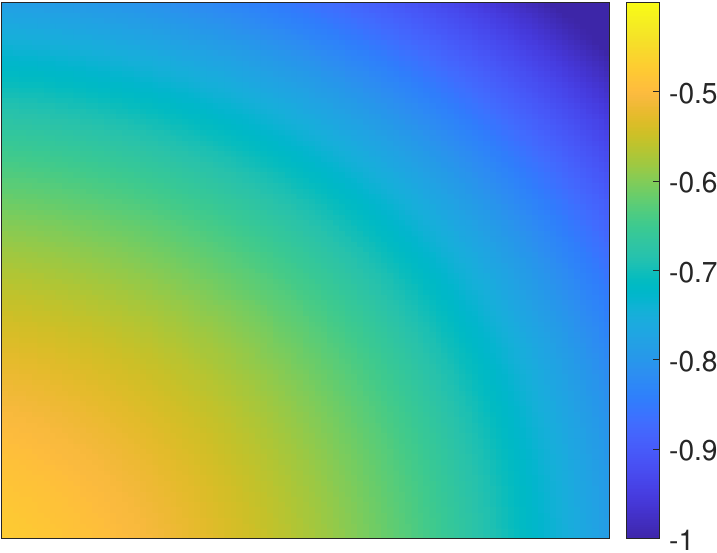} & \includegraphics[height=4cm]  {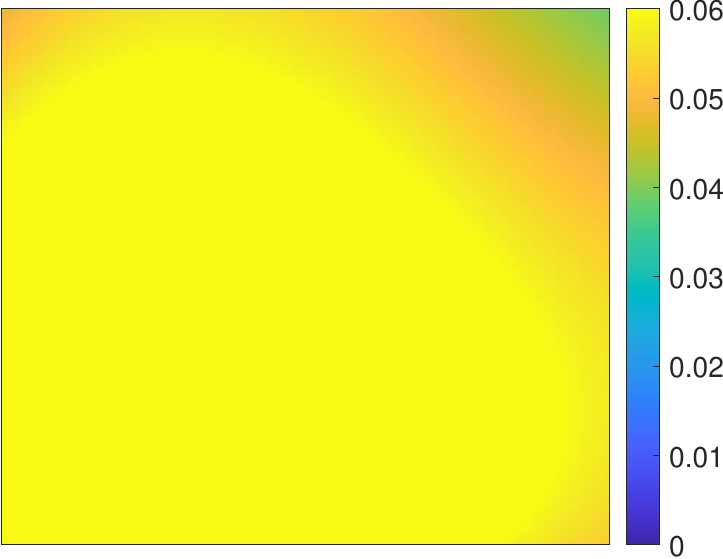}\\
\includegraphics[height=4cm]  {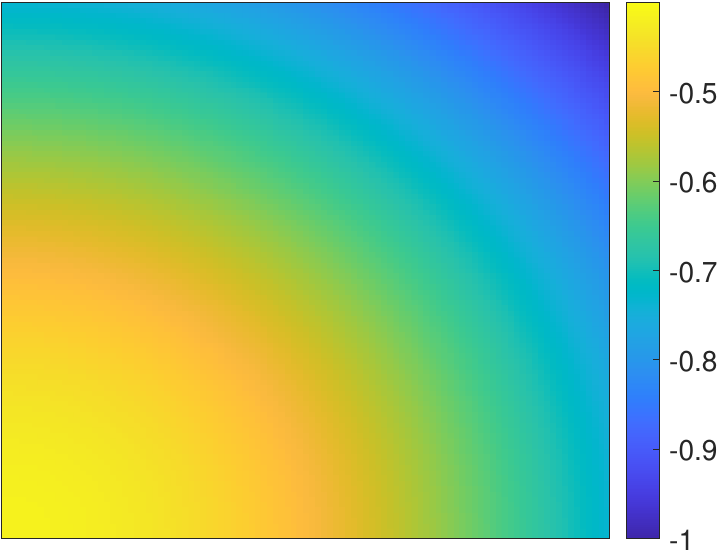} & \includegraphics[height=4cm]  {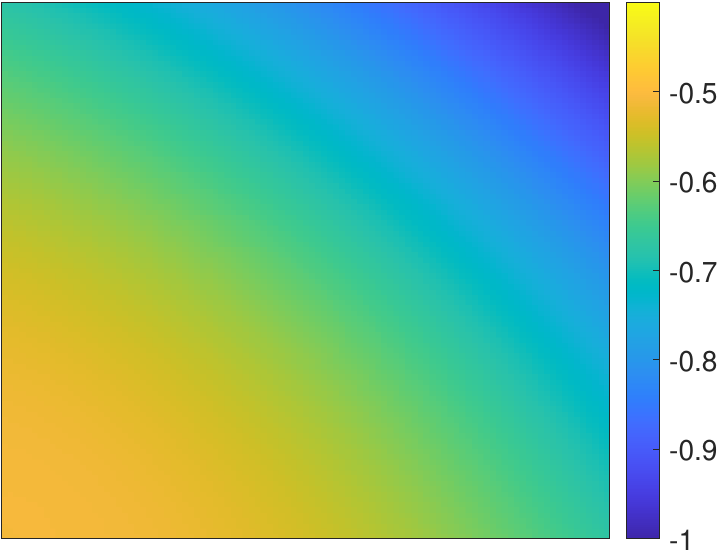} & \includegraphics[height=4cm]  {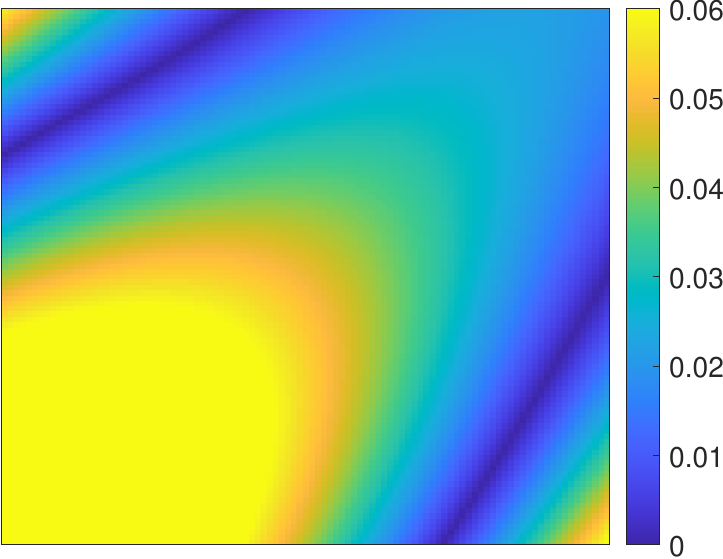}\\
(a) exact & (b) predicted  & (c) error
\end{tabular}
\caption{\label{fig:exam2:pLap} The DNN approximations for Example \ref{exam2:pLap}(ii), slices at $x_i=\frac{1}{2}(i=3,4,\cdots,10)$ (top: IDRM, middle: WAN, bottom: PINN).}
\end{figure}

\begin{figure}[hbt!]
\centering\setlength{\tabcolsep}{0pt}
\begin{tabular}{cc}
\includegraphics[height=6cm]  {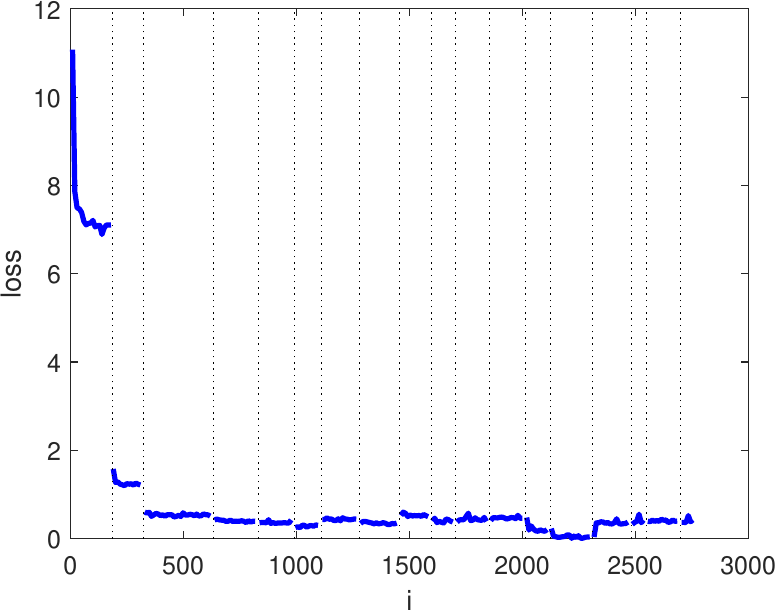} &  \includegraphics[height=6cm]  {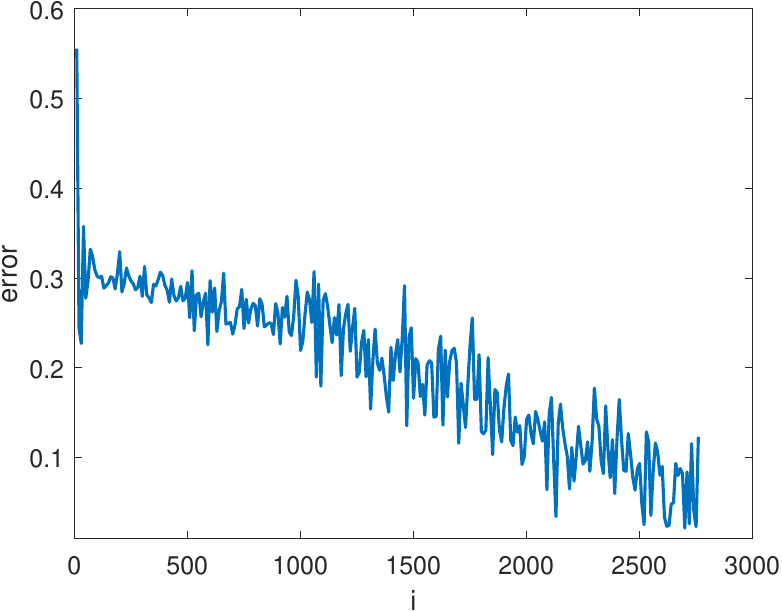}\\
(a) $L_\sigma^k$ vs $i$ & (b) $e$ vs $i$
\end{tabular}
\caption{\label{fig:exam:pLap-ite} The training dynamics of IDRM for Example \ref{exam2:pLap}(ii): {\rm(a)} the loss $L_\sigma^k$ versus the iteration index, {\rm(b)} the error $e$ versus the iteration index.}
\end{figure}

For case (ii), we take $X=W^{1,p}(\Omega)$, and the potential $\Phi_\mu(u)$ to be
$$\Phi_\mu(u)=\frac{1}{p}\int_\Omega|\nabla u|^p\,{\rm d}{x} + \int_0^1 \int_\Omega u(\vec{b}\cdot \nabla (tu))\,{\rm d}x{\rm d}t + \mu\|u\|_{W^{1,p}(\Omega)}^2.$$
Since $p=1.5<2$, computing the duality product $\outerp{v}{\mathcal{A}(u)}$ directly involves a singular integral, with a singularity at the origin. To resolve the issue, numerically we approximate it using the finite difference as
$$\outerp{v}{\mathcal{A}(u)}\approx\dfrac{1}{p}\int_\Omega\big(|\nabla (u+\delta v)|^p-|\nabla u|^p+ v(\vec{b}\cdot \nabla u)\big)\,{\rm d}{x},$$
where $\delta>0$ is small, fixed at $\delta=0.01$ below. We use IDRM, WAN and PINN to solve the problem. We employ a DNN architecture 10-16-32-32-16-1 to approximate the solution $u$ in IDRM and PINN, and set $\lambda_0=1$ and $\mu=2$ in the empirical IDRM loss $\widehat{L}^k(u_\theta)$.  The point-wise error of IDRM is small in the domain $\Omega$, cf. Fig. \ref{fig:exam2:pLap}. In contrast, applying either PINN or WAN directly yields large errors near the singularity point of the solution $u$, and thus fails to produce accurate approximations due to the low regularity of the solution $u$. The relative errors are shown in Table \ref{table:exam2:pLap}. The training dynamics in Fig. \ref{fig:exam:pLap-ite} shows that the loss $\widehat{L}_\sigma^k(u_\theta)$ decays fast in the first loop, but the decay slows down greatly as the iteration proceeds, which is consistent with the observation for Example \ref{exam1:elli}. Likewise, the error also suffers from oscillations in later loops.

\begin{table}[hbt!]
\centering
\begin{threeparttable}
\caption{\label{table:exam2:pLap} The parameters and results of different methods for Example \ref{exam2:pLap}(ii). The reported results are computed for the last iteration of the training for each method.}
\centering
\begin{tabular}[5pt]{c|c|c|c|c|c}
\toprule
Method &  Network structure & Learning rate & Penalty $\sigma$ & Iterations & Relative error \\
\midrule
IDRM & 10-16-32-32-16-1 & 5.0e-3 & 40.0 & 2.8k & 2.86e-2 \\
\hline
WAN & \makecell{trial net: 10-50-50-50-1\\
test net: 10-50-50-50-50-50-1} & 5.0e-3 & 40.0 & 3.3k & 1.03e-1\\
\hline
PINN & 10-16-32-16-1 & 5.0e-3 & 40.0 & 2.7k & 6.71e-2\\
\bottomrule
\end{tabular}
\end{threeparttable}
\end{table}

The next example is a high-dimensional parabolic problem.
\begin{Example}\label{exam3:heat}
Let $\Omega=(0,1)^{10}$, and fix $T>0$. Consider the following quasi-linear diffusion problem:
\begin{align}
\left\{\begin{aligned}
u_t &= \nabla\cdot(u\nabla u) + f,&& \mbox{in }\Omega\times(0,T],\\
u&=g,&& \mbox{on }\partial\Omega\times(0,T],\\
u&=u_0,&& \mbox{on }\Omega\times\{0\},
\end{aligned}\right.
\end{align}
with the exact solution $u$ given by $u({x},t)=(\sum_{i=1}^{10}x_i+1+t + \frac{1}{5}\sin(5\pi t))^{1/2}$.
\end{Example}

\begin{figure}[hbt!]
\centering\setlength{\tabcolsep}{2pt}
\begin{tabular}{ccc}
\includegraphics[height=4cm]  {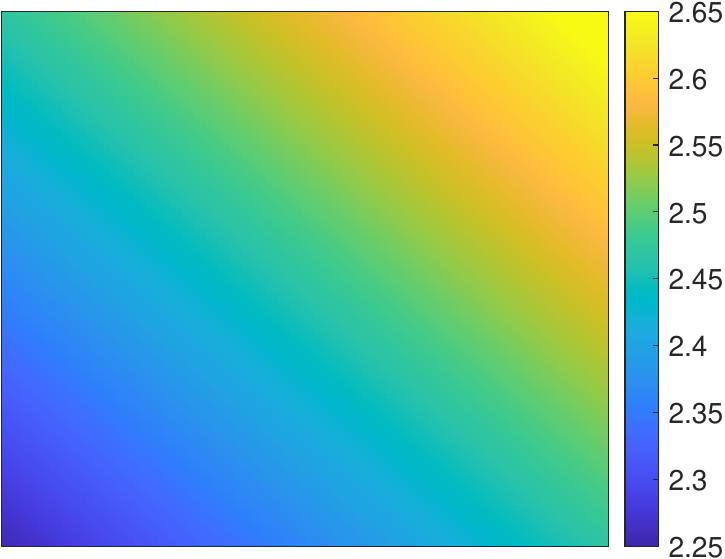} & \includegraphics[height=4cm]  {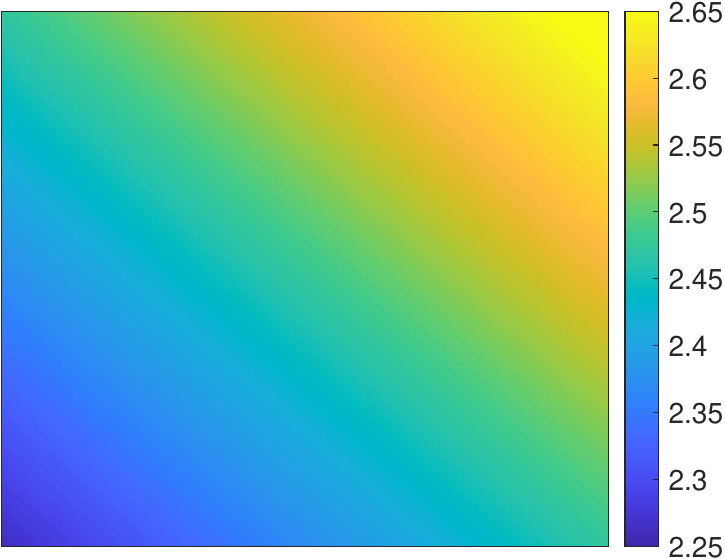} & \includegraphics[height=4cm]  {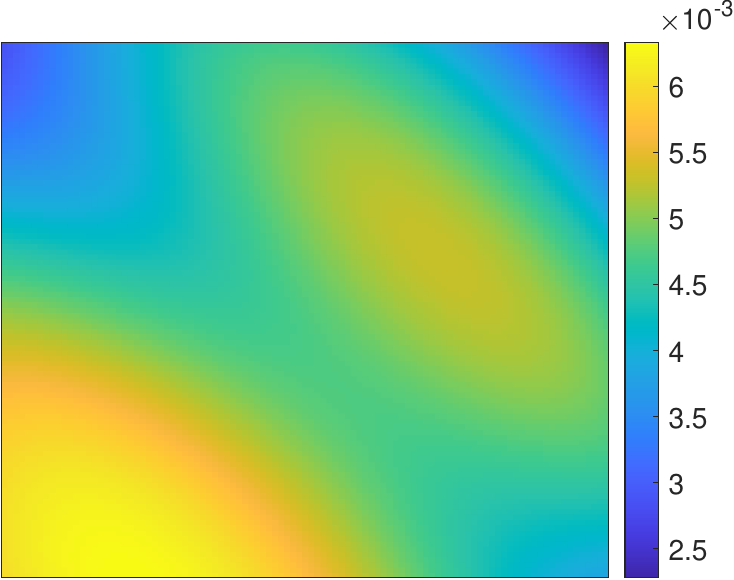}\\
(a) exact & (b) predicted  & (c) error
\end{tabular}
\caption{\label{fig:exam:heat} The DNN approximation for Example \ref{exam3:heat},  slices at $x_i=\frac{1}{2}(i=3,4,\cdots,10)$, at $t=0.3$.}
\end{figure}

\begin{figure}[hbt!]
\centering\setlength{\tabcolsep}{2pt}
\begin{tabular}{cc}
\includegraphics[height=5cm]  {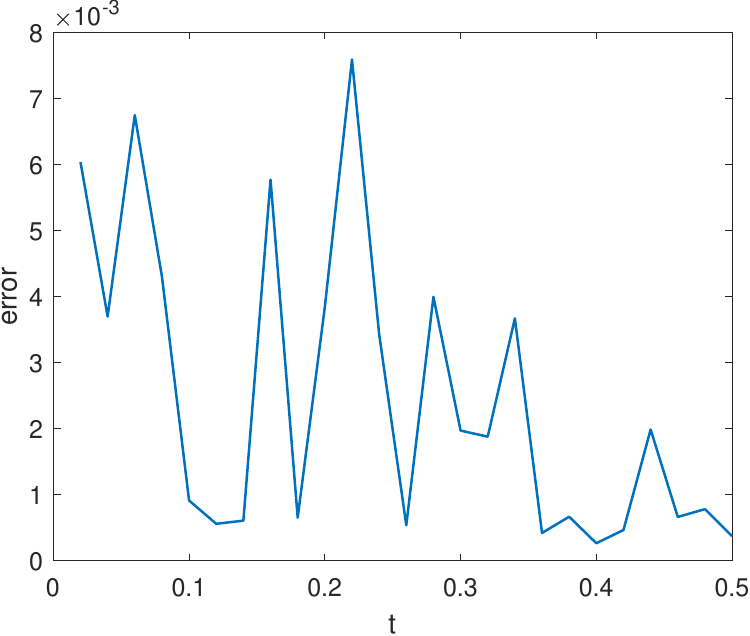} &  \includegraphics[height=5cm]  {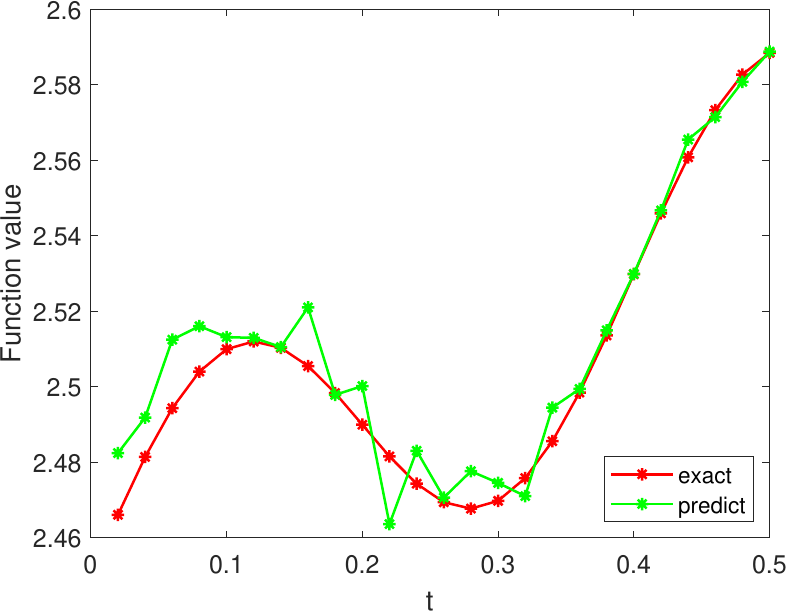}\\
(a) $e$ vs $t$ & (b) $e_0$ vs $t$
\end{tabular}
\caption{\label{fig:exam:heat-t} {\rm(a)} The evolution of the relative error $e_{\rm idrm}$, {\rm(b)} the variation of the exact and predicted function values at $\boldsymbol{x}=(1/2,1/2,1/2,\cdots,1/2)$.}
\end{figure}

We discretize the parabolic problem in time using the backward Euler method: we divide the time interval $[0,T]$ into $N$ equal sub-intervals with the partition $0=t_0<t_1<t_2<\cdots<t_N=T$, $t_k=k\Delta t$, $k=0,1,\cdots,N$, $\Delta t=T/N$ and let $u_k\equiv u({x},t_k)$. Then the time semi-discrete problem reads
$$\dfrac{u_{k+1}-u_k}{\Delta t} = \nabla\cdot(u_{k+1}\nabla u_{k+1}) + f({x}, t_{k+1}),\quad k=0,1,\cdots,N-1,$$
or equivalently
\begin{equation}\label{exam:para:timedis}
-\Delta t\nabla\cdot(u_{k+1}\nabla u_{k+1}) + u_{k+1} = \Delta t f({x}, t_{k+1}) + u_k,\quad k=0,1,\cdots,N-1.
\end{equation}
Then{\color{blue},} we employ IDRM to solve \eqref{exam:para:timedis}, and obtain the time stepping solutions $\{\widehat{u}_k\}_{k=1}^N$.
We fix $T=0.5$, employ a 10-20-40-40-20-1 DNN architecture to approximate the time-stepping solution $u_k$, and set $\lambda_0=1$
and $\mu=0.01$ in the loss $\widehat{L}_\sigma^k(u_\theta)$. Fig. \ref{fig:exam:heat} shows the exact solution $u_N$, prediction $\widehat{u}_N$ and pointwise
error $|u_N-\widehat u_N|$ at the terminal time $T$. Fig. \ref{fig:exam:heat-t}(a) shows the evolution of the relative error $e_{\rm idrm}$ with the time $t$ and Fig. \ref{fig:exam:heat-t}(b) shows the exact and predicted function values at one point.

The last example is about the three-dimensional Navier-Stokes problem.
\begin{Example}\label{exam4:stoke}
Let $\Omega=(0,1)^3$. Consider the three-dimensional Navier-Stokes equations:
\begin{align}
\left\{\begin{aligned}
-\mu\Delta \boldsymbol{u} + \nabla p + (\boldsymbol{u}\cdot\nabla) \boldsymbol{u} &= \boldsymbol{f},&& \mbox{in }\Omega,\\
\nabla\cdot \boldsymbol{u} &=0 ,&& \mbox{in }\Omega,\\
\boldsymbol{u}&=\boldsymbol{g},&& \mbox{on }\partial\Omega,
\end{aligned}\right.
\end{align}
with the viscosity coefficient $\mu = 0.1$, the pressure $p({x}) = x_1x_2x_3$ and the velocity field $\boldsymbol{u}({x})=(x_1x_3^{\frac{2}{3}}-\frac{2}{3}x_1x_2x_3^{-\frac{1}{3}},\frac{2}{3}x_1x_2x_3^{-\frac{1}{3}}-x_2x_3^{\frac{2}{3}},x_2x_3^{\frac{2}{3}}-x_1x_3^{\frac{2}{3}})^\top$.
\end{Example}

\begin{figure}[hbt!]
\centering \setlength{\tabcolsep}{2pt}
\begin{tabular}{ccc}
\includegraphics[height=4cm]{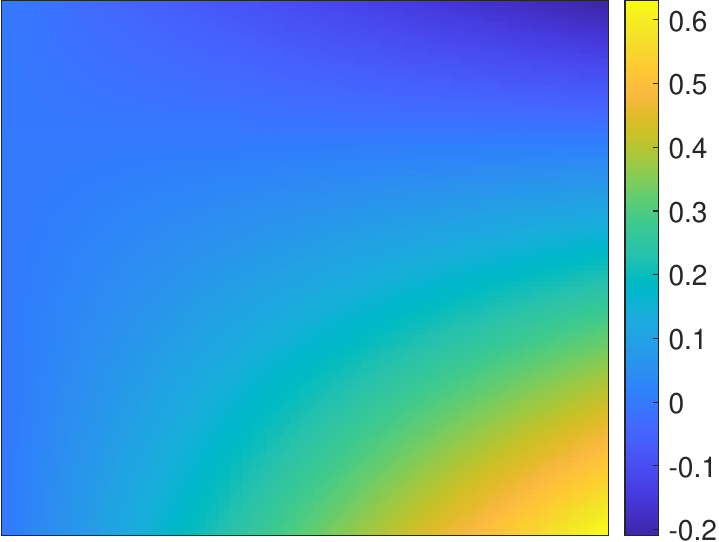}& \includegraphics[height=4cm]{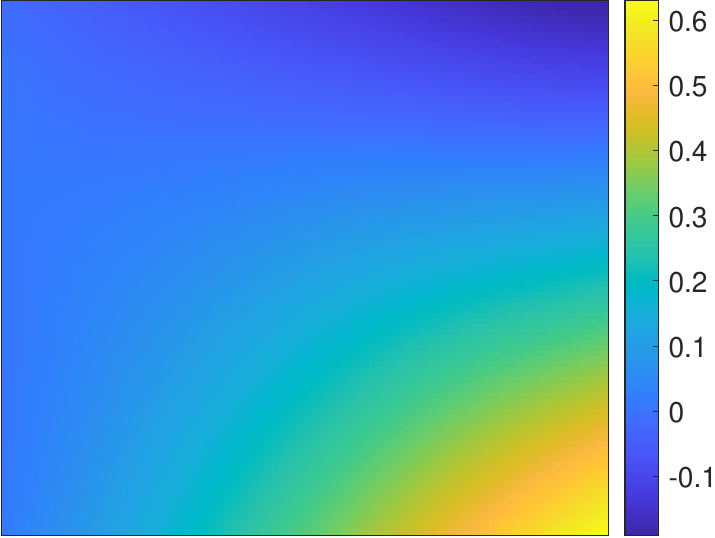}& \includegraphics[height=4.2cm]{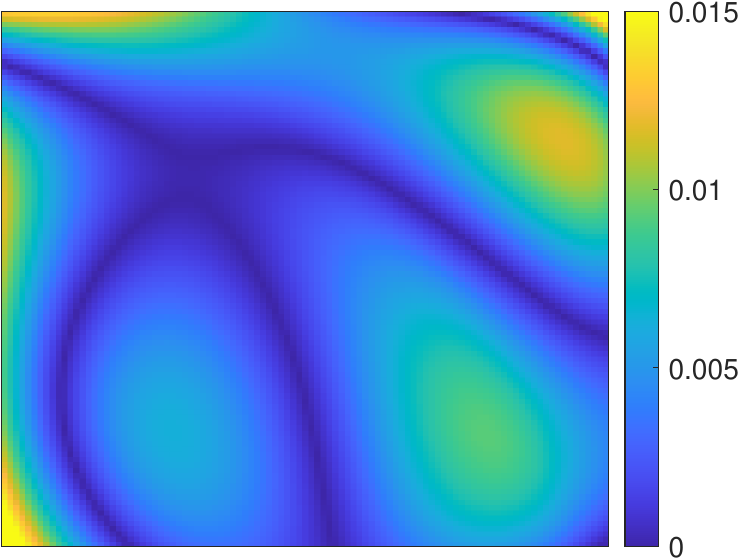}\\
\includegraphics[height=4cm]{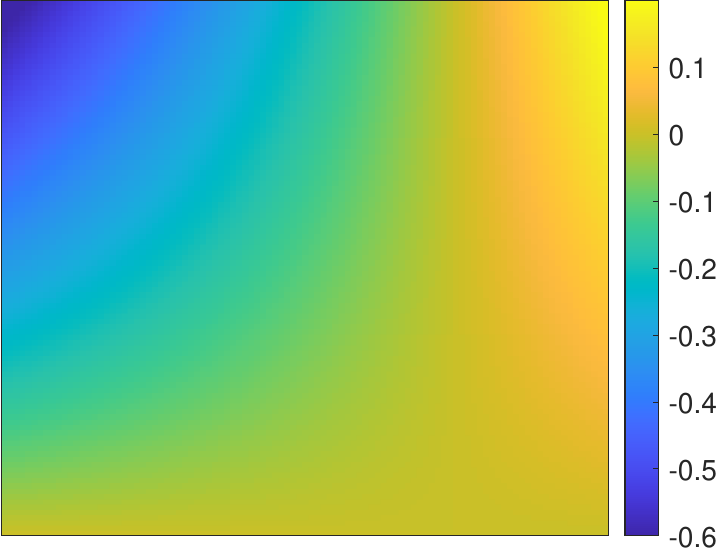} & \includegraphics[height=4cm]{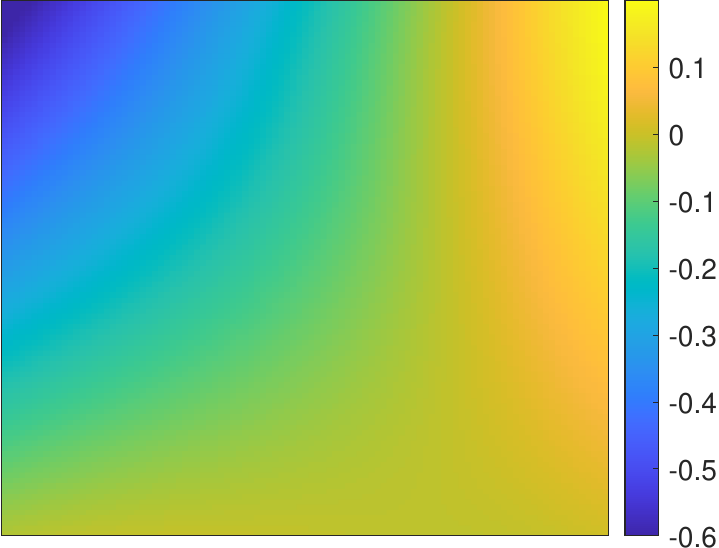}& \includegraphics[height=4.2cm]{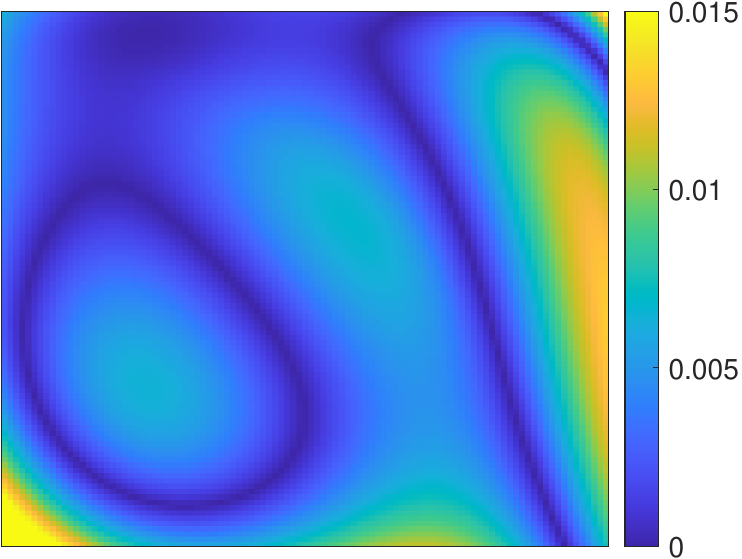}\\
\includegraphics[height=4cm]{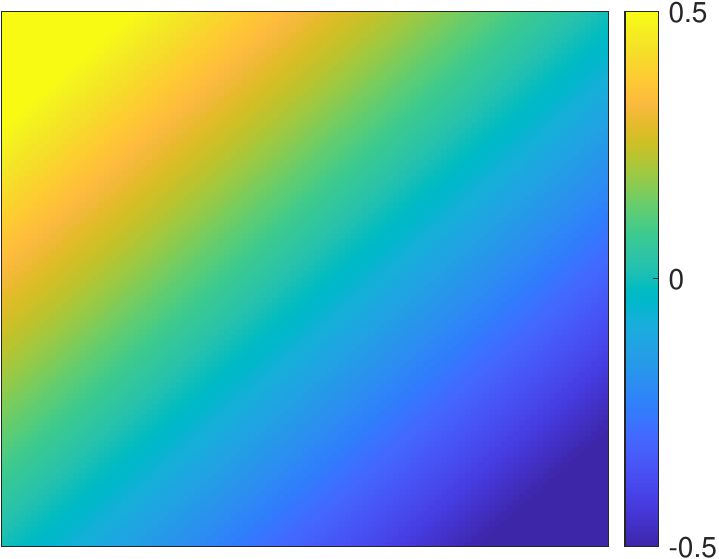} &
\includegraphics[height=4cm]{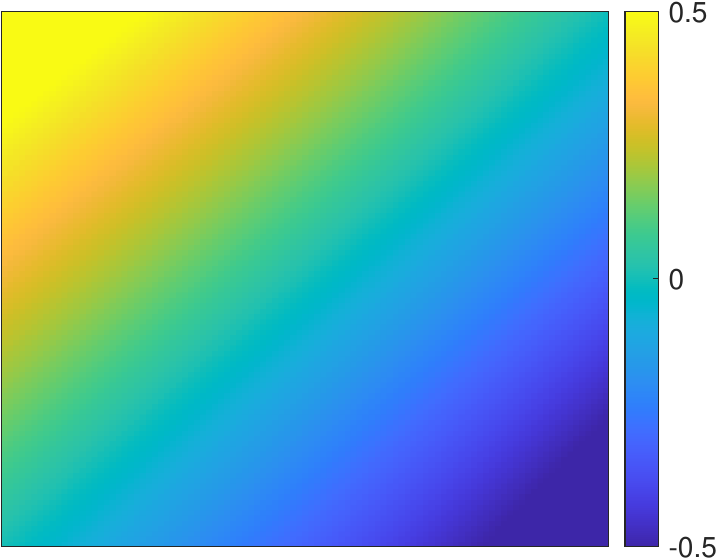} & \includegraphics[height=4.2cm]{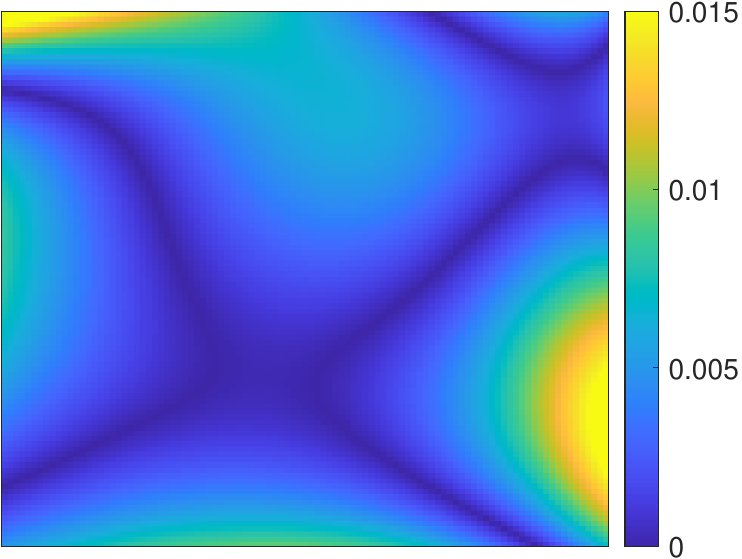}\\
exact & predicted & error \\
\end{tabular}
\caption{The DNN approximation of Example \ref{exam4:stoke}, from top to bottom: $u_1$, $u_2$, $u_3$, at the slice $x_3 = \frac{1}{2}$. \label{fig:exam4:stoke}}
\end{figure}

Recall that $\outerp{\boldsymbol{v}}{\mathcal{A}(\boldsymbol{u})}=\int_\Omega\mu\nabla\boldsymbol{u}:\nabla\boldsymbol{v}+ (\boldsymbol{u}\cdot\nabla) \boldsymbol{u}\cdot\boldsymbol{v}\,{\rm d}{x}$ with $\boldsymbol{v}\in (H_0^1(\Omega))^3$ and $\nabla\cdot\boldsymbol{v}=0$. To preserve the divergence-free condition $\nabla\cdot\boldsymbol{v}=0$, we use the rotation of a DNN to approximate the velocity field $\boldsymbol{u}({x})$, i.e., find a $\boldsymbol{u}_{k+1}\in \nabla\times\mathcal{N}$ such that
${L}^{k}(\boldsymbol{u}_{k+1})\leq \min_{\boldsymbol{u}\in \boldsymbol{u}_{k}+V_{\mu}}{L}^{k}(\boldsymbol{u})+\epsilon_M$
at each step of Algorithm \ref{alg:alg1}. We employ the DNN architecture $\mathcal{N}$ to be 3-10-10-1 and set $\lambda_0=1.0$ and $\mu=0.0$ in the loss $\widehat{L}^k(\boldsymbol{u}_\theta)$. Since the domain $\Omega$ is low-dimensional, we employ a uniform mesh with a grid size $h=0.01$ and the composite trapezoidal rule to approximate the integrals, which gives an accurate quadrature formula. To obtain a more accurate approximation, we adopt a path-following procedure \cite{AllgowerGeorg:2003,HuJinZhou:2022} on the loss:
$L^k_{\sigma_k}(\boldsymbol{u})=L^k(\boldsymbol{u}) + \sigma_k\|\boldsymbol{u}- T\boldsymbol{u}\|_{L^2(\partial\Omega)}^2$. We take a small value $\sigma_1=10.0$, and  after each loop update $\sigma$ geometrically: $\sigma_{k+1}=\rho\sigma_k$ with $\rho=1.5$. The pointwise error of the final reconstruction is shown in Fig. \ref{fig:exam4:stoke} with a relative error $e_{\rm idrm}=(\text{2.85e-2}, \text{2.76e-2}, \text{1.90e-2})^\top$ for the three components, indicating the good accuracy of the DNN approximation.
% See Fig. \ref{fig:exam:stoke-ite} for the training dynamics. The loss has some oscillations in last few iterations due to the employment of path-following strategy.

% \begin{figure}[hbt!]
% \centering\setlength{\tabcolsep}{0pt}
% \begin{tabular}{cc}
% \includegraphics[height=6cm]  {lossfunction-4.eps} &  \includegraphics[height=6cm]  {error-4.eps}\\
% (a) $L_{\sigma_k}^k$ vs $i$ & (b) $e$ vs $i$
% \end{tabular}
% \caption{\label{fig:exam:stoke-ite} The dynamics of the training process for IDRM for Example \ref{exam4:stoke}: {\rm(a)} the decay of the loss function versus the iteration index, {\rm(b)} the error $e$ versus the iteration index.}
% \end{figure}

\section{Proofs of Theorems \ref{thm:modelconv} and \ref{thm:modelerr}}\label{sec:proof}
In this section we present the technical proofs of Theorems \ref{thm:modelconv} and \ref{thm:modelerr}.
\subsection{Proof of Theorem \ref{thm:modelconv}}\label{ssec:modelconv}
To prove Theorem \ref{thm:modelconv}{\color{blue}, we} use extensively duality map and its properties.
\begin{Definition}
The duality map of $X$ is given for any $u\in X$, $u\neq 0$ by the only $j(u)\in X^{\prime}$ such that
\begin{equation*}
\xpnorm{j(u)}=\xnorm{u}\quad \mbox{and}\quad\outerp{u}{j(u)}=\xnorm{u}\xpnorm{j(u)}.
\end{equation*}
\end{Definition}
Since $X=W^{1,p}(\Omega)$ is uniformly convex and smooth, the duality map $j(u)$ is well-defined \cite[Lemma 3.2]{2013Roubiek}.
Moreover, $j(u)$ has the following properties.
\begin{Proposition}\label{prop:dual}
For any $u,v\in X$, we have
\begin{align*}
\outerp{u-v}{j(u)-j(v)}&\geq (p-1)\xnorm{u-v}^{2}.
\end{align*}
Furthermore, for any $u,v\in B_X(0,R)$, there holds
\begin{equation*}
\xpnorm{j(u)-j(v)}\leq C(p,R)\xnorm{u-v}^{p-1}.
\end{equation*}
\end{Proposition}
\begin{proof}
 The first inequality can be found in \cite[Corollary 2]{XU19911127}. We prove only the second inequality.
For any $u,v\in B_X(0,R)\subset X\equiv W^{1,p}(\Omega)$, since $\xpnorm{j(u)}=\xnorm{u}$, we deduce
\begin{equation*}
\begin{aligned}
\xpnorm{j(u)-j(v)}=&\xpnorm{j(u)-\frac{\xnorm{u}}{\xnorm{v}}j(v)+\left(\frac{\xnorm{u}}{\xnorm{v}}-1\right)j(v)}\\
\leq&\xnorm{u}\xpnorm{\frac{1}{\xnorm{u}}j(u)-\frac{1}{\xnorm{v}}j(v)}+\left|\xnorm{u}-\xnorm{v}\right|=:{\rm I}+{\rm II}.
\end{aligned}
\end{equation*}
Note that both $j(\tilde u):=\frac{1}{\xnorm{u}}j(u)$ and $j(\tilde v):=\frac{1}{\xnorm{v}}j(v)$ belong to the unit ball $\partial B(0,1)\subset X'$, and further, $X'$ is $p^{*}$-uniformly convex (with $p^{*}=\frac{p}{p-1}$) \cite[Proposition 5.6]{MishraMolinaro:2023}. Consequently,
\begin{equation*}
\inf_{\tilde{u},\tilde{v}\in\partial B(0,1)} \frac{2-\xpnorm{j(\tilde{u})+j(\tilde{v})}}{\xpnorm{j(\tilde{u})-j(\tilde{v})}^{p^{*}}}\geq C(p).
\end{equation*}
Meanwhile, since $\xnorm{\tilde u}=\xnorm{\tilde{v}}=1$, we have
\begin{align*}
{}&\xnorm{\tilde{u}-\tilde{v}}\xpnorm{j(\tilde{u})-j(\tilde{v})}
\geq \outerp{\tilde{u}-\tilde{v}}{j(\tilde{u})-j(\tilde{v})}\\
= &2-\outerp{\tilde{u}}{j(\tilde{v})}-\outerp{\tilde{v}}{j(\tilde{u})}\\
=& 4-\big(\xnorm{\tilde{u}}+\outerp{\tilde{u}}{j(\tilde{v})}\big)-\big(\xnorm{\tilde{v   }}+\outerp{\tilde{v}}{j(\tilde{u})}\big)\\
= &4-\outerp{\tilde{u}}{j(\tilde{u})+j(\tilde{v})}-\outerp{\tilde{v}}{j(\tilde{u})+j(\tilde{v})}.
\end{align*}
Moreover, using the identity $\|\tilde u\|_X=\|\tilde v\|_X=1$ again, we have
\begin{align*}
\xnorm{\tilde{u}-\tilde{v}}\xpnorm{j(\tilde{u})-j(\tilde{v})}\geq 2\left(2-\xpnorm{j(\tilde{u})+j(\tilde{v})}\right).
\end{align*}
Therefore, we derive
\begin{equation*}
{\rm I}=\xnorm{u}\xpnorm{j(\tilde{u})-j(\tilde{v})}\leq C(p)\xnorm{u}\xnorm{\tilde{u}-\tilde{v}}^{p-1}
\leq C(p,R)\xnorm{u-v}^{p-1}.
\end{equation*}
Furthermore, for $p\in(1,2]$, the triangle inequality gives
\begin{equation*}
{\rm II}\leq \xnorm{u-v}=\xnorm{u-v}^{2-p}\xnorm{u-v}^{p-1}\leq \left(\xnorm{u}+\xnorm{v}\right)^{2-p}\xnorm{{u}-{v}}^{p-1}
\leq (2R)^{2-p}\|u-v\|^{p-1}_X.
\end{equation*}
Combining these two estimates gives
the inequality $
\xpnorm{j(u)-j(v)}\leq C(p,R)\xnorm{u-v}^{p-1}$.
\end{proof}

We now show that the loss $L^{k}(u)$ has a unique minimizer in $\lambda_k^{-1}(u_k+V_{\mu})$, with the set $V_\mu$ defined in \eqref{eqn:Vmu}. A stationary point $v$ of $L^k(u)$ satisfies $DL^k(v)=0$:
\begin{equation}\label{eqn:stationary}
\mathcal{D}_{\mu}\left[\lambda_{k}\left(v-u_{k}\right)\right]+\mathcal{A}(u_{k})-f=0,
\end{equation}
where $\mathcal{D}_{\mu}\equiv D\Phi_{\mu}$ denotes the G\^ateaux derivative of the potential $\Phi_{\mu}$, given by
\begin{equation}\label{def:lmu}
\mathcal{D}_{\mu}(u) = D\Phi_{\mu}(u) = \int_{0}^{1}(\mathcal{A}(tu)+tD\mathcal{A}(tu)^*u)\,{\rm d} t + \mu j(u).
\end{equation}
It suffices to show that equation \eqref{eqn:stationary} has a unique solution in ${\lambda_k}^{-1}(u_k+V_{\mu})$. Now we show the local monotonicity of $\mathcal{D}_{\mu}(v)$ inside $V_{\mu}$.

\begin{Lemma}\label{lemma:well-D}
Under Assumption \ref{assump:bound}, $\mathcal{D}_{\mu}\in \mathcal{M}(X,X^{\prime})$ is a uniformly monotone and continuous operator in $V_{\mu}$: for any $u,v\in V_{\mu}$, there hold
\begin{align*}
\outerp{u-v}{\mathcal{D}_{\mu}(u)-\mathcal{D}_{\mu}(v)}&\geq\xnorm{u-v}^\rho,\\
\xpnorm{\mathcal{D}_{\mu}(u)-\mathcal{D}_{\mu}(v)}&\leq C(p,R,\mu)\xnorm{u-v}^{p-1}.
\end{align*}
\end{Lemma}
\begin{proof}
Define $\mathcal{H}(u)\in X^{\prime}$ by
\begin{equation*}
\mathcal{H}(u):=\int_{0}^{1}t\Sigma(tu)u\,{\rm d}t=\int_{0}^{1}t\left(D\mathcal{A}(tu)-D\mathcal{A}(tu)^*\right)u\,{\rm d}t.
\end{equation*}
Let $h(t)=t\mathcal{A}(tu)$. Then $h'(t)=\mathcal{A}(tu)+tD\mathcal{A}(tu)u$. By the fundamental theorem of calculus, we have
\begin{equation*}
\mathcal{A}(u)=h(1)-h(0)=\int_0^1h'(t)\,{\rm d}t=\int_0^1(\mathcal{A}(tu)+tD\mathcal{A}(tu)u)\,{\rm d}t.
\end{equation*}
Then it follows from the definitions of $\mathcal{H}(u)$ and $\mathcal{D}_\mu$ that
\begin{align}
&\mathcal{D}_{\mu}(u)=\int_{0}^{1}(\mathcal{A}(tu)+tD\mathcal{A}(tu)^*u)\,{\rm d} t + \mu j(u)\nonumber\\
=&\int_{0}^{1}\left(\mathcal{A}(tu)+tD\mathcal{A}(tu)u+tD\mathcal{A}(tu)^*u-tD\mathcal{A}(tu)u\right)\,{\rm d} t + \mu j(u)\nonumber\\
=&\mathcal{A}(u)-\mathcal{H}(u) + \mu j(u). \label{eqn:Dmu}
\end{align}
The skew symmetry property $\Sigma(tv)^*=-\Sigma(tv)$ implies that for any $w\in X$,
\begin{equation*}
\outerp{w}{\Sigma(tv)w}=\outerp{w}{\Sigma(tv)^*w}=-\outerp{w}{\Sigma(tv)w},
\end{equation*}
i.e., $\outerp{w}{\Sigma(tv)w}=0$. Thus, we can express $\outerp{u-v}{\mathcal{H}(u)-\mathcal{H}(v)}$ as
\begin{align*}
&\outerp{u-v}{\mathcal{H}(u)-\mathcal{H}(v)}=\int_{0}^{1}t\outerp{u-v}{\Sigma(tu)u-\Sigma(tv)v}\,{\rm d}t\\
=&\int_{0}^{1}t(\outerp{-v}{\Sigma(tu)u}-\outerp{u}{\Sigma(tv)v})\,{\rm d}t\\
=&\int_{0}^{1}t(\outerp{u}{\Sigma(tu)v}-\outerp{u}{\Sigma(tv)v})\,{\rm d}t\\
%=&\int_{0}^{1}t(\outerp{u}{\Sigma(tu)(v-u)}-\outerp{u}{\Sigma(tv)(v-u)})\,{\rm d}t\\
=&\int_{0}^{1}t\outerp{u}{[\Sigma(tu)-\Sigma(tv)](v-u)}\,{\rm d}t.
\end{align*}
Since $u,v\in V_\mu$ are bounded, Assumption \ref{eqassp}(iv) implies the continuity of $\Sigma:X\to X'$:
\begin{equation*}
\|{\Sigma(tu)-\Sigma(tv)}\|_{L(X,X')}\leq Nt\xnorm{u-v},
\end{equation*}
and consequently,
\begin{align*}
\|\left[\Sigma(tu)-\Sigma(tv)\right](v-u)\|_{X^{\prime}}\leq \|\Sigma(tu)-\Sigma(tv)\|_{L(X,X^{\prime})}\|v-u\|_X\leq Nt\|u-v\|_{X}^2.
\end{align*}
Consequently,
\begin{equation}\label{eqn:Dcts}
\outerp{u-v}{\mathcal{H}(u)-\mathcal{H}(v)}\leq N\int_{0}^{1}t^{2}\xnorm{u}\xnorm{u-v}^{2}\,{\rm d}t=\frac{N}{3}\xnorm{u}\xnorm{u-v}^{2}.
\end{equation}
Combining this estimate with Proposition \ref{prop:dual} and Assumption \ref{eqassp}(iii) gives
\begin{align}
{}\outerp{u-v}{\mathcal{D}_{\mu}(u)-\mathcal{D}_{\mu}(v)}=&\outerp{u-v}{\mathcal{A}(u)-\mathcal{A}(v)}-\outerp{u-v}{\mathcal{H}(u)-\mathcal{H}(v)} \nonumber\\&+\mu\outerp{u-v}{j(u)-j(v)}\nonumber\\
\geq& \xnorm{u-v}^{\rho}+(\mu(p-1)-\tfrac{N}{3}\xnorm{u})\xnorm{v-u}^{2}. \label{eqn:monoOfDmu}
\end{align}
Thus, $\mathcal{D}_\mu$ is monotone for $\xnorm{u}\leq\tfrac{3\mu(p-1)}{N}$.
Next, we prove the continuity of $\mathcal{D_\mu}$, for which it suffices to show the continuity of $\mathcal{H}$.
In fact, in view of the identity
\begin{align*}
&\mathcal{H}(u)-\mathcal{H}(v)=\int_{0}^{1}t[\Sigma(tu)u- \Sigma(tv)v]\,{\rm d}t=\int_{0}^{1}t[\Sigma(tu)u -\Sigma(tu)v+\Sigma(tu)v-\Sigma(tv)v]\,{\rm d}t,
\end{align*}
the boundedness of $u,v$ and the continuity of $\Sigma$, cf. Assumption \ref{eqassp}(iv), we deduce
\begin{align*}
\|\mathcal{H}(u)-\mathcal{H}(v)\|_{X'}
\leq&\int_{0}^{1}t[\|\Sigma(tu)\|_{L(X,X^{\prime})}\xnorm{u-v}+\|\Sigma(tu)-\Sigma(tv)\|_{L(X,X^{\prime})}\xnorm{v}]\,{\rm d}t\\
\leq&\int_{0}^{1}t[tRN\xnorm{u-v}+tRN\xnorm{u-v}]\,{\rm d}t\leq \frac{2}{3}RN\xnorm{u-v}.
\end{align*}
This estimate, the identity \eqref{eqn:Dmu}, the continuity of $\mathcal{A}$ and $j$ (cf. Proposition \ref{prop:dual}) and the boundedness of $u$ and $v$ yield
$$\begin{aligned}
&\xpnorm{\mathcal{D}_\mu(u)-\mathcal{D}_\mu(v)}\leq \xpnorm{\mathcal{A}(u)-\mathcal{A}(v)}+\xpnorm{\mathcal{H}(u)-\mathcal{H}(v)}+\mu\xpnorm{j(u)-j(v)}\\
\leq&M(R)\xnorm{u-v}^{p-1}+\tfrac{2}{3}RN\xnorm{u-v}+C(p,R)\mu\xnorm{u-v}^{p-1} \leq C(p,R,\mu)\xnorm{u-v}^{p-1}.
\end{aligned}$$
This completes the proof of the lemma.
\end{proof}

\begin{Proposition}\label{prop:exist}
For $\mu\geq\frac{N\xpnorm{\mathcal{A}(u_0)-f}}{3(p-1)}$, problem \eqref{eqn:stationary} has a unique solution $v_{k+1}$ in ${\lambda_k}^{-1}(u_k+V_{\mu})$.
\end{Proposition}

\begin{proof}
The proof is similar to that for globally monotone operators.
Set $g=f-\mathcal{A}(u_{k})$, and let $\{V_{m}\}_{m=1}^\infty$ be a sequence of finite-dimensional subspaces of $X$ such that $\cup_{m=1}^\infty V_m = X$ and $w^{m}\in V_m$ be its Galerkin solution,  i.e.,
\begin{equation}\label{eqn:Galerkin}
\outerp{v}{\mathcal{D}_{\mu}\left(w^{m}\right)-g}=0,\quad \forall v\in V_{m}.
\end{equation}
The existence of $w^{m}\in V_m$ follows directly from Brouwer fixed-point theorem; see \cite[Theorem 2.6]{2013Roubiek} for a detailed proof.
It follows from \eqref{eqn:Galerkin} (with $v=w^m$) and \eqref{eqn:monoOfDmu} (with $u=0$ and $v=v^m$) that
\begin{equation*}
\xpnorm{g}\geq\frac{\outerp{w^{m}}{g}}{\xnorm{w^{m}}}
=\frac{\outerp{w^{m}}{\mathcal{D}_{\mu}\left(w^{m}\right)}}{\xnorm{w^{m}}}
\geq\xnorm{w^{m}}^{\rho}+\mu(p-1)\xnorm{w^{m}}^{2}.
\end{equation*}
Thus, the sequence $\{w^{m}\}\subset V_\mu$ is uniformly bounded in $X$, and by the continuity of $\mathcal{A}$, $\{\mathcal{D}_{\mu}(w^{m})\}_{m=1}^\infty$ is uniformly bounded in $X'$. Thus, $\{w^m\}$ contains a subsequence, still denoted by $\{w^m\}$, that weakly converges to some $w^{\infty}\in V_{\mu}$ in $X$, and $\{\mathcal{D}_{\mu}(w^{m})\}$ contains a subsequence, still denoted by $\{\mathcal{D}_{\mu}(w^{m})\}$, converging weakly to some $f^{\infty} \in X'$.
Now for any $v \in \cup_{m}V_{m}$,
\begin{align*}
&{}\limsup_{m}\outerp{w^{m}}{\mathcal{D}_{\mu}(w^{m})}
=\limsup_{m}\left[\outerp{w^{m}-v}{\mathcal{D}_{\mu}(w^{m})}+\outerp{v}{\mathcal{D}_{\mu}(w^{m})}\right].
%=&\limsup_{m}\left[\outerp{w^{m}-v}{g}+\outerp{v}{\mathcal{D}_{\mu}(w^{m})}\right]\leq \outerp{w^{\infty}-v}{g}+\outerp{v}{f^{\infty}}.
\end{align*}
Since $v\in \cup_m V_m$, there exists a sequence $\{v^m\in V_m\}$ such that $v^m\to v$ in $X$. Then{\color{blue},} the Galerkin orthogonality \eqref{eqn:Galerkin} implies
\begin{align*}
&\limsup_{m}\outerp{w^{m}-v}{\mathcal{D}_{\mu}(w^{m})} = \limsup_m \left[\outerp{w^{m}-
v^m}{\mathcal{D}_{\mu}(w^{m})} + \outerp{v^m-v}{\mathcal{D}_{\mu}(w^{m})}\right]\\
\leq& \limsup_m \outerp{w^{m}-v^m}{g} + \limsup_m \|v^m-v\|_X\|\mathcal{D}_\mu(w^m)\|_{X'} = \outerp{w^\infty-v}{g}.
\end{align*}
Meanwhile, the weak convergence of $\mathcal{D}_\mu(w^m)$ to $f^\infty$ in $X'$ implies
\begin{equation*}
\limsup_{m}\outerp{v}{\mathcal{D}_{\mu}(w^{m})}=\outerp{v}{f^\infty}.
\end{equation*}
%Notice that $\outerp{w^{\infty}-v}{g}+\outerp{v}{f^{\infty}}$ is continuous with respect to $\xnorm{\cdot}$.
By the density of $\cup_{m}V_{m}$ in $X$, the inequality holds for any $v\in X$.
Taking $v=w^{\infty}$ gives
$$\limsup_{m}\outerp{w^{m}}{\mathcal{D}_{\mu}(w^{m})}\leq \outerp{w^{\infty}}{f^{\infty}}.$$
Further, the local monotonicity of $\mathcal{D}_\mu$ gives
\begin{align*}
&\liminf_{m} \outerp{w^{m}}{\mathcal{D}_{\mu}(w^{m})}\\
=&\liminf_{m} \left[\outerp{w^{m}}{\mathcal{D}_{\mu}(w^{m})}+\outerp{w^{\infty}}{\mathcal{D}_{\mu}(w^{\infty})}-\outerp{w^{m}}{\mathcal{D}_{\mu}(w^{\infty})}\right]\\
=&\liminf_{m} \left[\outerp{w^{m}-w^{\infty}}{\mathcal{D}_{\mu}(w^{m})-\mathcal{D}_{\mu}(w^{\infty})}+\outerp{w^{\infty}}{\mathcal{D}_{\mu}(w^{m})}\right]\\
\geq &\liminf_{m}\outerp{w^{\infty}}{\mathcal{D}_{\mu}(w^{m})}=\outerp{w^{\infty}}{f^{\infty}}.
\end{align*}
Then by the squeeze theorem,
\begin{align*}
\lim_{m}\outerp{w^{m}}{\mathcal{D}_{\mu}(w^{m})}= \outerp{w^{\infty}}{f^{\infty}}.
\end{align*}
It remains to show $f^\infty=\mathcal{D}_\mu(w^\infty)$. For any $w\in V_{\mu}$, let
$w^{\varepsilon}=w^{\infty}+\varepsilon w$. By the continuity of $\mathcal{D}_{\mu}$ from
Lemma \ref{lemma:well-D}, we have
\begin{align*}
&{}\outerp{w}{\mathcal{D}_{\mu}(w^{\infty})-f^{\infty}}
=\lim_{\varepsilon\rightarrow 0} \outerp{w}{\mathcal{D}_{\mu}(w^{\varepsilon})-f^{\infty}}
=\lim_{\varepsilon\rightarrow 0}\varepsilon^{-1}\outerp{w^{\varepsilon}-w^{\infty}}{\mathcal{D}_{\mu}(w^{\varepsilon})-f^{\infty}}\\
=&\lim_{\varepsilon\rightarrow 0}\varepsilon^{-1}\left[\outerp{w^{\varepsilon}}{\mathcal{D}_{\mu}(w^{\varepsilon})}+\outerp{w^{\infty}}{f^{\infty}}-\outerp{w^{\varepsilon}}{f^{\infty}}-\outerp{w^{\infty}}{\mathcal{D}_{\mu}(w^{\varepsilon})}\right]\\
=&\lim_{\varepsilon\rightarrow 0}\varepsilon^{-1}\left\{\lim_{m}\left[\outerp{w^{m}}{\mathcal{D}_{\mu}(w^{m})}
-\outerp{w^{\varepsilon}}{\mathcal{D}_{\mu}(w^{m})}
-\outerp{w^{m}}{\mathcal{D}_{\mu}(w^{\varepsilon})}\right]+\outerp{w^{\varepsilon}}{\mathcal{D}_{\mu}(w^{\varepsilon})}\right\}\\
=&\lim_{\varepsilon\rightarrow 0}\varepsilon^{-1}\lim_{m}\outerp{w^{m}-w^{\varepsilon}}{\mathcal{D}_{\mu}(w^{m})-\mathcal{D}_{\mu}(w^{\varepsilon})}\geq 0.
\end{align*}
That is, $f^{\infty} = \mathcal{D}_{\mu}(w^{\infty})$. Finally by the Galerkin relation \eqref{eqn:Galerkin} and the density of $\cup_{m}V_{m}$ in $X$, we deduce that $w^{\infty}$ is the unique solution.
\end{proof}

Now we show that the unique solution $v_{k+1}$ of problem \eqref{eqn:stationary}, i.e.,
\begin{equation}\label{eqn:vk+1sta}
\mathcal{D}_{\mu}\left[\lambda_{k}\left(v_{k+1}-u_{k}\right)\right]+\mathcal{A}(u_{k})-f=0,
\end{equation}
$v_{k+1}$ minimizes $L^k(u)$. Note that the Bregman distance $B_{\mathcal{D}_\mu}(v,u)$ of $\Phi_\mu$ at $u$ is defined by \cite{BREGMAN1967200}
\begin{equation*}
B_{\mathcal{D}_\mu}(v,u) = \Phi_\mu(v)-\Phi_\mu(u)-\langle v-u,\mathcal{D}_\mu(u)\rangle_{X,X'}.
\end{equation*}

\begin{Proposition}\label{prop:minimizer}
Under Assumption \ref{assump:bound}, $L^{k}(u)$ has a unique minimizer $v_{k+1}$ in the set $\lambda_k^{-1}(u_k+V_{\mu})$.
\end{Proposition}
\begin{proof}
It suffices to show that $L^k(u)$ is locally convex and continuous, for which we prove the local convexity and continuity of the potential $\Phi_\mu(u)$. In fact, for any $u, v\in V_{\mu}$,
$$ \Phi_{\mu}(v)-\Phi_{\mu}(u)=\int_{0}^{1}\frac{\,{\rm d}}{\,{\rm d}t}\Phi_{\mu}\left[u+t(v-u)\right]\,{\rm d}t=\int_{0}^{1}\outerp{v-u}{\mathcal{D}_{\mu}\left[u+t(v-u)\right]}\,{\rm d}t.$$
Thus the Bregman distance $B_{\mathcal{D}_\mu}(v,u)$ is represented by
\begin{align*}
B_{\mathcal{D}_\mu}(v,u)=\int_{0}^{1}t^{-1}\outerp{t(v-u)}{\mathcal{D}_{\mu}\left[u+t(v-u)\right]-\mathcal{D}_{\mu}(u)}\,{\rm d}t.
\end{align*}
By Lemma \ref{lemma:well-D}, we have
\begin{align}
%&\Phi_{\mu}(v)-\Phi_{\mu}(u)-\outerp{v-u}{\mathcal{D}_{\mu}(u)}\\
B_{\mathcal{D}_\mu}(v,u)\geq&\int_{0}^{1}t^{\rho-1}\xnorm{u-v}^{\rho}\,{\rm d}t=\rho^{-1}\xnorm{u-v}^{\rho}.\label{ineq:Phiconvex}
\end{align}
Also the following inequality is direct from Lemma \ref{lemma:well-D}:
\begin{equation}\label{eqn:Lcts}
{}\outerp{u-v}{\mathcal{D}_{\mu}(u)-\mathcal{D}_{\mu}(v)}\leq \xnorm{u-v}\xpnorm{\mathcal{D}_{\mu}(u)-\mathcal{D}_{\mu}(v)}\leq C(p,R,\mu)\xnorm{u-v}^{p}.
\end{equation}
Thus, we deduce the continuity of the potential $\Phi_\mu(u)$:
\begin{equation*}
\begin{aligned}
B_{\mathcal{D}_\mu}(v,u)
\leq&\int_{0}^{1}t^{p-1}C(p,R,\mu)\xnorm{u-v}^{p}\,{\rm d}t
=C(p,R,\mu)\xnorm{u-v}^{p}.
\end{aligned}
\end{equation*}
Then the convex optimization theory \cite{boydvandenberghe2004} gives the desired result.
\end{proof}

Next we relate the error to the dual potential $\Phi_{\mu}^{*}$. Roughly Lemma \ref{lemma:equiva} shows the equivalence of the norm $\xnorm{u-u^*}$ and the dual potential $\Phi_\mu^*$.
\begin{Lemma}\label{lemma:equiva}
For every $u\in B_X(0,R)\cap V_\mu$, there holds
\begin{equation*}
C(p,\rho,R,\mu)\xnorm{u-u^*}^{\frac{\rho(\rho-1)}{p-1}}\leq \Phi_{\mu}^{*}\left[\mathcal{A}\left(u\right)-f\right]\leq C(p,\rho,R,\mu)\xnorm{u-u^*}^{\frac{p(p-1)}{\rho-1}}.
\end{equation*}
\end{Lemma}
\begin{proof}
Let $\mathcal{F}=\{f\in X^{\prime}: \exists u\in V_\mu \text{ such that } \mathcal{D}_\mu(u)=f\}$. Proposition \ref{prop:exist} implies that $\outerp{u}{f}-\Phi_{\mu}(u)$ has a unique maximizer $\tilde{u}=\mathcal{D}_{\mu}^{-1}(f)\in V_\mu$ for any $f\in\mathcal{F}$. Thus for every $u\in V_\mu$, upon letting $g=\mathcal{A}(u)-f (=\mathcal{A}(u)-\mathcal{A}(u^*))$ and noting $\Phi_\mu(0)=0$, we have
\begin{align*}
\Phi_{\mu}^{*}[g]&=\outerp{\mathcal{D}_{\mu}^{-1}[g]}{g}-\Phi_{\mu}(\mathcal{D}_{\mu}^{-1}[g])\\
&=\outerp{\mathcal{D}_{\mu}^{-1}[g]}{g}+(\Phi_{\mu}({0})-\Phi_{\mu}(\mathcal{D}_{\mu}^{-1}[g])).
\end{align*}
The local convexity and continuity of $\Phi_\mu(u)$ in Proposition \ref{prop:minimizer} imply
\begin{align*}
\Phi_{\mu}({0})-\Phi_{\mu}(\mathcal{D}_{\mu}^{-1}[g]) &\geq-\outerp{\mathcal{D}_{\mu}^{-1}[g]}{g}+\rho^{-1}\xnorm{\mathcal{D}_{\mu}^{-1}[g]}^\rho,\\
\Phi_{\mu}({0})-\Phi_{\mu}(\mathcal{D}_{\mu}^{-1}[g])&\leq-\outerp{\mathcal{D}_{\mu}^{-1}[g]}{g}+C(p,R,\mu)\xnorm{\mathcal{D}_{\mu}^{-1}[g]}^{p}.
\end{align*}
Combining these two inequalities leads to
\begin{equation}\label{ineq:equivalence}
\rho^{-1}\xnorm{\mathcal{D}_{\mu}^{-1}[g]}^\rho\leq \Phi_{\mu}^{*}[g]\leq C(p,R,\mu)\xnorm{\mathcal{D}_{\mu}^{-1}[g]}^{p}.
\end{equation}
Upon noting $\mathcal{D}_\mu^{-1}(0)=0$, by the second estimate in Lemma \ref{lemma:well-D} and \eqref{eqn:Acoer}, we deduce
\begin{align}
\xnorm{\mathcal{D}_{\mu}^{-1}[g]}&=\xnorm{\mathcal{D}_{\mu}^{-1}[g]-\mathcal{D}_\mu^{-1}({0})}\geq C(p,R,\mu)\xpnorm{g}^{\frac{1}{p-1}}\geq C(p,R,\mu)\xnorm{u-u^*}^{\frac{\rho-1}{p-1}}. \label{ineq:wellpose-D-1}
\end{align}
Similarly, by the first estimate in Lemma \ref{lemma:well-D} and the continuity of $\mathcal{A}$ in Assumption \ref{eqassp}(i), we obtain
\begin{align}
\xnorm{\mathcal{D}_{\mu}^{-1}[g]}&=\xnorm{\mathcal{D}_{\mu}^{-1}[g]-\mathcal{D}_\mu^{-1}({0})}\leq \xpnorm{g}^{\frac{1}{\rho-1}}\leq M(R)^{\frac{1}{\rho-1}}\xnorm{u-u^*}^{\frac{p-1}{\rho-1}}.\label{ineq:wellpose-D-2}
\end{align}
Combining this with \eqref{ineq:equivalence} gives the desired assertion.
\end{proof}

Now we can state the proof of Theorem \ref{thm:modelconv}.
\begin{proof}
Note that $v_{k+1}$ minimizes $L^k(u)$. Let $r_{k+1}=f-\mathcal{A}\left(v_{k+1}\right)$ and $p_k=f-\mathcal{A}\left(u_k\right)$. The starting point of the analysis is the following inequality $
\xnorm{u_{k+1}-u^*}\leq\xnorm{u_{k+1}-v_{k+1}}+\xnorm{v_{k+1}-u^*}$. The rest of the lengthy proof is divided into five steps.

\medskip\noindent \textit{Step 1}. Bound $\|u_{k+1}-v_{k+1}\|$.
It follows from Algorithm \ref{alg:alg1} that ${L}^{k}(u_{k+1})\leq {L}^{k}(v_{k+1})+\epsilon_M$. This inequality, the expression of $L^k(u)$ and the identity \eqref{eqn:vk+1sta} (i.e., $\mathcal{D}_\mu(\lambda_k(v_{k+1}-u_k)) =p_k$) yield
\begin{equation*}
%\Phi_\mu[\lambda_k(u_{k+1}-u_k)]-\Phi_\mu[\lambda_k(v_{k+1}-u_k)]+\outerp{\lambda_{k}(u_{k+1}-v_{k+1})}{p_k}
B_{\mathcal{D}_\mu}(\lambda_k(u_{k+1}-u_k),\lambda_k(v_{k+1}-u_k)) \leq \epsilon_M.
\end{equation*}
By the inequality \eqref{ineq:Phiconvex}, we have
\begin{equation*}
%\Phi_\mu[\lambda_k(u_{k+1}-u_k)]-\Phi_\mu[\lambda_k(v_{k+1}-u_k)]-\outerp{\lambda_k(u_{k+1}-v_{k+1})}{\mathcal{D}_{\mu}(\lambda_k(v_{k+1}-u_k))}
B_{\mathcal{D}_\mu}(\lambda_k(u_{k+1}-u_k),\lambda_k(v_{k+1}-u_k))
\geq \lambda_k^\rho\rho^{-1}\xnorm{u_{k+1}-v_{k+1}}^\rho.
\end{equation*}
Consequently, $
\epsilon_M\geq
\rho^{-1} \lambda_k^\rho\xnorm{u_{k+1}-v_{k+1}}^\rho$, i.e.,
\begin{equation}\label{ineq:uv}
\xnorm{u_{k+1}-v_{k+1}}\leq \lambda_k^{-1}\rho^{\frac{1}{\rho}}\epsilon_M^{\frac{1}{\rho}}.
\end{equation}

\noindent\textit{Step 2.}
Bound $\Phi_{\mu}^{*}(r_{k+1})-\Phi_\mu^*(p_k)$. We appeal to Lemma \ref{lemma:equiva}.
The definition of $\Phi^*_\mu$ gives
$$\begin{aligned}
&\Phi_{\mu}^{*}(r_{k+1})-\Phi_\mu^*(p_k)\\
=&\outerp{\mathcal{D}_{\mu}^{-1}(r_{k+1})}{r_{k+1}}-\Phi_{\mu}\left(\mathcal{D}_{\mu}^{-1}(r_{k+1})\right)-(\outerp{\mathcal{D}_{\mu}^{-1}(p_k)}{p_k}-\Phi_{\mu}\left(\mathcal{D}_{\mu}^{-1}(p_k)\right))\\
=&\outerp{\mathcal{D}_{\mu}^{-1}(r_{k+1})-\mathcal{D}_{\mu}^{-1}(p_k)}{r_{k+1}}+\outerp{\mathcal{D}_{\mu}^{-1}(p_k)}{r_{k+1}-p_k}-\left(\Phi_{\mu}\left(\mathcal{D}_{\mu}^{-1}(r_{k+1})\right)-\Phi_{\mu}\left(\mathcal{D}_{\mu}^{-1}(p_k)\right)\right).
\end{aligned}
$$
The local convexity of $\Phi_\mu$ in \eqref{ineq:Phiconvex} implies
\begin{align*}
&\Phi_{\mu}\left(\mathcal{D}_{\mu}^{-1}(r_{k+1})\right)-\Phi_{\mu}\left(\mathcal{D}_{\mu}^{-1}(p_k)\right)\\
\geq &\outerp{\mathcal{D}_{\mu}^{-1}(r_{k+1})-\mathcal{D}_{\mu}^{-1}(p_k)}{p_k}+\rho^{-1}\xnorm{\mathcal{D}_{\mu}^{-1}(r_{k+1})-\mathcal{D}_{\mu}^{-1}(p_k)}^\rho\\
\geq &\outerp{\mathcal{D}_{\mu}^{-1}(r_{k+1})-\mathcal{D}_{\mu}^{-1}(p_k)}{p_k}.
\end{align*}
Consequently,
\begin{align*}%\label{ineq:error1+2}
\Phi_{\mu}^{*}(r_{k+1})-\Phi_\mu^*(p_k)\leq& \outerp{\mathcal{D}_{\mu}^{-1}(r_{k+1})-\mathcal{D}_{\mu}^{-1}(p_k)}{r_{k+1}-p_k}
+\outerp{\mathcal{D}_{\mu}^{-1}(p_k)}{r_{k+1}-p_k}=: {\rm I}+{\rm II}.
\end{align*}
By Lemma $\ref{lemma:well-D}$ and the continuity of $\mathcal{A}$ from Assumption \ref{eqassp}(i), there holds with $C=C(p,\rho,R,\mu)$,
\begin{align}
{\rm I}&\leq C\xnorm{\mathcal{D}_{\mu}^{-1}(r_{k+1})-\mathcal{D}_{\mu}^{-1}(p_k)}^{p}
\leq C\xpnorm{r_{k+1}-p_k}^{\frac{p}{\rho-1}}\nonumber\\
&= C\xpnorm{\mathcal{A}(u_k)-\mathcal{A}(v_{k+1})}^{\frac{p}{\rho-1}} \leq
C\xnorm{u_k-v_{k+1}}^{\frac{p(p-1)}{\rho-1}}.\label{ineq:error1con}
\end{align}
Then the identity $\mathcal{D}_{\mu}^{-1}(p_k)=\lambda_k(v_{k+1}-u_k)$ and the monotonicity of $\mathcal{A}$ in Assumption \ref{eqassp}(iii) imply
$$ \begin{aligned}
{\rm II}=&\lambda_k\outerp{v_{k+1}-u_k}{r_{k+1}-p_k}\\
=&-\lambda_k\outerp{v_{k+1}-u_k}{\mathcal{A}(v_{k+1})-\mathcal{A}(u_k)}\leq -\lambda_k\xnorm{v_{k+1}-u_k}^{\rho}.
\end{aligned}
$$
Combining the last two estimates and using \eqref{ineq:equivalence} leads to
\begin{align*}
\Phi_{\mu}^{*}(r_{k+1})-\Phi_\mu^*(p_k)&\leq  C(p,\rho,R,\mu)\xnorm{v_{k+1}-u_k}^{\frac{p(p-1)}{\rho-1}}-\lambda_k\xnorm{v_{k+1}-u_k}^{\rho}\\
&=C(p,\rho,R,\mu)\lambda_k^{-\frac{p(p-1)}{\rho-1}}\xnorm{\mathcal{D}_{\mu}^{-1}(p_k)}^{\frac{p(p-1)}{\rho-1}}-\lambda_k^{-(\rho-1)}\xnorm{\mathcal{D}_{\mu}^{-1}(p_k)}^{\rho}\\
&\leq C(p,\rho,R,\mu)\lambda_k^{-\frac{p(p-1)}{\rho-1}}\Phi_\mu^*(p_k)^{\frac{p(p-1)}{\rho(\rho-1)}}-C(p,R,\mu)\lambda_k^{-(\rho-1)}\Phi_\mu^*(p_k)^{\frac{\rho}{p}}.
\end{align*}
Now by the choice $\lambda_k=C(p,\rho.R,\mu)(-L^k(v_{k+1}))^{\alpha}=C(p,\rho,R,\mu)\Phi_\mu^*(p_k)^\alpha:=c_{\rm s}\Phi_\mu^*(p_k)^\alpha$ (cf. \eqref{eqn:lam=Phi*}), we obtain for some $c_{\rm p}=C(p,\rho,R,\mu)>0$
\begin{equation}\label{ineq:phistar1}
\begin{aligned}
\Phi_{\mu}^{*}(r_{k+1})-\Phi_\mu^*(p_k)\leq -c_{\rm p}(\Phi_\mu^*(p_k))^{\beta}.
\end{aligned}
\end{equation}

\medskip\noindent\textit{Step 3}. Bound $|\Phi_{\mu}^{*}(r_{k+1})-\Phi_{\mu}^{*}(p_{k+1})|$.
We may repeat the argument by replacing $p_k$ with $p_{k+1}$ till the step \eqref{ineq:error1con} and get
$$\begin{aligned}
&|\Phi_{\mu}^{*}(r_{k+1})-\Phi_{\mu}^{*}(p_{k+1})|\\
\leq &\left|\outerp{\mathcal{D}_{\mu}^{-1}(r_{k+1})-\mathcal{D}_{\mu}^{-1}(p_{k+1})}{r_{k+1}-p_{k+1}}\right|+\left|\outerp{\mathcal{D}_{\mu}^{-1}(p_{k+1})}{r_{k+1}-p_{k+1}}\right|\\
\leq &C(p,\rho,R,\mu)\xnorm{v_{k+1}-u_{k+1}}^{\frac{p(p-1)}{\rho-1}}+\xnorm{\mathcal{D}_{\mu}^{-1}(p_{k+1})}M(R)\xnorm{v_{k+1}-u_{k+1}}.
\end{aligned}$$
By the estimate \eqref{ineq:wellpose-D-2} and Assumption \ref{assump:bound}, we have
\begin{equation*}
\xnorm{\mathcal{D}_\mu^{-1}(p_{k+1})}=\xnorm{\mathcal{D}_\mu^{-1}(f-\mathcal{A}(u_{k+1}))}\leq C(\rho,R)||u^*-u_{k+1}||_X^{\frac{p-1}{\rho-1}}\leq C(\rho,R).
\end{equation*}
The last two estimates and \eqref{ineq:uv} imply
\begin{align}
|\Phi_{\mu}^{*}(r_{k+1})-\Phi_{\mu}^{*}(p_{k+1})|\leq &C(p,\rho,R,\mu)\xnorm{v_{k+1}-u_{k+1}}^{\frac{p(p-1)}{\rho-1}}+C(\rho,R)\xnorm{v_{k+1}-u_{k+1}}\nonumber\\
\leq & C(p,\rho,R,\mu)\big(\lambda_k^{-\frac{p(p-1)}{\rho-1}}\epsilon_M^\frac{p(p-1)}{\rho(\rho-1)}+\lambda_k^{-1}\epsilon_M^{\frac{1}{\rho}}\big).
\label{ineq:phistar2}
\end{align}
We denote the constant by $c_{\rm r}$ and the upper bound on $\lambda_k^{-\frac{p(p-1)}{\rho-1}}\epsilon_M^\frac{p(p-1)}{\rho(\rho-1)}+\lambda_k^{-1}\epsilon_M^{\frac{1}{\rho}}$ by $\varepsilon$. Then the inequality $\Phi_\mu^*(p_{k+1})-\Phi_\mu^*(p_k) \leq |\Phi_\mu^*(p_{k+1})-\Phi_\mu^*(r_{k+1})|+ \Phi_\mu^*(r_{k+1})-\Phi_\mu^*(p_k)$ and the estimates \eqref{ineq:phistar1} and \eqref{ineq:phistar2} imply that the sequence $\{\Phi_\mu^*(p_k)\}_{k=1}^\infty$ satisfies
\begin{equation}\label{eqn:decay-Phi}
\Phi_{\mu}^{*}(p_{k+1})\leq c_{\rm r}\varepsilon+\Phi_{\mu}^{*}(p_{k})-c_{\rm p}\Phi_\mu^*(p_k)^{\beta}.
\end{equation}
This recursion forms the basis for proving assertions (i) and (ii) below. Note that the conditions $\rho\geq2$ and $p\in(1,2]$ imply
$\rho^2-p(\rho-1)\geq\rho^2-2(\rho-1)=(\rho-1)^2+1\geq2.$ This and the condition $0<p(p-1)-(\rho-1)^2\leq 1$ also ensure $\beta>0$.

\medskip\noindent \textit{Step 4}. Prove the monotone decay of $\Phi_\mu^*(p_k)$.
 The recursion \eqref{eqn:decay-Phi} indicates that $\Phi_\mu^*(p_k)$ is decreasing whenever
$c_{\rm r}\varepsilon<c_{\rm p}(\Phi_\mu^*(p_k))^{\beta}$.  The condition $\Phi_\mu^*(p_k)=-L^k(v_{k+1})>(2c_{\rm t})^{\frac{\rho(\rho-1)}{p(p-1)}}\epsilon_M$ (from the definition of $k^*$) and the identity $\beta+\frac{\alpha p(p-1)}{\rho-1}=\frac{p(p-1)}{\rho(\rho-1)}$ imply
\begin{equation*}
\Phi_\mu^*(p_k)^{-\frac{\alpha p(p-1)}{\rho-1}}\epsilon_M^{\frac{p(p-1)}{\rho(\rho-1)}}< (2c_{\rm t})^{-1} \Phi_\mu^*(p_k)^{\beta},
\end{equation*}
whereas the condition $\Phi_\mu^*(p_k)=-L^k(v_{k+1})>(2c_{\rm t}\epsilon_M^{\frac{1}{\rho}})^{\frac{1}{\alpha+\beta}}$ gives
\begin{equation*}
\Phi_\mu^*(p_k)^{-\alpha}\epsilon_M^{\frac{1}{\rho}} < (2c_{\rm t})^{-1} \Phi_\mu^*(p_k)^{\beta}.
\end{equation*}
Thus the following inequality holds:
$\Phi_\mu^*(p_k)^{-\frac{\alpha p(p-1)}{\rho-1}}\epsilon_M^{\frac{p(p-1)}{\rho(\rho-1)}}+\Phi_\mu^*(p_k)^{-\alpha}\epsilon_M^{\frac{1}{\rho}}< c_{\rm t}^{-1}\Phi^*(p_k)^{\beta}$.
With the choice $\lambda_k=c_{\rm s}\Phi_\mu^*(p_k)^\alpha$, it is equivalent to
\begin{equation*}
\lambda_k^{-\frac{p}{\rho-1}}\epsilon_M^\frac{p}{\rho(\rho-1)}+\lambda_k^{-1}\epsilon_M^{\frac{1}{\rho}}< c_{\rm p}c_{\rm r}^{-1}\Phi_\mu^*(p_k)^{\beta}.
\end{equation*}
This directly implies the preasymptotic monotone decay of $\Phi^*(p_k)$ in (i).

\medskip\noindent \textit{Step 5}. Bound the convergence rate. If $-L^k(v_{k+1})\geq (2c_{\rm t}\eta^{-1})^{\frac{\rho(\rho-1)}{p(p-1)}}\epsilon_M$ for some $\eta\in(0,1)$ (independent of $k$). In view of the identity \eqref{eqn:lam=Phi*}, the condition on $L^k(v_{k+1})$ implies
\begin{equation*}
\Phi_\mu^*(p_k)^{-\beta-\frac{\alpha p(p-1)}{\rho-1}}\epsilon_M^{\frac{p(p-1)}{\rho(\rho-1)}}\leq (2c_{\rm t})^{-1}\eta\Phi_\mu^*(p_k)^{\beta},
\end{equation*}
whereas the condition $\Phi_\mu^*(p_k)=-L^k(v_{k+1})>(2c_{\rm t}\eta^{-1}\epsilon_M^{\frac{1}{\rho}})^{\frac{1}{\alpha+\beta}}$ yields
\begin{equation*}
\Phi_\mu^*(p_k)^{-\alpha}\epsilon_M^{\frac{1}{\rho}} < (2c_{\rm t})^{-1}\eta \Phi_\mu^*(p_k)^{\beta}.
\end{equation*}
Repeating the preceding argument yields $c_{\rm r}\varepsilon\leq \eta c_{\rm p}(\Phi_\mu^*(p_k))^{\beta}$, and then
\begin{equation*}
\Phi_{\mu}^{*}(p_{k+1})\leq \Phi_{\mu}^{*}(p_{k})-(1-\eta)c_{\rm p}\Phi_\mu^*(p_k)^{\beta}.
\end{equation*}
This recursion allows deriving the convergence rate, and we treat the two cases $p=2$ and $p\in(1,2)$ separately. Note that
since $p\leq 2$ and $\rho\geq2$, we have $$\beta=\frac{(p-1)(\rho^2-p(\rho-1))}{\rho(p(p-1)-(\rho-1)^2)}=\frac{\rho^2-p(\rho-1)}{\rho(p-\frac{(\rho-1)^2}{p-1})}\geq \frac{\rho^2-p(\rho-1)}{\rho(2-(\rho-1)^2)}\geq \frac{\rho^2-p(\rho-1)}{\rho}\geq \frac{\rho^2-2(\rho-1)}{\rho}\geq 1,$$
and $\beta=1$ if and only if $p=\rho=2$.\\
Case (ii.a): $p=2$. Then we have
$\Phi_{\mu}^{*}(p_{k+1})\leq \eta c_{\rm p}\Phi_{\mu}^{*}(p_{k})$.
Applying the estimate repeatedly yields
\begin{equation*}
\Phi_{\mu}^{*}(p_{k})\leq (\eta c_{\rm p})^k\Phi_{\mu}^{*}(p_{0})
\end{equation*}
This and Lemma \ref{lemma:equiva} give the result of part (ii.a).\\
Case (ii.b): $p\in(1,2)$. When $\beta>1$, we have
\begin{align*}
\Phi_{\mu}^{*}(p_{k+1})\leq \Phi_{\mu}^{*}(p_{k})-(1-\eta)c_{\rm p}(\Phi_\mu^*(p_k))^{\beta}.
\end{align*}
Letting $a_k=\Phi_{\mu}^{*}(p_{k})$ gives $ a_{k+1}\leq a_k-(1-\eta)c_{\rm p}a_k^{\beta}.$
Then for the sequence $\{a_k\}_{k=1}^\infty$, we obtain
\begin{equation*}
a_{k+1}^{1-\beta}\geq a_k^{1-\beta}(1-(1-\eta)c_{\rm p}a_k^{\beta-1})^{1-\beta}.
\end{equation*}
This and the inequality $(1-x)^{1-\beta}\geq 1+(\beta-1)x$ for $x\in(0,1)$ imply
\begin{equation*}
a_{k+1}^{1-\beta}\geq a_k^{1-\beta}+(1-\eta)c_{\rm p}(\beta-1).
\end{equation*}
So it follows that
$a_k\leq ((1-\eta)c_{\rm p}(\beta-1)k+a_0^{1-\beta})^{\frac{1}{1-\beta}}$, and hence,
\begin{equation*}
\Phi_{\mu}^{*}(p_{k})\leq \left(\frac{\Phi_{\mu}^{*}(p_{0})^{\beta-1}}{(1-\eta)C(p,\rho,R,\mu)k\Phi_{\mu}^{*}(p_{0})^{\beta-1}+1}\right)^{\frac{1}{\beta-1}}\leq C(p,\rho,R,\mu,\xnorm{u_0-u^*})((1-\eta)k)^{-\frac{1}{\beta-1}}.
\end{equation*}
Combining this estimate with Lemma \ref{lemma:equiva} leads to the desired result.
\end{proof}

\begin{Corollary}\label{cor:plarge}
If $\mathcal{A}=\widetilde{\mathcal{A}}+\overline{\mathcal{A}}$, with $\widetilde{\mathcal{A}}=D\Phi$ for some convex $\Phi:X\rightarrow \mathbb{R}^{+}$ and linear $\overline{\mathcal{A}}$, then the conclusion of Theorem \ref{thm:modelconv} holds for any $p\in(1,+\infty)$ with $\mu =0$.
\end{Corollary}
\begin{proof}
Under the given conditions, $\Sigma=D\mathcal{A}-D\mathcal{A}^*\in \mathcal{M}(X,X^{\prime})$ is skew linear. Thus,
$\outerp{u-v}{\mathcal{H}(u)-\mathcal{H}(v)}=\int_{0}^{1}t\outerp{u-v}{\Sigma(tu)u-\Sigma(tv)v}\,{\rm d}t=0$,
and $\mathcal{D}_{0}$ is a globally monotone operator:
\begin{align*}
{}&\outerp{u-v}{\mathcal{D}_{0}(u)-\mathcal{D}_{0}(v)}
=\outerp{u-v}{\mathcal{A}(u)-\mathcal{A}(v)}-\outerp{u-v}{\mathcal{H}(u)-\mathcal{H}(v)}
\geq \xnorm{u-v}^{\rho}.
\end{align*}
Therefore, Lemma \ref{lemma:well-D} holds.
The rest of the proof is identical with that for $p\leq 2$.
\end{proof}

\subsection{Proof of Theorem \ref{thm:modelerr}}\label{ssec:learningerror}

We first decompose the learning error $\epsilon_M$ into the approximation error $\mathcal{E}_{\rm app}$ and the statistical error $\mathcal{E}_{\rm stat}$.
\begin{Lemma}\label{lem:decomp}
For the learning error $\epsilon_M$, we have
\begin{equation*}
\mathbb{E}_{\mathbb{X},\mathbb{T}}\left[\epsilon_M\right]\le C(R)\lambda_{k}\underbrace{ \inf_{u_{\theta} \in \mathcal{N}}\xnorm{u_{\theta}-v_{k+1}}}_{\mathcal{E}_{\rm app}}+2 \underbrace{\mathbb{E}_{\mathbb{X},\mathbb{T}}\sup _{u_{\theta} \in \mathcal{N}}\left|L^{k}(u_{\theta})-\widehat{L}^{k}(u_{\theta})\right|}_{\mathcal{E}_{\rm stat}}.
\end{equation*}
\end{Lemma}
\begin{proof}
For any $w_{\theta}\in \mathcal{N}$, by the minimizing property of $\widehat{u}_{k+1}$ to the loss $\widehat{L}^k(u_\theta)$ over the set $\mathcal{N}$, we have
\begin{align*}
\epsilon_M:=& L^{k}(\widehat{u}_{k+1})-L^{k}(v_{k+1})
=\big[L^{k}(\widehat{u}_{k+1})-\widehat{L}^{k}(\widehat{u}_{k+1})\big]+\big[\widehat{L}^{k} (\widehat{u}_{k+1})-\widehat{L}^{k}(w_{\theta})\big]\\
&+\big[\widehat{L}^{k}(w_{\theta})-L^{k}(w_{\theta})\big]
+\big[L^{k}(w_{\theta})-L^{k}(v_{k+1})\big]\\
\leq &\big[L^{k}(\widehat{u}_{k+1})-\widehat{L}^{k}(\widehat{u}_{k+1})\big]+\big[\widehat{L}^{k} (w_{\theta})-L^{k}(w_{\theta})\big]+\big[L^{k}(w_{\theta})-L^{k}(v_{k+1})\big]\\
\leq &2\sup _{u_{\theta} \in \mathcal{N}}\big|L^{k}(u_{\theta})-\widehat{L}^{k}(u_{\theta})\big|+\big[L^{k}(w_{\theta})-L^{k}(v_{k+1})\big].
\end{align*}
Then taking infimum over $w_\theta\in\mathcal{N}\cap u_k+V_\mu$, taking the expectation with respect to $\mathbb{X}$ and $\mathbb{T}$, and noting the continuity of the loss $L^k(u)$ lead to
\begin{equation*}
\mathbb{E}_{\mathbb{X},\mathbb{T}}[\epsilon_M]\leq 2\mathbb{E}_{\mathbb{X},\mathbb{T}}\sup _{u_{\theta} \in \mathcal{N}}\left|L^{k}(u_{\theta})-\widehat{L}^{k}(u_{\theta})\right|+C(R)\lambda_{k}\inf_{w_\theta\in\mathcal{N}\cap u_{k}+V_{\mu}}\xnorm{w_{\theta}-v_{k+1}}.
\end{equation*}
Since $v_{k+1}\in u_{k}+V_{\mu}$, it can be bounded by $\mathcal{E}_{\text{app}}$.
\end{proof}

For the approximation error $\mathcal{E}_{\rm app}$, we use the approximation capability of tanh DNNs in Sobolev spaces  \cite[Proposition 4.8]{GuhringRaslan2021}.
\begin{Proposition}\label{prop:approx}
For any $v\in W^{1+r,p}(\Omega)$ with $r>1$, there exists a DNN $w_{\theta}\in \mathcal{N}(L,{W}, {W}^{3(1+1/p)+(2r+2)/d})$ with $L = C\log(d + 1 + r)$ and $\nu>0$ is arbitrarily small such that
\begin{equation*}
\|v-w_{\theta}\|_{W^{1, p}(\Omega)}\leq C(d, r, p)W^{-(r - \nu)/d}\|v\|_{W^{1+r, p}(\Omega)}.
\end{equation*}
\end{Proposition}

To analyze the statistical error $\mathcal{E}_{\text{stat}}$, we use Rademacher complexity
\cite{Anthony_Bartlett_1999,Bartlett2003RademacherAG}, which measures the complexity of a
function class by its correlation with Rademacher random variables.
\begin{Definition}
The Rademacher complexity of a function class $\mathcal{F}$ defined on the domain $D$
associate with random samples $\mathbb{X}:=\{X_i\}_{i=1}^{N}$ on $D$ is defined by
\begin{equation*}
\mathfrak{R}_N(\mathcal{F}) = \mathbb{E}_{\mathbb{X},\boldsymbol{\sigma}}\left[\sup_{u\in \mathcal{F}}\frac{1}{N}\left|\sum_{i=1}^N \sigma_i u(X_i)\right|\right].
\end{equation*}
where $\boldsymbol{\sigma}:=\{\sigma_i\}_{i=1}^N$ are $N$ i.i.d.  Rademacher random variables with $\mathbb{P}(\sigma_i = 1) = \mathbb{P}(\sigma_i = -1) = \frac{1}{2}.$
\end{Definition}

We have the following bound on the statistical error $\mathcal{E}_{\rm stat}$.
\begin{Proposition}\label{prop:sta}
Let the domain $\Omega\subset (-1,1)^d$. Then under Assumption \ref{assump:boundaf}, the statistical error $\mathcal{E}_{\rm stat}$ is bounded by
\begin{equation*}
\mathbb{E}\sup_{u_{\theta} \in \mathcal{N}}\left|L^{k}(u_{\theta})-\hat{L}^{k}(u_{\theta})\right|
\leq C(|\Omega|,\ell_{a},d)\lambda_{k}LW^{\max\{3L-2,4\}}B^{3L}\sqrt{\frac{\log \left(d LWB N\right)}{N}}.
\end{equation*}
\end{Proposition}
\begin{proof} Let
$\nabla\mathcal{F}=\{\partial_{x_{i}}f: f\in\mathcal{F}, i=1,\ldots,d\}.$
By the triangle inequality, we have
\begin{align*}
&\mathbb{E}_{\mathbb{X},\mathbb{T}}\sup_{u_{\theta} \in \mathcal{N}}|L^{k}(u_{\theta})-\widehat{L}^{k}(u_{\theta})|\\
\leq & \lambda_k\mathbb{E}_{\mathbb{X},\mathbb{T}}\sup_{u \in \mathcal{N}}\left|\int_{\Omega\times(0,1)}I_{1}^{k}((x, t), u(x))\,{\rm d}(x,t)-\frac{|\Omega|}{N}\sum_{i=1}^{N}I_{1}^{k}((X_{i}, T_{i}), u(X_i))\right|\\
&+\mu\lambda_k^2\mathbb{E}_{\mathbb{X}}\sup_{u \in \mathcal{N}}\left|\left[\int_{\Omega}I_2^k(x,u(x))\,{\rm d}x\right]^{2/p}-\left[\frac{|\Omega|}{N}\sum_{i=1}^{N}I_2^k(X_i,u(X_i))\right]^{2/p}\right|\\
&+\lambda_k\mathbb{E}_{\mathbb{X}}\sup_{u \in \mathcal{N}}\left|\int_{\Omega}I_{3}^{k}(x, u(x))\,{\rm d}x-\frac{|\Omega|}{N}\sum_{i=1}^{N}I_{3}^{k}(X_{i}, u(X_i))\right|.
\end{align*}
Next, we bound the three components, denoted by ${\rm I}$, ${\rm II}$ and ${\rm III}$. Let $\widetilde{\mathbb{X}}=\{\widetilde{X}_{i}\}_{i=1}^N$ and $\widetilde{\mathbb{T}}=\{\widetilde{T}_{i}\}_{i=1}^N$ be independent copies of $\mathbb{X}$ and $\mathbb{T}$. For the term ${\rm I}$, using the convexity of $\sup$ and Jensen's inequality,
\begin{align*}
{\rm I}
=&\frac{|\Omega|}{N}\mathbb{E}_{\mathbb{X},\mathbb{T}}\sup_{u \in \mathcal{N}}\left|\mathbb{E}_{\widetilde{\mathbb{X}},\widetilde{\mathbb{T}}}\sum_{i=1}^{N}I_{1}^{k}((\widetilde{X}_{i}, \widetilde{T}_{i}), u(\widetilde{X}_i))-\sum_{i=1}^{N}I_{1}^{k}((X_{i}, T_{i}), u(X_i))\right|\\
\leq&\frac{|\Omega|}{N}\mathbb{E}_{\mathbb{X},\mathbb{T},\widetilde{\mathbb{X}},\widetilde{\mathbb{T}}}\sup_{u \in \mathcal{N}}\left|\sum_{i=1}^{N}I_{1}^{k}((\widetilde{X}_{i}, \widetilde{T}_{i}), u(\widetilde X_i))-\sum_{i=1}^{N}I_{1}^{k}((X_{i}, T_{i}), u(X_i))\right|.
\end{align*}
The independence between $\mathbb{X,T}$ and $\mathbb{\widetilde{X},\widetilde{T}}$ implies
\begin{align*}
{\rm I}\leq &\frac{|\Omega|}{N}\mathbb{E}_{\mathbb{X},\mathbb{T},\widetilde{\mathbb{X}},\widetilde{\mathbb{T}},\boldsymbol{\sigma}}\sup_{u \in \mathcal{N}}\left|\sum_{i=1}^{N}\sigma_{i}I_{1}^{k}((\widetilde{X}_{i}, \widetilde{T}_{i}), u(\widetilde{X}_i))-\sum_{i=1}^{N}\sigma_{i}I_{1}^{k}((X_{i}, T_{i}), u(X_i))\right|\\
\leq &\frac{2|\Omega|}{N}\mathbb{E}_{\mathbb{X},\mathbb{T},\boldsymbol{\sigma}}\sup_{u \in \mathcal{N}}\left|\sum_{i=1}^{N}\sigma_{i}I_{1}^{k}((X_{i}, T_{i}), u(X_i))\right|.
\end{align*}
By substituting the explicit form of $I_1^k$ and the triangle inequality, we get
\begin{align*}
{\rm I}\leq&\frac{2|\Omega|}{N}\mathbb{E}_{\mathbb{X},\mathbb{T},\boldsymbol{\sigma}}\sup_{u \in \mathcal{N}}\left|\sum_{i=1}^{N}\sigma_{i}\sum_{m=0}^{1}\nabla^{m}(u-u^{k})(X_{i})\cdot a_{m}(X_{i},T_{i}\lambda_{k}(u-u_{k})(X_i),T_{i}\lambda_{k}\nabla (u-u_{k})(X_i))\right|\\
=&\frac{2|\Omega|}{N}\mathbb{E}_{\mathbb{X},\mathbb{T},\boldsymbol{\sigma}}\sup_{u \in \mathcal{N}}\left|\sum_{m=0}^{1}\sum_{i=1}^{N}\sigma_{i}\nabla^{m}(u-u^{k})(X_{i})\cdot a_{m}(X_{i},T_{i}\lambda_{k}(u-u_{k})(X_i),T_{i}\lambda_{k}\nabla (u-u_{k})(X_i))\right|\\
\leq &\frac{2|\Omega|}{N}\sum_{m=0}^{1}\mathbb{E}_{\mathbb{X},\mathbb{T},\boldsymbol{\sigma}}\sup_{u \in \mathcal{N}}\left|\sum_{i=1}^{N}\sigma_{i}\nabla^{m}(u-u^{k})(X_{i})\right|\sup_{1\leq i\leq N}\left\|a_{m}(x,T_{i}\lambda_{k}(u-u_{k}),T_{i}\lambda_{k}\nabla (u-u_{k})\right\|_{L^{\infty}(\Omega)}.
\end{align*}
{\color{blue}Then,} by Assumption \ref{assump:boundaf} and using the definition of Rademacher complexity{\color{blue}, we can} obtain
\begin{align*}
{\rm I}\leq &\frac{2|\Omega|}{N}\sum_{m=0}^{1}\mathbb{E}_{\mathbb{X},\mathbb{T}, \boldsymbol{\sigma}}\sup_{u \in \mathcal{N}}\left|\sum_{i=1}^{N}\sigma_{i}\nabla^{m}(u-u^{k})(X_{i})\right|\ell_{a}\lambda_{k}^{p-1}\sqrt{d}\|u-u_{k}\|_{W^{1,\infty}(\Omega)}^{p-1}\\
\leq &2|\Omega|\ell_{a}\lambda_{k}^{p-1}\|u-u_{k}\|_{W^{1,\infty}(\Omega)}^{p-1}\sqrt{d}\sum_{m=0}^{1}\mathfrak{R}_N(\nabla^{m}\mathcal{N}).
\end{align*}
Using the $L^p(\Omega)$ norm of the function $u$ and its gradient, we can similarly bound the term ${\rm II}$ as
\begin{align*}
{\rm II}
\leq&\mathbb{E}_{\mathbb{X}}\sup_{u \in \mathcal{N}}\|u\|_{W^{1,\infty}(\Omega)}^{2-p}\left|\int_{\Omega}I_2^k(x,u(x))\,{\rm d}x-\frac{|\Omega|}{N}\sum_{i=1}^{N}I_2^k(X_i,u(X_i))\right|\\
\leq&2|\Omega|\sup_{u \in \mathcal{N}}\|u\|_{W^{1,\infty}(\Omega)}^{2}\sum_{m=0}^{1}\mathfrak{R}_N(\nabla^{m}\mathcal{N}).
\end{align*}
Last, we can bound the term $\rm III$ directly by
\begin{align*}
{\rm III}
\leq&\frac{2\lambda_{k}|\Omega|}{N}\mathbb{E}_{\mathbb{X},\boldsymbol{\sigma}}\sup_{u \in \mathcal{N}}\left|\sum_{i=1}^{N}\sigma_{i}\sum_{m=0}^{1}\nabla^{m}(u-u^{k})(X_{i})\cdot\left[a_{m}(X_{i},u_{k}(X_i),\nabla u_{k}(X_i))-f_{m}(X_{i})\right]\right|\\
\leq &2|\Omega|\left(\ell_{a}\sqrt{d}\|u_{k}\|_{W^{1,\infty}(\Omega)}+\|f_{0}\|_{L^{\infty}(\Omega)}+\|f_{1}\|_{L^{\infty}(\Omega)}\right)\lambda_{k}\sum_{m=0}^{1}\mathfrak{R}_N(\nabla^{m}\mathcal{N}).
\end{align*}
Note that the $W^{1,\infty}(\Omega)$ norm of $u_k\in\mathcal{N}(L,W,B)$ is bounded by
$\|u_k\|_{W^{1,\infty}(\Omega)}\leq WB+\sqrt{d}W^{L-1}B^L$
(see \cite[Lemma 5.9]{Jiao2024} and \cite[Lemma 3.4]{JinLiLu:2022}) and that the Rademacher complexity $\mathfrak{R}_N(\nabla^{m}\mathcal{N})$, $m=0,1$, is bounded by
\begin{align*}
\mathfrak{R}_N(\mathcal{N})&\leq C \frac{W^2LB}{\sqrt{N}} \sqrt{\log \left(LWBN\right)} \quad\mbox{and}
\quad\mathfrak{R}_N(\nabla\mathcal{N})\leq C \frac{\sqrt{d}W^{L}LB^{L}}{\sqrt{N}} \sqrt{\log \left(dLWBN\right)}.
\end{align*}
These estimates follow from uniform bounds on the function class $\mathcal{N}$ and $\nabla\mathcal{N}$, Lipschitz continuity of DNN functions with respect to $\theta$ and Dudley lemma; see, e.g., \cite[Theorem 5.1]{Jiao2024} and \cite[Proposition 3.2]{JinLiLu:2022} for related computation.
Combining the preceding estimates leads to the desired assertion.
\end{proof}

\section{Conclusion}
In this work we have developed a novel neural network based solver for a general class of elliptic problems, termed as iterative deep Ritz method (IDRM), which is applicable to a broad range of monotone elliptic problems  under weak regularity assumptions. We have also provided a convergence analysis of the algorithm. We present several numerical experiments to validate the convergence of the approach and to illustrate its practical applicability. The numerical results show that it can achieve good accuracy, and is well suited for problems requiring high smoothness or lacking a Ritz variational formulation. Numerically it improves existing neural techniques, e.g., physics informed neural networks and deep Ritz method. In sum, IDRM represents a new flexible approach with solid theoretical underpinnings. Potential future works include enhancing its computational efficiency, reducing the learning error and analyzing the optimization error. It is also of much interest to extend the proposed IDRM to elliptic PDEs with discontinuous coefficients (e.g., interface problems), for which the solution may exhibit jumps across the interfaces. This extension would require suitably modifying the neural network architecture according to the discontinuity interface.

\bibliographystyle{abbrv}
\bibliography{ref}

\end{document}